\input amstex

\loadeufm
\loadmsbm
\loadeufm

\documentstyle{amsppt}
\input amstex
\catcode `\@=11
\def\logo@{}
\catcode `\@=12
\magnification \magstep1
\NoRunningHeads
\NoBlackBoxes
\TagsOnLeft

\def \={\ = \ }
\def \+{\ +\ }
\def \-{\ - \ }

\def \b|{\big |}

\def \g1{\Gamma_1}

\def \nfp{\demo\nofrills{Proof:\usualspace\usualspace }}

\def\rarr#1#2{\smash{\mathop{\hbox to .5in{\rightarrowfill}}
 	 \limits^{\scriptstyle#1}_{\scriptstyle#2}}}

\def\larr#1#2{\smash{\mathop{\hbox to .5in{\leftarrowfill}}
	  \limits^{\scriptstyle#1}_{\scriptstyle#2}}}

\def\swarr#1#2 {\llap{$\scriptstyle #1$}  \swarrow
  	\vcenter to .5in{}\rlap{$\scriptstyle #2$}}

\topmatter
\title Diophantine Geometry over Groups VIII: 
\centerline{Stability}
\endtitle
\author
\centerline{ 
Z. Sela${}^{1,2}$}
\endauthor
\footnote""{${}^1$Hebrew University, Jerusalem 91904, Israel.}
\footnote""{${}^2$Partially supported by an Israel academy of sciences fellowship.}
\abstract\nofrills{}
This paper is the eighth in a sequence 
on the structure of
sets of solutions to systems of equations in  free and hyperbolic groups, projections
of such sets (Diophantine sets), and the structure of definable sets over  free and hyperbolic groups.
In the eighth paper we use a modification of the sieve procedure, that was presented in [Se6]
as part of the  quantifier elimination procedure, 
 to prove that free and torsion-free (Gromov) hyperbolic groups are stable.
\endabstract
\endtopmatter

\document

\baselineskip 12pt

In the first 6 papers in the sequence on Diophantine geometry over groups we
studied sets of solutions to systems of equations in a free group, and
developed basic techniques and objects that are required for the analysis of sentences and
elementary sets that are defined over a free group.  The techniques we developed, 
enabled us to present an iterative procedure that analyzes $EAE$ sets defined
over a free group (i.e., sets defined using 3 quantifiers), and shows that
every such set is in the Boolean algebra generated by $AE$ sets 
([Se6],41), hence, we obtained a quantifier elimination over a free group.

In 1983 B. Poizat [Po1] proved that free groups are not super-stable (W. Hodges pointed out to us that
this was also known to Gibone around 1976).
In this paper we use our analysis of definable sets, and the geometric structure they admit as a
consequence from
our quantifier elimination procedure, together
with the tools and the techniques that are presented in the previous  papers in
the sequence, to prove that free groups are stable (Theorem 5.1  - for a definition of a stable theory see [Pi] or
the beginning of section 5).  
Since in [Se8] it was shown that the structure
of definable sets and the tools that were developed for the analysis of them generalize to non-elementary,
torsion-free hyperbolic groups, the argument that we use for proving the stability of a free group
generalizes to an arbitrary non-elementary, torsion-free hyperbolic group (Theorem 5.2).

The stability of free and hyperbolic groups gives a linkage between negative curvature in Riemannian and coarse
geometry and in geometric group theory and stability theory. With stability it is possible
to continue the study of
the first order theories of free and hyperbolic groups using well-developed objects and notions from
model theory. Furthermore, following Shelah, 
logicians often view stability as the border line between "controlled" and
"wild" structures. From certain points of view, and in certain aspects,  this border line is reflected 
in group theory (see  [Po2],[Po3]). 
Negatively curved groups are stable. For non-positively curved groups we don't really
know, but we suspect that there should be unstable non-positively curved groups. For other classes of
groups the question of stability is still wide open.

To prove the stability of free and hyperbolic groups, we start by analyzing a special class of definable
sets that we call $minimal$ $rank$. These sets are easier to analyze than general definable sets, and
in section 1 we prove that minimal rank definable sets are in the Boolean algebra generated by equational sets
(recall that equational sets and theories were defined by G. Srour. For a definition see 
the beginning of section 1 and
[Pi-Sr]). 

In section 2 we slightly modify the sieve procedure that was presented in [Se6] (and used
for quantifier elimination) to prove that Diophantine sets are equational. The equationality
of  Diophantine sets is essentially equivalent to the termination of the sieve procedure for quantifier
elimination in [Se6], and it is  a key in obtaining 
stability for general definable sets in the sequel. In section 3 we present a basic object that we
use repeatedly in proving stability - $Duo$ $limit$ $groups$ (definition 3.1), and their $rectangles$
(definition 3.2). We further prove a boundedness property of duo limit groups and their rectangles
(Theorems 3.3), 
that is not required in the sequel, but still motivates our approach to stability. 

In section 4 we use duo limit groups and their rectangles, together with the sieve procedure and the
equationality of Diophantine sets, to prove the stability of some families of definable sets, that are in a sense
the building blocks of general definable sets (over a free group).
These include the set of values of the defining parameters of a rigid and solid 
limit groups,
for which the rigid (solid) limit group
has precisely $s$ rigid (strictly solid families of) specializations for some fixed integer $s$
 (see section 10 in [Se1] and section 1 in [Se3] for these notions). 

In section 5 we use the geometric structure of a general definable set that was proved using the
sieve procedure in [Se6], together with the stability of the families of definable sets that are 
considered in section 4, to prove the stability of a general definable set over a free group, hence, to
obtain the stability of a free group (Theorem 5.1). Using the results of [Se8] we further generalize our results to
a non-elementary, torsion-free (Gromov) hyperbolic group (Theorem 5.2).


The objects, techniques and arguments that we use in proving stability, are all based on the work on
Tarski's problems, and in particular on the sieve procedure for quantifier elimination ([Se1]-[Se6]). Parts of
the arguments require not only familiarity with the main objects that are presented in these papers, but
also with the procedures that are used in them. We give the exact references wherever we apply these procedures,
or use previously defined notions. 

Quite a few people have assisted us along the course of this work. In particular we would like to thank
G. Cherlin, W. Hodges, O. Belegradek, A. Pillay, B. Zilber, and especially E. Hrushovski for their help and
suggestions. Dave Gabai has encouraged us to revise this paper, and Eliyahu Rips read it thoroughly and made
us double its length. I am grateful to both of them.

\vglue 1.5pc
\centerline{\bf{\S1. The Minimal (Graded) Rank Case}}
\medskip
Our aim in this paper is to prove that free and hyperbolic groups are stable.
Before treating the stability of these groups, we study
a subcollection of definable sets, that we called
$minimal$ $rank$ (in section 1 of [Se5]), 
and prove that these sets are in the Boolean algebra generated by equational
sets (and hence are in particular stable).
 
Recall that a Diophantine set over a free group, $F_k=<a>$, is a projection of a variety, i.e., it is defined as:
$$D(p) \ = \ \{ \, p \ | \ \exists x \ \Sigma (x,p,a)=1 \, \}.$$
With the (set of solutions to the) system of equations, $\Sigma(x,p,a)=1$, one can 
associate canonically finitely many
limit groups (see theorem 7.2 in [Se1]), $L_1(x,p,a),\ldots,L_t(x,p,a)$. If we denote the parameter (free
variables)
subgroup $P=<p>$, then the Diophantine set, $D(p)$, is determined by the finitely many homomorphisms,
$h_i:P \to L_i$, $i=1,\ldots,t$.

\vglue 1pc
\proclaim{Definition 1.1} 
A Diophantine set, $D(p)$,  is called
$minimal$ $rank$ if the targets $L_i$, in the homomorphisms: $h_i:P \to L_i$, $i=1,\ldots,t$, 
that determine the Diophantine set, $D(p)$,  
admit no restricted
epimorphism onto a free product of the coefficient group and an infinite  cyclic group,
$F_k*<t>=<a>*<t>$.
A definable set is called minimal rank, if it is contained in the union of finitely many
minimal rank Diophantine sets.

A parametric family of Diophantine sets is defined as: 
$$D(p,q) \ = \ \{ \, (p,q) \ | \ \exists x \ \Sigma (x,p,q,a)=1 \, \}$$
(where the variables $q$ are considered to be the parameters of the family, and for each value of the variables 
$q$ the fiber is a Diophantine set). 

The parametric family, $D(p,q)$, is called $minimal$ $rank$ if the targets $L_i(x,p,q,a)$, in the homomorphisms:
$u_i:<p,q> \to L_i(x,p,q,a)$, that determine the family, $D(p,q)$, admit no restricted epimorphism
onto a free product of the coefficient group and an infinite cyclic group,
$F_k*<t>=<a>*<t>$, that maps the subgroup $<q>$ into the coefficient group $F_k=<a>$.
A parametric family of definable sets is called minimal rank, if it is contained in the union of finitely many
minimal rank parametric families of Diophantine sets.
\endproclaim

Minimal (graded) rank sets were
 treated separately in our procedure for quantifier elimination ([Se5]-[Se6]), and
it was indicated there that our procedure for quantifier elimination for minimal (graded) rank
formulas  is far easier
than it is for general formulas (see section 1 in [Se5] for the analysis of minimal rank sets).
 
In order to prove that minimal rank families of definable sets are contained in a Boolean algebra 
of equational sets, 
we introduce a collection of (minimal rank) $equational$ $sets$ for which: 
\roster
\item"{(i)}"  the Boolean algebra generated by the  collection of equational sets contains
the collection of minimal rank families of definable sets.

\item"{(ii)}" if $\varphi(p,q)$ is (the formula that defines) an equational set, then 
there exists a constant $N_{\varphi}$,
so that  for every sequence of values $\{q_i\}_{i=1}^{m}$, for which the sequence of sets
that corresponds to the intersections:
$ \{ \, \wedge_{i=1}^{j} \varphi(p,q_i) \, \}_{j=1}^m$ is a strictly decreasing sequence, satisfies:
$m \leq N_{\varphi}$ (O. Belegradek has pointed out to us that this is the definition of 
equationality that one needs to use in case the underlying model is not necessarily saturated).
\endroster

To define
the subcollection of equational sets, and  prove the descending chain
condition that they satisfy, we study the Boolean algebra of minimal rank definable sets gradually. 
\roster
\item"{(1)}"  Diophantine sets - we show that minimal rank parametric families of Diophantine sets are equational.

\item"{(2)}" Rigid limit groups are defined in section 10 of [Se1], and their rigid values
are analyzed in sections 1-2 of [Se3].  In theorem 2.5 in [Se3] it is proved  that given a rigid
limit group, there exists a global bound on the number of rigid values that are associated
with any possible value of
the defining parameters. 

With a given minimal rank rigid limit group $Rgd(x,p,q,a)$ (where $<p,q>$ is the
parameters group),  we associate a natural existential formula,
$\varphi(p,q)$, that specifies those values of the defining parameters $p,q$ for
which $Rgd(x,p,q,a)$ admits at least $m$ rigid values, for some fixed 
integer $m$. We prove the existence of a collection of equational formulas, so that
the Boolean algebra generated by this collection contains all the formulas
$\varphi(p,q)$, that are  associated with all minimal rank rigid limit groups and an arbitrary integer $m$.

\item"{(3)}" Solid limit groups are defined in section 10 of [Se1], and their strictly solid
families of  specializations 
are analyzed in sections 1-2 of [Se3] (see definition 1.5 in [Se3]). In parallel with rigid limit groups, it
is proved in theorem 2.9 in [Se3] that given a solid limit group, there exists a global bound on the number of
strictly solid families that are associated with any possible value of the defining parameters (strictly 
solid families of a solid limit group are defined in definition 1.5 in [Se3]). 

With a given minimal rank solid limit group $Sld(x,p,q,a)$ 
  we associate a natural $EA$ formula, 
$\varphi(p,q)$, that specifies those values of the defining parameters $p,q$ for
which $Sld(x,p,q,a)$ admits at least $m$ strictly solid families of values, 
for some fixed 
integer $m$. As for rigid limit groups, we show the existence of a 
collection of equational formulas, so that
the Boolean algebra generated by this collection contains all the formulas
$\varphi(p,q)$, that are associated with minimal rank solid limit groups and an arbitrary integer $m$.

\item"{(4)}" Given a graded resolution (that terminates in either a rigid or a solid limit group),
and a finite collection of (graded) closures of that graded resolution,
we define a natural
formula, $\alpha(p,q)$ (which is in the Boolean algebra of $AE$ formulas), that specifies those values of  
the defining parameters for which the given set of closures forms a covering closure of the given graded
resolution (see definitions 1.15 and 1.16 in [Se2] for a closure and a covering closure).
We show the existence of a 
collection of equational sets, so that
the Boolean algebra generated by this collection contains all the sets that are defined by the formulas
$\alpha(p,q)$ that are associated with all the graded resolutions for which their terminal
rigid or solid limit group is of minimal (graded) rank.

\item"{(5)}"  Finally, we show the existence of a subcollection of 
equational sets that generates the Boolean algebra of minimal rank parametric families of definable sets.
\endroster

\vglue 1pc
\proclaim{Theorem 1.2} Let $F_k=<a_1,\ldots,a_k>$ be a non-abelian free group, and let:
$$D(p,q) \ = \ \{ \, (p,q) \ | \ \exists x \ \Sigma (x,p,q,a)=1 \, \}$$
be a minimal rank parametric family of Diophantine sets that is defined over $F_k$ 
(where the variables $q$ are considered to be the parameters of
the family). 

\noindent
Then $D(p,q)$ is equational.
\endproclaim

\nfp 
We need to show that $D(p,q)$  is equational, i.e.,
that there exists an integer $N_{D}$,
so that  every sequence of values, $\{q_i\}_{i=1}^{m}$, for which the sequence of intersections:
$ \{ \, \cap_{i=1}^{j} D(p,q_i) \, \}_{j=1}^m$ is a strictly decreasing sequence, satisfies:
$m \leq N_{D}$.

Let: $L_1(x,p,q,a),\ldots,L_t(x,p,q,a)$ be the finite  collection
of maximal limit groups that is canonically associated with the system of equations $\Sigma(x,p,q,a)=1$
(see theorem 7.2 in [Se1] for the existence of this canonical finite collection). Since we assume that $D(p,q)$ is
a minimal rank family of  Diophantine sets,
each of the limit groups $L_i(x,p,q,a)$ is of minimal rank, when viewed
as a graded limit group with respect to the parameter subgroup $<q>$ (i.e., $L_i(x,p,q,a)$ admits no
restricted epimorphism onto a free group $F_k*F$ where $F$ is non-trivial free group, and
the subgroup $<q>$ is mapped into the coefficient
group $F_k$). 

To prove the existence of a bound $N_{D}$, we associate with the set $D(p,q)$
a universal finite diagram. The construction of the diagram is based on the sieve procedure for quantifier
elimination in the minimal rank case, that is presented in section 1 of [Se5]. Once the universal diagram
is constructed, equationality of the original family of Diophantine sets, $D(p,q)$, will be deduced, by uniformly
bounding the lengths of certain (decreasing) paths along the constructed diagram. In particular, the equationality
constant, $N_{D}$, can be computed from the diagram.

We start the construction of the universal finite diagram with each of the maximal limit groups, 
$L_1(x,p,q,a),\ldots,
L_t(x,p,q,a)$, in parallel. With a limit group, $L_i(x,p,q,a)$, viewed as a graded limit
group with respect to the parameter subgroup $<q>$, we associate its 
strict graded Makanin-Razborov
diagram (for the construction of the strict Makanin-Razborov diagram, see proposition 1.10
in [Se2]. The modification of a graded Makanin-Razborov diagram to a strict diagram is identical to the ungraded case,
and the strict graded Makanin-Razborov diagram is used repeatedly in the 
quantifier elimination procedure, e.g., in the proof of theorem 1.4 in [Se5]). 
With each resolution in the graded strict Makanin-Razborov diagram, we further associate
its singular locus (the singular locus of a graded resolution collects all the rigid or strictly solid
values of the rigid or solid terminal limit group of the graded resolution, for which the fiber
of specializations that is associated with such value is degenerate - see section 11 in [Se1]
for the exact definition, stratification,  and the construction of the singular locus), 
and the strict graded resolutions that are associated with
each of the strata in the singular locus. 

Altogether we have a finite collection of strict graded resolutions, those that appear in the strict 
graded Makanin-Razborov diagrams of the groups $L_i$, and those that are associated with the strata
in their singular loci. 
We conclude the first step of the
construction of the diagram, by associating the (graded) completion  with each of the
graded resolutions in our finite collection, that we denote, $Comp(x,p,z,q,a)$ (see definition 1.12  in [Se2] for
the completion of a strict  resolution), and
with each graded completion we associate its complexity, according to 
definition 1.16 in [Se5]. These (finitely many) completions form the first level of the universal diagram. 

We continue to the construction of the second level of the diagram with each of the completions 
$Comp(x,p,z,q,a)$ in parallel. With each such completion we associate the collection of
all the values, $(x^0_1,x^0_2,p_0,z_0,q_1^0,q_2^0,a)$, for which:
\roster
\item"{(1)}" $(x_1^0,p_0,z_0,q_1^0,a)$ is a specialization of the completion, $Comp(x,p,z,q,a)$.

\item"{(2)}" $(x_2^0,p_0,q_2^0,a)$ is a specialization of at least one of  the maximal limit groups, 
$L_i(x,p,q,a)$,
that is (canonically) associated with the system of equations $\Sigma(x,p,q,a)$ (that defines the 
Diophantine set $D(p,q)$).
\endroster
By the standard arguments that are presented in section 5 of [Se1], with this collection of
values we can
canonically associate a canonical finite collection of maximal limit groups, $M_j(x_1,x_2,p,z,q_1,q_2,a)$, which we
view as graded limit groups with respect to the parameter subgroup $<q_1,q_2>$ (note that the finite collection
of limit groups, $\{M_j\}$, is dual to the Zariski closure of the given collection of values).

Since we assume that each of the limit groups, $L_i(x,p,q,a)$, is of minimal rank,
each of the completions, $Comp(x,p,z,q,a)$, is of minimal rank as well (as the limit groups
that are associated with the various levels of these completions are quotients of the limit groups, $L_i$).

\noindent
The limit groups, $M_j(x_1,x_2,p,z,q_1,q_2,a)$, are constructed from specializations of the completions,
$Comp(x,p,z,q,a)$, and the limit groups, $L_i(x,p,q,a)$.
Since both the completions, $Comp(x,p,z,q,a)$, and the limit groups, $L_i(x,p,q,a)$, are of minimal rank,
so are the limit groups, $M_j$, i.e. each of the limit groups, $M_j$, admits no epimorphism onto a free
group, $F_k*F$, where $F$ is a non-trivial free group, and the subgroup, $<q_1,q_2>$, 
is mapped into the coefficient group $F_k$. To analyze the values of the defining parameters, $(p,q)$, that extend to 
values of the constructed limit groups, $M_j$, we need the following theorem.

\vglue 1pc
\proclaim{Theorem 1.3} Let $L(x,p,q,a)$ be a minimal rank graded limit group (graded with respect to the parameter subgroup
$<q>$),  let $LRes(x,p,q,a)$ be a graded resolution of
$L(x,p,q,a)$, and let $Comp(LRes)(x,p,z,q,a)$ be the completion of the resolution, $LRes$. 
With the resolution $LRes$ we can associate
a complexity according to definition 1.16 in [Se5].
 
Let $U(u,p,q,a)$ be a minimal rank limit 
group, and let $S(u,x,p,z,q_1,q_2,a)$ be a limit group that is obtained as a limit of a sequence of values:
$\{(u(n),x(n),p(n),z(n),q_1(n),q_2(n),a)\}_{n=1}^{\infty}$, where the tuples, $(x(n),p(n),z(n),q_1(n),a)$, 
are specializations of the completion,
$Comp(LRes)(x,p,z,q,a)$, and the tuples, $(u(n),p(n),q_2(n),a)$, are specializations of the limit group, $U(u,p,q,a)$.

With the graded limit group, $S(u,x,p,z,q_1,q_2,a)$, it is possible to associate finitely many strict
graded resolutions:
$VRes_1(v,x,p,z,q_1,q_2,a),\ldots,VRes_{\ell}(v,x,p,z,q_1,q_2,a)$, that are graded with respect to the parameter subgroup,
$<q_1,q_2>$, for which:
\roster
\item"{(1)}" every value, $(p_0,q_1^0,q_2^0)$, of the variables $p,q_1,q_2$, 
that can be extended to a specialization of the limit group, 
$S(u,x,p,z,q_1,q_2,a)$, can be extended to a value, $(v_0,x_0,p_0,z_0,q_1^0,q_2^0,a)$, 
that factors through at least one of the 
resolutions, $VRes_i$,
$i=1,\ldots,\ell$.

\item"{(2)}" the complexity of each of the resolutions, $VRes_i$, is bounded by the complexity of the resolution that
we have started
with, $LRes$ (where the complexity of a minimal rank resolution is the one presented in definition 1.16 in [Se5]).

\item"{(3)}" if the complexity of a resolution, $VRes_i$, is equal to the complexity of the graded resolution, $LRes$,
then the completion of $VRes_i$, $Comp(VRes_i)$, has the same structure as a graded closure of the completion of $LRes$, $Comp(GRes)$ (see definition 
1.14 in [Se2] for a closure of a completion). That is $Comp(VRes_i)$ is obtained from $Comp(LRes)$ by possibly adding roots to abelian
vertex groups in abelian decompositions that are associated with the various levels of $Comp(LRes)$, and replacing the terminal 
rigid or solid limit group of $Comp(LRes)$ (which is graded with respect to the parameter subgroup $<q>$), with a rigid or solid
limit group with respect to the parameter subgroup $<q_1,q_2>$.
\endroster
\endproclaim

\nfp The construction of such a finite set  of resolutions,
is precisely the  construction that is conducted in the general step of the sieve procedure for quantifier
elimination in the  minimal rank case, in section 1 of [Se5] (see the proof of  theorem 1.22 in [Se5]).

\line{\hss$\qed$}

Both the completions, $Comp(x,p,z,q,a)$, and the limit groups, $M_j(x_1,x_2,p,z,q_1,q_2,a)$, that we have associated
with the Diophantine set, $D(p,q)$, are
of minimal rank, and the limit groups, $M_j$, are obtained from a collection of specializations of
a minimal rank completion, $Comp(x,p,z,q,a)$, by imposing on them an additional (minimal
rank) Diophantine
conditions. Hence, the assumptions of theorem 1.3 are satisfied, with $Comp(x,p,z,q,a)$, $L_i(x,p,q,a)$, and $M_j$  in place of $Comp(LRes)$,
$U(u,p,q,a)$, and $S$, in the statement of theorem 1.3.  
By the conclusion of the theorem, with each of the limit groups, $M_j$, we can associate finitely many 
strict graded resolutions that satisfy properties (1)-(3) in the statement of theorem 1.3.


Therefore, some of the strict graded resolutions that are associated with  a limit group, $M_j$,
have the structure of (graded) closures of the completion, $Comp(x,p,z,q,,a)$, from which $M_j$ was constructed
(part (3) in theorem 1.3),
and the other resolutions have strictly
smaller complexity than the complexity of $Comp(x,p,z,q,a)$. Those constructed resolutions that have
strictly smaller complexity than their associated completion, $Comp(x,p,z,q,a)$, or the structure of proper closures
of the completion, $Comp(x,p,z,q,a)$, i.e., the structure of closures  that contain non-trivial roots of elements in abelian 
vertex groups that are associated with with (abelian decompositions of)
$Comp(x,p,a,q,a)$,  form the second level of the universal diagram that is associated with
the Diophantine set, $D(p,q)$. We add a directed edge from each of the completions, $Comp(x,p,z,q,a)$, that
form the first level of the diagram, to each of the
graded resolutions  that are associated with it in the second level of the diagram.

We continue to the third step of the construction of the diagram with each of the graded resolutions in the second 
level, and we continue in parallel.
Given such a (graded) resolution, we repeat  the same operations that
we have conducted in the second step. Given a  graded resolution in the second level of the diagram, 
we take its completion, 
and look at 
all the specializations of that completion,  for which there exists a
value $(x_3^0,q_3^0)$, so that the combined value, $(x_3^0,p^0,q_3^0,a)$, is a specialization
of one of the 
limit groups, $L_i(x,p,q,a)$, that are associated with the system of equations $\Sigma(x,p,q,a)$
(that was used to define the Diophantine set $D(p,q)$). 
With the collection of these values, we canonically associate a finite collection of
maximal limit groups (that is associated with the Zariski closure of the given collection
of values according to section 5 of [Se1]). Each such maximal limit group has to be of minimal rank
(it does not admit an epimorphism onto a free group $F_k*F$ (where $F$ is non-trivial) 
that maps the parameter subgroup $<q_1,q_2,q_3>$ 
into the coefficient group $F_k$). 

By theorem 1.3 with the obtained  (minimal rank) maximal
limit groups and the completions of the resolutions in the second level from which they were constructed,
we associate 
 a finite collection of minimal rank graded resolutions (with respect to the parameter subgroup
$<q_1,q_2,q_3>$).
By part (2) of theorem 1.3, the complexity of each of the constructed
resolutions is bounded by the complexity of the resolution in the second level of the diagram from which it was constructed.
Furthermore, by part (3) of theorem 1.3, in case of equality in the complexities
of a graded resolution that appears in the second level of the diagram, and a constructed graded resolution that
was constructed from it, the constructed resolution has to have the structure of a closure of the graded resolution from the second level
(the structure of a closure in the sense of part (3) in theorem 1.3).

Those constructed resolutions that have
strictly smaller complexity than the resolution in the second level from which they were constructed,
or those constructed resolutions that have the structure of (proper) closures
of the associated resolutions in the second level, in which non-trivial roots were added to abelian vertex groups that
are associated with the resolution from the second level,
form the third level of the universal diagram that is associated with
the Diophantine set, $D(p,q)$. We add a directed edge from each of the resolutions in the second level to any
of the (finitely many) graded resolutions that were constructed from  it in the third level of the diagram.

We continue the construction iteratively and repeat the same operations at each step. Given a graded resolution that appears
in level $n$ of the diagram, we associate with it its completion. Then we collect all the specializations that
factor through this completion, and satisfy an additional Diophantine condition, i.e., their restrictions to the variables $p$,
extend to values that factor
through one of the finitely maximal limit groups, $L_i$,  that are associated with the system of equations $\Sigma$, that
was used to define the Diophantine set, $D(p,q)$, and these values of the limit groups $L_i$ restrict to values $q_{n+1}^0$
of the variables $q$ in the generating set of the limit groups $L_i$. 
We associate with  the collection of the  combined values (the value of the completion of the graded resolution in level $n$,
and the corresponding value of some limit group $L_i$) its 
canonical collection of maximal limit groups. 
By part (2) of theorem 1.3, the complexity of each of the constructed
resolutions is bounded by the complexity of the resolution in the second level of the diagram from which it was constructed.
By part (3) of theorem 1.3, in case of equality in the complexities
of a graded resolution that appears in the $n$-th level of the diagram, and a constructed graded resolution that
was constructed from it, the constructed resolution has to have the structure of a closure of the graded resolution from the $n$-th level
(the structure of a closure in the sense of part (3) in theorem 1.3).

We continue to level $n+1$ of the diagram only with those constructed resolutions that have
strictly smaller complexity than the resolution in the $n$-th level from which they were constructed,
or those constructed resolutions that have the structure of proper closures
of the associated resolutions in $n$-th  level. 
We add a directed edge from each of the resolutions in the $n$-th level to any
of the (finitely many) graded resolutions that were constructed from  it in the $n+1$-th level of the diagram.

The diagram that we constructed is locally finite, hence, we may apply Konig's lemma to
prove that its construction terminates.
The complexities of graded resolutions along a path in the diagram are non-increasing. By
theorem 1.18 in [Se5],  a strict reduction in the complexities of  successive resolutions along a path in the diagram
can occur only at finitely
many levels. Given a resolution in the diagram, and a subpath (that have the structures) of proper closures of it
(in the sense of part (3) in theorem 1.3),
its successive resolutions along the path are obtained from it by imposing one of (fixed) finitely
many Diophantine conditions.

A graded resolution that has the structure of a proper closure of its preceding one along a path in the diagram,
is obtained from the completion of its preceding one by adding proper roots to some of the  abelian vertex groups that are associated
with the preceding completion. 
By theorem 1.3, given a completion along a path, there are finitely many 
graded resolutions that are associated with it in the next level of the diagram, and in particular finitely many graded
resolutions that have the structure of
proper closures of the original completion.

Therefore, given a completion along a path in the constructed diagram, there is a global bound (that depends only
on the completion and the finitely many Diophantine conditions) on the index  
of abelian supergroups of abelian vertex groups  that are associated with the given  completion, for all the graded
resolutions that have the structure of proper closures of the
given completion along the given path. Hence, there is a bound on the length of a subpath that starts
with the completion and continues from it with a sequence of resolutions that have the structure of proper closures of it.
The finiteness  of  subpaths of proper closures, together with the finiteness of
the number of levels with a complexity reduction, along a given path in the diagram, 
imply that every path in the diagram has to be finite. Therefore, by Konig's lemma,
the constructed diagram is finite. 

Note that the obtained diagram is  a directed forest, where at each vertex we placed a (strict) 
graded resolution, or alternatively its completion. Furthermore, the constructed diagram is universal, which in particular means that
given an arbitrary sequence of values of the defining parameters: $q_1,q_2,\ldots,q_m$, we can analyze the structure of the intersections:
$ \{ \, \cap_{i=1}^{j} D(p,q_i) \, \}_{j=1}^m$ using the constructed diagram.

By theorems 2.5, 2.9 and 2.13 of [Se3] the number of rigid or strictly
solid families of values of a rigid or solid limit group, 
that are  associated with a given value of the defining parameters,
is uniformly bounded by a bound that depends only on the rigid or solid limit
group (and not on the specific value of the defining parameters). 

Let $depth$ be the number of levels in the universal diagram that we have
associated with the Diophantine set, $D(p,q)$, and let $w$ be the maximal number
of vertices in a single level of the diagram. At each vertex in the diagram
we have placed a graded resolution (or alternatively, a graded completion of that
resolution). Each such graded resolution terminates in either a rigid or a
solid limit group, and by theorems 2.3, 2.9, and 2.13 in [Se3],  with each such
terminal rigid or solid limit group, there is a corresponding global bound on
the number of rigid or strictly solid families of values that is
associated with any possible value of its defining parameters. Let $b$ be the maximum
of all the global bounds that are associated with the terminal rigid or solid limit
groups of all the graded resolutions that are associated with the vertices in the
constructed universal diagram.

In the sequel we will often need the following notion:

\vglue 1pc
\proclaim{Definition 1.4} Let $GRes(x,p,q,a)$ be a graded resolution that terminates in either a
rigid limit group, $Rgd(x,p,q,a)$, or a solid limit group, $Sld(x,p,q,a)$. Let $Comp(x,p,z,q,a)$ be the
completion of $GRes$. A $fiber$ of the graded resolution, $GRes$, or of its completion, $Comp$, is the set
of specializations of the completion, $Comp$, that extends a given rigid or strictly solid value of 
the terminal rigid or solid limit group, $Rgd$ or $Sld$, of the resolution $GRes$.

A $q$-$fiber$ of the graded resolution, $GRes$, or of its completion, $Comp$, is the bounded collection of
fibers that extends a given value $q_0$ of the defining parameters $q$.
\endproclaim

Let $q_1^0,q_2^0,\ldots$, be a given sequence of values (in the
coefficient group $F_k$) of the (free) variables $q$ in the Diophantine set, $D(p,q)$.
First, we look at $q_1$ as parameters. There are at most
$w$ graded resolutions in the first level of the constructed universal diagram,
and  there are at most $b$ fibers that are associated with each of these graded
resolutions and with the specialization $q_1^0$. Hence, there are at most $w b$ fibers
of the graded resolutions in the first level of the diagram that are associated with
$q_1^0$.

If $D(p,q_1^0) \, \cap \, D(p,q_2^0) \ = \ D(p,q_1^0)$, there is no change. If 
 $D(p,q_1^0) \, \cap \, D(p,q_2^0) \ \neq  \ D(p,q_1^0)$, the intersection of the
two Diophantine sets is strictly contained in $D(p,q_1^0)$. The Diophantine set $D(p,q_1^0)$ is a 
finite union of at most $wb$ fibers of the graded resolutions in the initial level of the diagram. Since
the intersection of the two Diophantine sets, $D(p,q_1^0)$ and $D(p,q_2^0)$,  
is strictly contained in $D(p,q_1^0)$, it is a finite union of fibers - 
a proper (possibly empty) subset of
the fibers that are associated with $D(p,q_1^0)$, and at least one of the fibers that  is associated with
$D(p,q_1^0)$, that is replaced by a
(possibly empty)  finite collection of fibers  that are associated with the pair,
$(q_1^0,q_2^0)$, and with some of the graded resolutions that appear in the second level
of the diagram. By the structure of the universal diagram, each fiber in the first
level can be replaced by at most $w  b$ fibers in the second level.

We repeat this argument iteratively. Each time a value $q_n^0$ is added, and
the corresponding intersection is a proper subset of the previous intersection,
at least one of the fibers that was associated with the intersection of the first $n-1$ values,
 is replaced by at most $w  b$ fibers in level that succeeds the level of that fiber (a fiber that is associated
with the last level of the diagram can only be  replaced by the empty set). As the digram has
$depth$ levels, it takes at most 
  $ depth \cdot {(wb)}^{depth-1}$ values of the variables $q_n$ (for which there is a strict reduction in the corresponding 
intersection) to be left with at most
  $ depth \cdot {(wb)}^{depth}$ 
fibers in the terminal level of the diagram, and at most an additional 
  $ depth \cdot {(wb)}^{depth}$ values of the variables $q_n$ (for which there is a strict reduction in the corresponding intersection)
to eliminate these fibers in the terminal level. Therefore, altogether there can be at most 
  $2 \cdot  depth \cdot {(wb)}^{depth}$ values of the variables $q_n$ for which there is a strict reduction in the  intersection:
$  \cap_{i=1}^{j} D(p,q_i) $, 
  which proves the equationality of the set $D(p,q)$ 
(where the equationality constant satisfies: $N_{D}= 
  2 \cdot depth \cdot {(wb)}^{depth}$).

\line{\hss$\qed$}

\vglue 1pc
\proclaim{Theorem 1.5} Let $F_k=<a_1,\ldots,a_k>$ be a non-abelian free group, and let
$Rgd(x,p,q,a)$ be a rigid limit group, 
 with respect to the parameter subgroup $<p,q>$.
Let $s$ be a positive integer, and let $NR_s$ be the set of values of the
defining parameters $<p,q>$ for which the rigid limit group, $Rgd(x,p,q,a)$, has
at least $s$ rigid values. 

There exists a collection of equational
sets, so that the Boolean algebra generated by this collection contains the sets $NR_s$
for every minimal rank rigid limit group $Rgd(x,p,q,a)$,
and every possible integer $s$.  
\endproclaim

\nfp We construct iteratively a collection of equational sets that generate a Boolean algebra that contains the
sets of the form $NR_s$. With a set of the form $NR_s$, we 
associate a minimal rank Diophantine set $D_1$, and show that $NR_s \cup D_1$ is equational. Clearly:
$$NR_s \, = \, ((NR_s \cup D_1) \,  \setminus \, D_1) \, \cup \, (D_1 \cap NR_s).$$
Since 
by theorem
1.2 the minimal rank Diophantine set $D_1$ is equational, to prove the theorem we further need to study the
set $D_1 \cap NR_s$. We study this set in the same way we treated the set $NR_s$. We further associate  a complexity
with the sets $NR_s$ and $NR_s \cap D_1$, and  argue
that the complexity of the set $D_1 \cap NR_s$ is strictly smaller than the complexity of the
original set $NR_s$. We continue iteratively. At each step, we add a (minimal rank) Diophantine correction to the 
remaining set from the previous step,
prove the equationality of the union of the remaining set and the Diophantine correction, and argue that the intersection of the
Diophantine correction and the remaining set from the previous step has strictly lower complexity. Finally, the reduction in
 complexity forces the iterative procedure to terminate,
hence, prove the theorem for the sets $NR_s$. 

We start with the construction of the set $D_1$ that is associated with the set $NR_s$. 
As a preparation to the definition of $D_1$, we look at the collection of all the tuples of values,
$(x_1^0,\ldots,x_s^0,p_0,q_0,a)$, for which for every index $i$, $1 \leq i \leq s$, $(x_i^0,p_0,q_0,a)$ is a rigid
value of the given rigid limit group $Rgd(x,p,q,a)$, and for every couple $i,j$,
$1 \leq i < j \leq s$, $x_i^0 \neq x_j^0$. By our standard arguments, with this collection
of values we can associate canonically a finite collection of maximal limit groups,
$T_j(x_1,\ldots,x_s,p,q,a)$.

We continue with each of the limit groups $T_j(x_1,\ldots,x_s,p,q,a)$ in parallel. With $T_j$
viewed as a
graded limit group with respect to the parameter subgroup $<q>$,
we associate its 
strict graded Makanin-Razborov
diagram (according to the construction of this diagram as it appears in proposition 1.10 in [Se2]).
With each resolution in the strict graded Makanin-Razborov diagram, we further associate
its singular locus, and the strict graded resolutions that are associated with
each  stratum in the singular locus.  With each of the obtained graded resolutions we
further associate its (graded) completion (according to definition 1.12 in [Se2]), and
with each graded completion we associate its complexity, according to 
definition 1.16 in [Se5].

We continue with each of the completions in parallel. Given such a completion,
 we look at all its specializations,  for which
either one of the values that are supposed to be rigid is flexible, or those for which
two rigid values that are supposed to be distinct coincide. Note that the conditions that we
impose on the specializations of the completions are all basic conditions, i.e., the specializations
are required to satisfy one of finitely many possible additional equations. With the collection of all such
specializations we can associate a canonical finite collection of (graded) limit groups. Each
such graded limit group is minimal rank  by our assumptions. Hence, 
  we can associate with it a finite
collection of resolutions according to theorem 1.3.

By theorem 1.3, some of the associated graded resolutions have the structure of  graded closures of
the original resolution, and the rest
have strictly smaller complexity than the completion that they were constructed from. Since the maximal limit
groups that we analyze are obtained from specializations of completions of the original
resolutions that satisfy one of finitely many additional basic conditions, each of the graded
resolutions of these limit groups that have the structure of a  graded closure of the completion from which it was constructed
(see part (3) of theorem 1.3),
has the precise structure of the completion that it was constructed from, i.e., no proper roots were
added to any of the abelian vertex groups that are associated with the completion from which the graded
resolution was constructed.

We omit the subcollection of resolutions that have the structure of graded closures from the list of associated graded resolutions
that we constructed. With each resolution
that has strictly smaller complexity, we associate its completion, and we set the Diophantine set $D_1$
to be the disjunction of all the Diophantine sets that are associated with completions of those resolutions
that are not of maximal complexity, i.e., resolutions that do not have the structure of graded closures.

\vglue 1pc
\proclaim{Remark} With the set $NR_s$ we have associated finitely many graded limit groups, $T_j$. With
these limit groups we have associated the resolutions in their strict graded Makanin-Razborov diagrams. 
By adding the set $D_1$ to the set $NR_s$, we fill all the fibers that are associated with these graded 
resolutions and contain at least one (in fact, generic) point from the set $NR_s$.
\endproclaim

\vglue 1pc
\proclaim{Proposition 1.6} The set $NR_s  \cup  D_1$ is equational.
\endproclaim

\nfp To prove the proposition, we associate with the set $NR_s \cup D_1$ a finite diagram, that
is constructed iteratively, in a similar way to the construction of the diagram that is associated with
a minimal rank Diophantine set that and was used in the proof of theorem 1.2.

We start the construction of the diagram with the collection of all the tuples of values,
$(x_1^0,\ldots,x_s^0,p_0,q_0,a)$, for which for every index $i$, $1 \leq i \leq s$, $(x_i^0,p_0,q_0,a)$ is a rigid
value of the given rigid limit group $Rgd(x,p,q,a)$, and the $x_i^0$'s are distinct,
and the collection of all tuples $(u_0,p_0,q_0,a)$ that are specializations  of one of the (finitely
many) completions that are associated with 
the Diophantine set $D_1$, $Comp(u,p,q,a)$. By our standard arguments, with this collection
of values we  associated canonically a finite collection of maximal limit groups,
$T_j(x_1,\ldots,x_s,p,q,a)$, and the completions that are associated with $D_1$, $Comp(u,p,q,a)$.

We continue with each of the limit groups $T_j(x_1,\ldots,x_s,p,q,a)$ in parallel. 
With a maximal limit group $T_j$,
viewed as a
graded limit group with respect to the parameter subgroup $<q>$,
we associate its 
strict graded Makanin-Razborov
diagram.
With each resolution in the strict graded Makanin-Razborov diagram, we further associate
its singular locus, and the graded strict resolutions that are associated with
each of the strata in the singular locus.  With each of the obtained strict graded resolutions we
further associate its (graded) completion, $Comp(x_1,\ldots,x_s,p,z,q,a)$, and
with each graded completion we associate its complexity, according to 
definition 1.16 in [Se5]. These graded resolutions and their completions, together with the graded resolutions
that are associated with the Diophantine set $D_1$, and their completions, $Comp(u,p,q,a)$, 
form the first level of the diagram
that we associate with $NR_s \cup D_1$.

In the second level of the diagram we need to place resolutions that will  assist us in analyzing
the intersections: $(NR_s \cup D_1)(p,q_1) \, \cap \, (NR_s \cup D_1)(p,q_2)$. These intersections can be
written as unions of sets of the form: $NR_s(p,q_1) \cap NR_s(p,q_2)$, $NR_s(p,q_1) \cap D(p,q_2)$, and
$D(p,q_1) \cap D(p,q_2)$.

We start the construction of the resolutions in the second level of the diagram with each of the completions, 
$Comp(x_1,\ldots,x_s,p,z,q,a)$, and each of the the completions that are associated with $D_1$, $Comp(u,p,q,a)$,
  in parallel. 
With each  completion, $Comp(x_1,\ldots,x_s,p,z,q,a)$, we associate the collection of
values, $(y_1^0,\ldots,y_s^0,x_1^0,\ldots,x_s^0,p_0,z_0,q_1^0,q_2^0,a)$ and 
$(u_2^0,x_1^0,\ldots,x_s^0,p_0,z_0,q_1^0,q_2^0,a)$, and with each completion that is associated with $D_1$,
$Comp(u,p,q,a)$,
we associate
the collection of values,
 $(y_1^0,\ldots,y_s^0,u_1^0,p_0,q_1^0,q_2^0,a)$ and 
$(u_2^0,u_1^0,p_0,q_1^0,q_2^0,a)$,  so that the restrictions of these
values satisfy the 
following conditions:
\roster
\item"{(1)}" $(x_1^0,\ldots,x_s^0,p_0,z_0,q_1^0,a)$ is a
specialization of the completion, 
$Comp(x_1,\ldots,x_s,p,z,q,a)$.

\item"{(2)}" the values,
 $(x_i^0,p_0,q_1^0,a)$, $i=1,\ldots,s$, are distinct rigid values of the rigid limit group, 
$Rgd(x,p,q,a)$, and so are the values, $(y_i^0,p_0,q_2^0,a)$, $i=1,\ldots,s$.

\item"{(3)}"  $(u_j^0,p_0,q_j^0,a)$, $j=1,2$, is a  specialization of one of the completions that is
associated with the set $D_1$, $Comp(u,p,q,a)$.
\endroster

\noindent
With this set of values we 
(canonically) associate a canonical finite collection of maximal limit groups 
(according to theorem 7.2 in [Se1]),  which we
view as graded limit groups with respect to the parameter subgroup $<q_1,q_2>$.

By our assumptions each of the completions, $Comp(x_1,\ldots,x_s,p,z,q,a)$, and each of the completions
that is associated with $D_1$, $Comp(u,p,q,a)$, is of minimal rank. 
Hence, we may apply theorem 1.3,
and associate with each of the (finitely many) graded limit
groups that is associated with the collection of values under consideration, 
a finite collection of minimal rank strict graded resolutions (with respect
to the parameter subgroup $<q_1,q_2>$). By theorem 1.3,
the complexity of each of these minimal 
rank resolutions is bounded above by the complexity of the resolution from which  the 
corresponding completion was constructed (the completion in the first level of the diagram
with which we have started the construction of the corresponding part of the second level),
$Comp(x_1,\ldots,x_s,p,z,q,a)$ or $Comp(u,p,q,a)$. By the same theorem, in case of equality in complexities
(between a constructed resolution and the completion it was constructed from),  the obtained resolution 
has to have the structure of  a graded
closure of the  completion from which it was constructed (see part
(3) of theorem 1.3 for the properties of that structure). 
Therefore, some of the obtained resolutions
have the structure of  (graded) closures of the completions,  
$Comp(x_1,\ldots,x_s,p,z,q,a)$ 
and $Comp(u,p,q,a)$, 
and the other resolutions have strictly
smaller complexity than the complexity of corresponding completion, $Comp(x_1,\ldots,x_s,p,z,q,a)$ 
or $Comp(u,p,q,a)$.

We continue to the third level only with those resolutions that have strictly smaller complexity than
the  completion from which they were constructed, or with resolutions that have the structure of  proper
closures of the completions from which they were constructed (see part (3) in theorem 1.3). 
Given such a (graded) resolution, we perform the same operations
that we have conducted in constructing the second level, i.e., we take its completion, and look at 
all the specializations of that completion, that satisfy the corresponding
(non-degeneration) rigidity conditions, and for which 
either there exists 
a value, $(t_1^0,\ldots,t_s^0,q_3^0)$, so that the combined values, $(t_i^0,p_0,q_3^0,a)$, are
distinct rigid values of $Rgd(x,p,q,a)$, or a value, $(u_0,p_0,q_3^0,a)$ that is a specialization of one
of the (finitely many) completions, $Comp(u,p,q,a)$, 
 that are associated with the Diophantine set $D_1$.

With the collection of these values, we canonically associate a finite collection of
maximal limit groups, that are all of minimal rank, and with them we associate finitely many 
(minimal rank) strict graded resolutions 
by applying theorem 1.3. 
By theorem 1.3, the complexity of each of the associated graded  
resolutions is bounded by the complexity of the corresponding resolution from the second level of the diagram
from which it was constructed.
We continue to the fourth level, only with those resolutions that have strictly smaller complexity
than the resolution from the second level from which they were constructed,
or with graded resolutions that have the structure of proper closures of 
the completion  of that resolution (i.e., graded resolutions that satisfy part (3) in theorem 1.3 and for
which proper roots were added to some of the abelian vertex group that are associated with the completion from
which they were constructed).

We continue the construction iteratively. Since the obtained diagram is locally finite, we may
apply Konig's lemma to prove the finiteness of the diagram. By theorem 1.32 in [Se5], a reduction in the
complexity of successive resolutions can occur only at finitely many steps along a path in the diagram.
By the same argument that we used in proving theorem 1.2, every subpath in the constructed diagram in
which a successive resolution is a proper closure of its predecessor has to be finite. Hence, every path 
in the constructed diagram is finite, and by Konig's lemma the entire diagram is finite.

At this stage we can deduce the equationality of the set $NR_s \cup D_1$ from the constructed diagram,
using a modification of the argument that was used in the proof of theorem 1.2. Recall that by theorems 2.5, 
2.9 and 2.13 of [Se3] the number of rigid or strictly
solid families of values of a rigid or solid limit group, 
that are  associated with a given value of the defining parameters,
is uniformly bounded by a bound that depends only on the rigid or solid limit
group (and not on the specific value of the defining parameters). 

Keeping our notation from the proof of theorem 1.2, let $depth$ be the number of levels in the 
diagram that we have
associated with the set, $NR_s \cup D_1$, and let $w$ be the maximal number
of vertices in a single level of the diagram. At each vertex in the diagram
we have placed a graded resolution. 
Each such graded resolution terminates in either a rigid or a
solid limit group, and by theorems 2.3, 2.9, and 2.13 in [Se3],  with each such
terminal rigid or solid limit group, there is a corresponding global bound on
the number of rigid or strictly solid families of values that is
associated with any possible value of the defining parameters. Let $b$ be the maximum
of all the global bounds that are associated with the terminal rigid or solid limit
groups of all the graded resolutions that are associated with the vertices in the
constructed universal diagram.

Let $q_1^0,q_2^0,\ldots$, be a given sequence of values (in the
coefficient group $F_k$) of the (free) variables $q$ in the set, $(NR_s \cup D_1)(p,q)$.
To prove equationality, we need to prove that the intersection, $\cap_{i=1}^{j} \, 
(NR_s \cup D_1)(p,q_i)$, strictly decreases for boundedly many indices $j$ (where the bound 
does not depend on the specific sequence $\{q_i\}$).  

First, we look at $q_1$ as parameters. There are at most
$w$ graded resolutions in the first level of the constructed universal diagram,
and  there are at most $b$ fibers that are associated with each of these graded
resolutions and with the value $q_1^0$. Hence, there are at most $w b$ fibers
of the graded resolutions in the first level of the diagram that are associated with
$q_1^0$.

The fibers that are associated with the value $q_1^0$, are either fibers of one of the completions,
$Comp(u,p,q,a)$,
or of one of the completions of the graded 
resolutions of the limit groups, $T_j(x_1,\ldots,x_s,p,q,a)$ (these collections of completions form the first level of the constructed
diagram).
Values in a  fiber of a completion, $Comp(u,p,q,a)$, are clearly in the Diophantine set,
$D_1(p,q)$. If a fiber in a graded resolution that is associated with a limit group,
$T_j(x_1,\ldots,x_s,p,q,a)$, contains a point, $(x_1^0,\ldots,x_s^0,p_0,q_1^0,a)$,
 for  which the the values, $(x_i^0,p_0,q_1^0,a)$, $i=1,\ldots,s$, are distinct  rigid values
of $Rgd(x,p,q,a)$, then the basic conditions that were imposed in constructing the Diophantine set $D_1$
do not hold for generic points in the fiber, hence, the basic conditions that were imposed in constructing $D_1$ may
hold only for points in the fiber that are contained in boundedly many fibers of graded resolutions that are associated
with $q_1^0$ and $D_1$, and have strictly smaller complexity than the original graded resolution of $T_j$. These last fibers are
contained in $D_1$, and therefore,
 the  entire fiber (or rather the restrictions of the
points in the fiber to the variables $(p,q)$)
is contained  in the definable set, $NR_s \cup D_1$.  

 If a fiber in a graded resolution that is associated with a limit group,
$T_j(x_1,\ldots,x_s,p,q,a)$, does not contain a  point, $(x_1^0,\ldots,x_s^0,p_0,q_1^0,a)$,
 for  which  the values, $(x_i^0,p_0,q_1^0,a)$, $i=1,\ldots,s$, are distinct  rigid specializations
of $Rgd(x,p,q,a)$, 
then we
omit this fiber from the (bounded) list of fibers that are associated with the value $q_1^0$. The set
$(NR_s \cup D_1)(p,q_1^0)$ is contained in the union of the (restrictions to the variables $(p,q)$ of points in the)
 remaining fibers. Therefore, after omitting all
such fibers, the set $(NR_s \cup D_1)(p,q_1^0)$ is precisely the (bounded) union of the remaining fibers.

We continue as we did in the proof of theorem 1.2.
If $((NR_s \cup D_1)(p,q_1^0)) \, \cap \, ((NR_s \cup D_1)(p,q_2^0)) \ = \ (NR_s \cup D_1)(p,q_1^0)$, there is no 
change, i.e., we remain with the same bounded collection of fibers that were associated with $q_1^0$.
 If 
 $((NR_s \cup D_1)(p,q_1^0)) \, \cap \, ((NR_s \cup D_1)(p,q_2^0)) \ \neq  \ (NR_s \cup D_1)(p,q_1^0)$, the intersection of the
two  sets is strictly contained in $(NR_s \cup D_1)(p,q_1^0)$. The  set $(NR_s \cup D_1)(p,q_1^0)$ is a 
finite union of at most $wb$ fibers of the graded resolutions in the first level of the diagram. Since
the intersection of the two  sets, $(NR_s \cup D_1)(p,q_1^0)$ and $(NR_s \cup D_1)(p,q_2^0)$,  
is strictly contained in $(NR_s \cup D_1)(p,q_1^0)$, it is a finite union of fibers - 
a proper (possibly empty) subset of
the fibers that are associated with $(NR_s \cup D_1)(p,q_1^0)$, and at least one of the fibers that  is associated with
$(NR_s \cup D_1)(p,q_1^0)$, that is replaced by a
(possibly empty)  finite collection of fibers  that are associated with the pair,
$(q_1^0,q_2^0)$, and with some of the graded resolutions that appear in the second level
of the constructed diagram. By the structure of the universal diagram, each fiber in the first
level can be replaced by at most $w  b$ fibers in the second level.

As we argue for fibers of graded resolutions in the first level that are associated with $q_1^0$, from the
bounded list of fibers that are associated with the pair $(q_1^0,q_2^0)$, we omit
fibers of graded resolutions in the second level that are associated with $(q_1^0,q_2^0)$, for which 
for generic values in these fibers (i.e., test sequences), either at least 
one of the values, $(x_i,p,q_1^0,a)$ or
$(y_i,p,q_1^0)$, $i=1,\ldots,s$, is flexible, or if some pair of these values is not distinct.

As we did in the proof of theorem 1.2, we repeat this argument iteratively. 
Each time a value $q_n^0$ is added, and
the corresponding intersection is a proper subset of the previous intersection,
at least one of the fibers that was associated with the intersection of the first $n-1$ values,
 is replaced by at most $w  b$ fibers in level that succeeds the level of that fiber (a fiber that is associated
with the last level of the diagram can only be  replaced by the empty set). As the digram has
$depth$ levels, the intersection:
$  \cap_{i=1}^{j} \ (NR_s \cup D_1)(p,q_i) $ 
can strictly decrease in at most 
  $2 \cdot depth \cdot {(wb)}^{depth}$ indices, which proves the equationality of the set $(NR_s \cup D_1)(p,q)$. 

\line{\hss$\qed$}

Proposition 1.6 proves that the set $NR_s \cup D_1$ is equational. To prove theorem 1.5
we continue iteratively. With the set $NR_s \cap D_1$ we associate a Diophantine set $D_2$,
precisely as we associated the Diophantine set $D_1$ with $NR_s$. 

Recall that $D_1$ was defined by finitely many completions, that we denote
$Comp(u,p,q,a)$. These completions were constructed by first collecting all the
non-degenerate values, $(x_1^0,\ldots,x_s^0,p,q,a)$, and associate with this collection of
values finitely many graded limit groups. Then we defined the completions that are associated with $D_1$
by further imposing a (basic) degeneration condition, applying theorem 1.3, and keeping only
those graded resolutions that have strictly smaller complexity than the completions that
they were constructed from (see the detailed description in the first part of the proof).
 
In order to define $D_2$, we start with the finitely many completions, $Comp(u,p,q,a)$, that
are associated with the Diophantine set $D_1$. We further look at all the values,
$(y_1^0,\ldots,y_s^0,u_0,p_0,q_0,a)$, where $(u_0,p_0,q_0,a)$ is a value of one of the
completions, $Comp(u,p,q,a)$, and the values, $(y_i^0,p_0,q_0,a)$, $i=1,\ldots,s$, are distinct
rigid values of the rigid limit group $Rgd(x,p,q,a)$. 

With this collection of values we canonically associate a finite collection of maximal limit
groups (according to theorem 7.2 in [Se1]), and we further apply theorem 1.3, and associate with 
 each limit group a finite collection of graded resolutions, such that the complexity of each 
graded resolution is bounded by the complexity of the completion that was used for its construction,
$Comp(u,p,q,a)$. 

At this point, we repeat what we did in constructing $D_1$. We look at all the values of
the completions of the obtained graded resolutions, for which either one of the values,
$(y_i^0,p_0,q_0,a)$, $i=1,\ldots,s$, is a flexible (non-rigid) value of $Rgd(x,p,q,a)$, or 
at least two of these (rigid) values are not distinct. These degenerate values of the
obtained completions satisfy one of finitely 
additional basic conditions (non-trivial equations), and with them we associate finitely many limit
groups, and by applying theorem 1.3, we further associate with them finitely many graded resolutions.

By theorem 1.3, the complexity of each of the constructed graded resolutions is bounded by the complexity
of the completion that was used in its construction,$ Comp(u,p,q,a)$.  As we did in the construction of
$D_1$, and since as in constructing $D_1$ the new completions were obtained by forcing additional basic
conditions,  we keep only those graded resolutions that have strictly smaller complexities than the completions,
$Comp(u,p,q,a)$, that they were constructed from.

We define $D_2$ to be  the Diophantine set that is the union of the Diophantine sets that are defined by the 
completions of those of the constructed graded resolutions that have strictly smaller complexity
than the completion, $Comp(u,p,q,a)$, that they were constructed from (note that this last completion was
associated with $D_1$.  
By construction, $D_2 \subset D_1$, and the complexities of the resolutions
that are associated with $D_2$, are   strictly smaller than the complexities of the corresponding resolutions that
are associated with $D_1$ (the definition of the complexity of a minimal rank resolution appears 
in definition 1.16 in [Se5]). 

By the same argument that was used to prove proposition 1.6, the
set $(NR_s \cap D_1) \cup D_2$ is equational. We continue to the third step with the set
$NR_s \cap D_1 \cap D_2 = NR_s \cap D_2$ and treat it exactly in the same way. By the
descending chain condition for complexities of minimal rank resolutions (cf. theorem 1.18 in [Se5]), this iterative
process terminates after finitely many steps, and the finite termination finally implies that the original set,
$NR_s$, is in the Boolean algebra of a collection of (minimal rank) equational sets, so theorem 1.5 follows.

\line{\hss$\qed$}

Essentially the same argument that was used to prove theorem 1.5 for the sets $NR_s$,
that are associated with minimal rank rigid graded limit groups,
can be used to prove a similar statement for sets of parameters for which a 
minimal rank solid limit
group admits at least $s$ strictly solid families of specializations. 

\vglue 1pc
\proclaim{Theorem 1.7} Let $F_k=<a_1,\ldots,a_k>$ be a non-abelian free group, and let
$Sld(x,p,q,a)$ be a  solid limit group,
 with respect to the parameter subgroup $<p,q>$.
Let $s$ be a positive integer, 
and let $NS_s$ be the set of values of the
defining parameters $<p,q>$ for which the solid limit group, $Sld(x,p,q,a)$, has
at least  $s$ strictly solid families of specializations.

There exists a collection of (minimal rank) equational
sets, so that the Boolean algebra generated by this collection contains the sets 
$NS_s$,
for every minimal rank solid limit group $Sld(x,p,q,a)$,
and every possible integer $s$. 
\endproclaim

\nfp The proof is similar to the proof of theorem 1.5. With the set $NS_s$
we associate a Diophantine set $D_1$, for which $NS_s \cup D_1$ is equational. We further argue 
that the intersection $NS_s \cap D_1$ is simpler than the original set $NS_s$.
We continue iteratively, precisely as we did in the proof of theorem 1.5.

To construct the Diophantine set $D_1$,
we look at the entire collection of values, $(x_1,\ldots,x_s,p,q,a)$, for
which the values, $(x_i,p,q,a)$, $1 \leq i \leq s$, belong to distinct strictly solid families. 
With this collection we
associate a canonical finite collection of maximal limit groups, that we view as graded
with respect to the parameter subgroup $<q>$. Since the solid limit group, $Sld(x,p,q,a)$, is of minimal rank, so
are all the maximal limit groups that we associated with the given collection of values. 
With these graded limit groups we further associate
the (graded) resolutions that appear in their strict graded Makanin-Razborov diagrams, and the
resolutions that are associated with the various strata in the singular loci of the diagram.

Given a (graded) completion, 
$Comp(x_1,\ldots,x_s,p,z,q,a)$, 
of one of these graded resolutions, we look at all the specializations
of the  completion for which either (at least) one of the
values, $(x_i,p,q,a)$, $i=1,\ldots,s$, that is supposed to be strictly solid is not strictly solid, or two such
values belong to the same strictly solid family. Note  that if such value 
is not strictly solid, or if two such values belong to the same strictly solid family, then
the (ambient) specialization of the given completion has to satisfy at least one of finitely
many (fixed) Diophantine conditions that are associated with the given solid
limit group, $Sld(x,p,q,a)$ (see definition 1.5 in [Se3] for these Diophantine conditions). 

With this collection of 
specializations of the (finitely many) completions, that are extended by values of extra variables that are  
added to demonstrate the validity of the Diophantine 
conditions they satisfy, we canonically associate a finite collection of graded limit groups. Each of these
maximal graded limit group has to
be of minimal rank, since the completions are of minimal rank, and by the structure of the additional Diophantine
conditions (see definition 1.5 in [Se3]). 
We further start with each of these maximal graded limit groups, and apply theorem 1.3 to
associate finitely many minimal rank
graded resolutions with each of the (finitely many) maximal graded limit groups that is associated with the
given collection of values. By theorem 1.3,
the complexity
of each of the constructed resolutions is bounded above by the complexity of the completion that was used for its construction,
and in case of equality in complexities, a constructed  resolution has to have the structure of a graded closure of
the completion from which it was constructed (part (3) in theorem 1.3).

The definition of the set $D_1$ in the solid case slightly differ from its definition in the rigid case.
In the solid case, we define the Diophantine set $D_1$ to be  the disjunction of the Diophantine sets that are associated with
the completions of all the constructed graded
resolutions  that either have strictly smaller complexity than
the completion they were constructed from,   or they are proper closures of the completion from which they were constructed
(i.e., non-trivial roots are added to abelian vertex groups that are associated with the completion from
which they were constructed).

By precisely the same argument that 
was used to prove proposition 1.6 (in the rigid case),
the set $NS_s \cup D_1$ is equational.

As in the proof of theorem 1.5, we continue  by analyzing the set
$NS_s \cap D_1$. With this set we associate a Diophantine set $D_2$ in a similar way to the construction
of the Diophantine set $D_1$. 
$D_2$ is a union of finitely many Diophantine sets, 
that are associated with   completions
of resolutions that have strictly smaller complexities than the complexities of the 
corresponding completions and closures
that define the set $D_1$, together with some proper closures of the completions and closures that define
$D_1$.  
By the same argument that was used to prove proposition 1.6,
$(NS_s \cap D_1) \cup D_2$ is equational.

We continue iteratively,
precisely as we did in proving theorem 1.5 in the
rigid case. As we argued in proving the termination of the construction of the diagram that was used in proving
the equationality of minimal rank Diophantine sets (theorem 1.2), given a completion that is associated with the Diophantine set
$D_1$, there is a  bound (that depends only on $D_1$) on the indices of supergroups of abelian vertex groups that
are associated with that completion, 
in all the completions of graded resolutions that are used to define any of the sets $D_n$, and are proper closures
of a completion that was used to define $D_1$.
 Hence, 
at some step $n_0$, all the completions and closures that define the   
Diophantine set, $D_{n_0}$, have strictly smaller complexity than the maximal complexity of the completions
and closures that define the Diophantine set $D_1$. Continuing with this argument iteratively, and combining it
with  the d.c.c.\ for complexities of minimal
rank resolutions ([Se5],1.18), guarantees that the iterative process (of corrections with minimal rank Diophantine sets)
for the analysis of the set $NS_s$
terminates after finitely many steps,
and the finite termination implies that the sets $NS_s$ are in the Boolean algebra of (minimal rank) equational sets.

\line{\hss$\qed$}

Theorem 1.2 proves that in the minimal  rank case
Diophantine sets are equational. Theorems 1.5 and 1.7 prove that 
sets of values of the defining parameters, for which a minimal rank rigid or solid limit group have at least $s$ rigid 
(strictly solid families of)
values, are in the Boolean algebra of equational sets. Before we analyze general minimal rank definable
sets, we need to analyze the (definable) 
set of values of the defining
parameters, for which a given (finite) collection of covers of a  graded resolution forms a covering
closure  (see definition 1.16 in [Se2] for a covering closure).

\vglue 1pc
\proclaim{Theorem 1.8} Let $F_k=<a_1,\ldots,a_k>$ be a non-abelian free group, let $G(x,p,q,a)$ be a graded
limit group (with respect to the parameter subgroup $<p,q>$), and let $GRes(x,p,q,a)$ be a 
graded strict resolution
of $G(x,p,q,a)$ that terminates in the rigid (solid) limit group, $Rgd(x,p,q,a)$ ($Sld(x,p,q,a)$). 
%
Suppose that the terminating rigid (solid) limit group, $Rgd(x,p,q,a)$ ($Sld(x,p,q,a)$), 
is of minimal (graded) rank.

Let $GCl_1(z,x,p,q,a),\ldots,GCl_t(z,x,p,q,a)$ be a given set of graded closures of $GRes(x,p,q,a)$. Then the
set of values of the parameters $<p,q>$ for which the given set of the (associated fibers of the)
graded closures forms a covering closure
of the (associated fibers of the) graded resolution $GRes(x,p,q,a)$, that we denote, $Cov(p,q)$, is in the Boolean algebra of equational sets.
\endproclaim

\nfp The proof is similar to the proofs of theorems 1.5 and 1.7. The set $Cov(p,q)$ is defined to be the set 
of values of the defining parameters, $<p,q>$, for which the fibers that are associated with the given 
(finite) set of closures and the given values of the parameters, form a covering closure of the fibers that 
are  associated with the graded
resolution, $GRes(x,p,q,a)$  (and the given values of the parameters). As in the proofs of theorems
1.5 and 1.7, with $Cov(p,q)$
we associate a minimal rank Diophantine set $D_1$, for which $Cov(p,q) \cup D_1$ is equational, and $D_1$
and $D_1 \cap Cov(p,q)$ are simpler than
$Cov(p,q)$, in a similar way to what was shown in theorems 1.5 and 1.7.

To analyze the set $Cov(p,q)$ and construct the Diophantine set $D_1$,
we look at the collection of values: 
$$(x_1^0,\ldots,x_s^0,y_1^0,\ldots,y_m^0,r_1^0,\ldots,r_s^0,p_0,q_0,a)$$ for
which:
\roster
\item"{(i)}" for the tuple $(p_0,q_0)$ there exist precisely $s$ distinct rigid (strictly solid families of) 
specializations
of the rigid (solid) limit group, $Rgd(x,p,q,a)$ ($Sld(x,p,q,a)$), and at least (total number of) 
$m$ distinct rigid 
and strictly solid families of specializations  of the terminal (rigid and solid) limit groups of the
closures: 
$GCl_1(z,x,p,q,a),\ldots,GCl_t(z,x,p,q,a)$.

\item"{(ii)}" in case the terminal limit groups of $GRes$ is rigid, the values, $(x_i^0,p_0,q_0,a)$, 
$i=1,\ldots,s$, denote  the distinct rigid 
values of $Rgd(x,p,q,a)$. In case the terminal limit group of $GRes$ is solid, the values, $(x_i^0,p_0,q_0,a)$, $i=1,\ldots,s$,  belong to 
the $s$ distinct strictly solid families of $Sld(x,p,q,a)$. 

\item"{(iii)}" the values, $(y_j^0,p_0,q_0,a)$, $j=1,\ldots,m$, are either distinct rigid values or belong to distinct strictly solid 
families of values of the terminal limit groups of the closures: $GCl_1,\ldots,GCl_t$.

\item"{(iv)}" the values $r_i^0$'s  are  added only in case the terminal limit group
of $GRes$ is solid. In this case the values, $r_i^0$,  demonstrate that the fibers that are 
associated with the given closures and the values, $y_1^0,\ldots,y_m^0$, form a covering closure of
the fibers that are  associated with the resolution $GRes$ and the (strictly solid) values, 
$x_1^0,\ldots,x_s^0$.   
These include values of primitive roots of the specializations of all the non-cyclic abelian groups, and edge groups,
in the abelian decomposition that is associated with the solid terminal limit group of $GRes$, $Sld(x,p,q,a)$, and
values of elements that demonstrate that multiples of these primitive roots up to the least common multiples
of the indices of the finite index subgroups that are associated with the graded closures, $GCl_1,\ldots,GCl_t$,
do belong to the fibers that are associated with the values, $y_1^0,\ldots,y_m^0$, and their
corresponding closures (cf. section 1 of [Se5] in which we added similar values of elements,  to indicate
that a proof statement is a valid proof statement).
\endroster

We look at the collection of all such values for all the possible values of $s$ and $m$ (note that 
$s$ and $m$ are bounded, since the number of rigid values of a rigid limit group, and 
the number of strictly solid families of values of a solid limit group, that are associated with a given
value of the defining parameters are globally bounded by theorems 2.5, 2.9 and 2.13 in [Se3]).

With this collection of values we
associate a canonical finite collection of maximal limit groups, that we view as graded
(limit groups) with respect to the parameter subgroup $<q>$. With these graded limit groups we associate
the (graded) resolutions that appear in their strict graded Makanin-Razborov diagrams, and the
strict 
resolutions that are associated with the various strata in the singular loci of the diagrams. 
Since we assumed that
the terminal limit group $Rgd(x,p,q,a)$ ($Sld(x,p,q,a)$) are of minimal (graded) rank, 
all the resolutions in these graded Makanin-Razborov diagrams are of minimal (graded) rank (i.e., all the limit groups 
that appear along these graded resolutions are of minimal (graded) rank).

Given a (graded) completion, 
$Comp(x_1,\ldots,x_s,y_1,\ldots,y_m,r_1,\ldots,r_s,z,p,q,a)$, 
of one of these graded resolutions, we look at all the specializations
of  the completion for which either:
\roster
\item"{(1)}" (at least) one of the
values, $x_1^0,\ldots,x_s^0,y_1^0,\ldots,y_m^0$, that is supposed to be rigid or strictly solid is not rigid
or not strictly solid. 

\item"{(2)}"  two of these values that are supposed to be rigid and distinct coincide, or two of these values that
are supposed to be values in distinct strictly
solid families

\noindent
Note  that (as in the proofs of theorems 1.5 and 1.7)
the condition that  a given value 
is not rigid or not strictly solid, or that two values are equal or belong to the same strictly solid family, 
translates into one of finitely many  Diophantine conditions that the specializations of the given 
completion have to satisfy (see definition 1.5 in [Se3] for the definition of these Diophantine conditions).

\item"{(3)}" a value of what is supposed to be a primitive root, $r_i^0$, has a root of order
that is not co-prime to the least common multiple of the indices of the  
finite index subgroups that are associated
with the corresponding graded closures, $GCl_1,\ldots,GCl_t$. Note that once again (as in part (2)), 
this condition translates into
one of finitely many Diophantine conditions that the ambient specializations of the given completion
have to satisfy.

\item"{(4)}" there exists an extra rigid (strictly solid family of) value(s) of the rigid limit
group $Rgd(x,p,q,a)$ ($Sld(x,p,q,a)$), in addition to those specified by the values, $x_1^0,\ldots,x_s^0$. 
\endroster

With this collection of 
specializations of the (finitely many) completions of graded resolutions in the constructed Makanin-Razborov
diagrams,
in addition with values of extra variables that are being 
added to demonstrate the the specializations of the completions satisfy one of the the Diophantine 
conditions (1)-(3), or the existence of an extra rigid or strictly solid value (condition (4)),
we canonically associate a finite collection of graded limit groups. We further apply
 theorem 1.3, and associate finitely many graded resolutions
with these graded limit groups that we denote $DGRes$.
By theorem 1.3, the complexity
of each of the associated resolutions, $DGRes$, is bounded above by the complexity of the completion from which
they were constructed,
and in case of equality in complexities, an associated  resolution $DGRes$ has to be a graded closure of
the completion from which it was constructed.



At this point we look at the subcollection of  graded closures of the original completions,
$Comp(x_1,\ldots,x_s,y_1,\ldots,y_m,r,z,p,q,a)$, that were constructed from values that are obtained
from specializations of one of these completions and values of an extra rigid or strictly 
solid value of the terminal rigid or solid limit group of the given resolution $GRes$,
i.e., that are constructed according to case (4). With each such graded closure, that
we denote, $BCl$, we
associate an additional collection of graded minimal rank resolutions.

We 
 collect all the  of specializations of each of the closures,
$BCl$,  for which 
 the restriction to the value of the  elements that represent the extra rigid or strictly solid value,  is either
flexible (i.e., not rigid) or it coincides with one of the rigid values, $x_1^0,\ldots,x_s^0$,
that are the restrictions of the specialization of the completion from which 
the closure $BCl$ was constructed, 
 or it
is not strictly solid or
belongs to one of the strictly solid families that are associated with the values, $x^0_1,\ldots,x^0_s$ (in
the solid case). Note that these degenerations of the extra rigid or strictly solid value can be 
enforced by one of finitely
many Diophantine conditions, as we did in cases (1) and (2).

With these values that are obtained from specializations of the closures, $BCl$, and values
of elements that demonstrate  that one of the 
extra Diophantine conditions that are imposed on these specializations is fulfilled, 
we canonically associate a finite
collection of maximal limit groups (according to theorem 7.2 in [Se1]). By theorem 1.3, with these limit groups
we can associate a finite collection of graded resolutions, that we denote $EFRes$. 
The complexities of these graded resolutions, $EFRes$, are
bounded by the complexities of the corresponding closures, $BCl$, and in case of equality, a corresponding
graded resolution is a graded closure of (the closure) $BCl$, hence, a graded closure of the  completion,
$Comp(x_1,\ldots,x_s,y_1,\ldots,y_m,r,z,p,q,a)$, from which $BCl$ was constructed.

With the set, 
$Cov(p,q)$, we first associated finitely many completions, 
$Comp(x_1,\ldots,x_s,y_1,\ldots,y_m,r,z,p,q,a)$, that were constructed from values that satisfy 
properties (i)-(iv). With these completions we further associated finitely many graded resolutions, $DGRes$,
by extending the specializations of these completions to values that satisfy one of the properties
(1)-(4). By theorem 1.3 the complexities of these graded resolutions, $DGRes$, are bounded by the complexity of
the completion from which they were constructed. We denoted those of the constructed graded resolutions that were
constructed from values that satisfy part (4), and are graded closures of the completions from which they
were constructed by $BCl$. With each graded closure $BCl$ we further associated a collection of graded
resolutions that we denoted $EFRes$, and in which the values that correspond  to the extra rigid or strictly solid
element are degenerate.

We define the Diophantine set $D_1$ to be the disjunction of all the Diophantine sets that are associated
with completions of (the constructed)  graded
resolutions, $DGRes$ and $EFRes$,  that have either strictly smaller complexity than
the completion they were constructed from, or they are proper graded  closures of the completions that they 
were constructed from (recall that a proper graded closure is a  closure in which proper roots were added to
some of the abelian vertex groups that are associated with the various levels of the completion from which the
closure was constructed).

Note that to analyze  the sets, $NR_s$ and $NS_s$, we started with their configuration limit groups and the 
completions of the resolutions in their Makanin-Razborov diagrams. The degeneracy of a configuration
homomorphism can be expressed by a basic condition (in the rigid case) or by a Diophantine condition (in the
solid case). If a non-proper closure of such a completion satisfies the non-degeneracy basic or Diophantine 
condition, then the entire fibers that are associated with this closure can be removed and ignored when we 
analyze the sets $NR_s$ or $NS_s$. However, when we analyze the set $Cov(p,q)$, it may be that a non-proper
closure satisfy the degeneracy condition (4), and still there will be values in the corresponding fiber
that restrict to values of the defining parameters $p,q$ that are in the set $Cov(p,q)$, i.e., values
for which the degeneracy condition (4) collapses. Precisely for this reason we need to construct the
resolutions, $EFRes$, and add the Diophantine sets that are associated with their completions to the 
(correcting) set
$D_1$. 

By the same argument that 
was used in proving proposition 1.6 and theorem 1.7,
the set $Cov(p,q) \cup D_1$ is equational. As in proving theorems 1.5 and 1.7,  
we continue  by analyzing the set
$Cov(p,q) \cap D_1$. The rest (iterative continuation) of the argument is 
identical to the one that is used in proving theorems 1.5
and 1.7.

\line{\hss$\qed$}

Proving that (minimal rank) Diophantine sets are equational, that 
(in minimal (graded) rank)
sets for which a rigid or solid limit group have at least $s$ rigid (strictly solid families of)
values, are in the Boolean algebra of equational sets, and that the set of values
of the defining parameters for which a given set of closures forms a covering closure of a
given graded resolution (assuming its terminating rigid or solid limit group is of minimal
(graded) rank), is in the Boolean algebra of equational sets, we are finally ready to prove the main theorem
of this section, i.e., that minimal rank definable sets are in the Boolean algebra of equational sets.

\vglue 1pc
\proclaim{Theorem 1.9} Let $F_k=<a_1,\ldots,a_k>$ be a non-abelian free group, and let $L(p,q)$ be a
minimal rank definable set (see definition 1.1). Then $L(p,q)$ is in the Boolean algebra of equational sets.
\endproclaim

\nfp To analyze the minimal (graded) rank set $L(p,q)$, we use the precise description of a definable set 
that was obtained using the sieve procedure for quantifier elimination that is presented in 
[Se5] and [Se6]. 
The quantifier elimination procedure is long and uses a long list of objects and terms that we can not
present here  in detail. In the minimal rank it is described in detail in section 1 of [Se5]. We use the 
terminology that is presented and used in this section in [Se5].

Recall that
with the set $L(p,q)$ the sieve procedure associates a finite collection of graded PS resolutions
that terminate in rigid
and solid limit groups (with respect to the parameter subgroup $<p,q>$), 
and with each such graded resolution it associates a finite collection of graded
closures of these resolutions that contains Non-Rigid, Non-Solid, Left, Root, Extra PS, and 
collapse extra PS resolutions
(see definitions  1.25-1.30 of [Se5] for the exact definitions of these resolutions). 

By the construction of the sieve
procedure, 
since the definable set $L(p,q)$ is assumed to be
of minimal (graded) rank, all the terminating rigid and solid limit groups of the PS resolutions that
are associated
with $L(p,q)$ by the sieve procedure are of minimal (graded) rank as well. 

As there are finitely many PS resolutions that are associated by the sieve procedure
with the definable set $L(p,q)$, with any given value of the defining parameters $p,q$ there can be at most
boundedly many  fibers that are associated with a given value of $p,q$ and with  one of the PS resolutions 
that are associated with $L(p,q)$ (see definition 1.4 for a fiber of a graded resolution). 

\noindent
By the sieve procedure, that eventually leads to
quantifier elimination over a free group, the definable set $L(p,q)$ is precisely the union of
those values of the defining parameters $p,q$, for which:
\roster
\item"{(1)}" there exists a fiber of one of the (finitely
many) PS resolutions that are associated with $L(p,q)$, and is associated 
with the given value of the parameters
$p,q$.

\item"{(2)}" this fiber is not covered by the bounded collection of fibers that are associated with the 
given value of $p,q$ and with the (finite) collection of   
Non-Rigid, Non-Solid, Left, Root and extra PS resolutions, minus
the fibers that are associated with the
collapse extra PS resolutions (see definition 1.16 in [Se2] for a covering closure).
\endroster

Let $PSRes_i$, $i=1,\ldots,r$, be  the finitely many PS resolutions that are associated with the given minimal rank definable set
$L(p,q)$. For each index $i$, $i=1,\ldots,r$, let $Rgd_i(x,p,q,a)$ ($Sld_i(x,p,q,a)$) be the terminal rigid (solid) limit group of $PSRes_i$.
With the PS resolution $PSRes_i$ and its terminal rigid or solid limit group $Rgd_i$ or $Sld_i$, we associate the definable
set, $NR^i_1(p,q)$ or $NS^i_1(p,q)$, that defines those values of the defining parameters $p,q$ that
extend to  rigid or  strictly solid values of $Rgd_i$ or $Sld_i$. By theorems 1.5 and 1.7 the sets $NR^i_1$ and $NS^i_1$ are
in the boolean algebra of equational sets.

With each of the PS resolutions, $PSRes_i$, 
the sieve procedure  associates a finite collection of graded closures of
it that contains Non-Rigid, Non-Solid, Left, Root, Extra PS, and 
collapse extra PS resolutions. With the graded resolution $PSRes_i$, and its given set of closures, we associate
a definable set $Cov_i(p,q)$, that defines those values of the defining parameters $p,q$ for which the associated
fibers of $PSRes_i$ that are associated with the value $p,q$ are covered by the fibers that are associated with the
given finite set of closures of it and with the value of $p,q$.  By theorem 1.8 $Cov_i(p,q)$ is in the Boolean algebra of
equational sets.

By the sieve procedure, as indicated by (1) and (2) above, the definable set $L(p,q)$ is the finite union:
$$\cup_{i=1}^r \, NR^i_1(p,q) \, (NS^i_1(p,q)) \, \setminus \, Cov_i(p,q)$$
In particular, $L(p,q)$ is a Boolean combination of the sets $NR^i_1$ ($NS^i_1$) and $Cov_i$. Since by theorems 1.5, 1.7 and 1.8, 
the sets, $NR^i_1$, $NS^i_1$ and $Cov_i(p,q)$, are all in the Boolean algebra of equational sets, so is their Boolean
combination $L(p,q)$, and theorem 1.9 follows.

\line{\hss$\qed$}

\vglue 1.5pc
\centerline{\bf{\S2. Diophantine Sets}}
\medskip
Our first step in approaching the stability
of free (and hyperbolic) groups, is proving that Diophantine sets are equational.
This was proved in theorem 1.2 in the minimal rank case, and is more involved though
still valid in general.

\vglue 1pc
\proclaim{Theorem 2.1} Let $F_k=<a_1,\ldots,a_k>$ be a non-abelian free group, and let
$$D(p,q) \ = \ \{ \, (p,q) \, ¦ \, \exists x \  \Sigma (x,p,q,a)=1 \, \}$$
be a Diophantine set defined over $F_k$. Then $D(p,q)$ is equational.
\endproclaim

\nfp With the system of equations $\Sigma(x,p,q,a)=1$ we associate its graded 
Makanin-Razborov diagram (with respect to the parameter subgroup $<p,q>$), 
and we look at the (finite) collection of rigid limit groups, $Rgd(x,p,q,a)$, and
solid limit groups, $Sld(x,p,q,a)$, along the diagram. By the properties of the
graded diagram, the Diophantine set $D(p,q)$ is precisely the collection of 
values of the parameter subgroup $<p,q>$, for which at least one of the 
rigid or solid limit groups along the graded Makanin-Razborov diagram of $\Sigma(x,p,q,a)$ admits
a rigid or a strictly solid value.

To prove the equationality of a general Diophantine set
$D(p,q)$, we associate with it a finite diagram, similar but somewhat different
to the one we associated with a minimal rank Diophantine set in proving theorem 1.2.
To prove the termination of the iterative procedure that is used for the construction of
the diagram, we apply the techniques that were used in proving the termination of the sieve procedure
that was used in obtaining quantifier
elimination in [Se6].

We start the construction of the diagram by collecting all the values of the
tuple, $(x,p,q,a)$, that are rigid or strictly solid values of one of the rigid
or solid limit groups that appear along the graded Makanin-Razborov diagram
of the system $\Sigma(x,p,q,a)$. With this collection of values, we associate its Zariski closure,
that by theorem 7.2 in [Se1] is dual to  a canonical finite
collection of maximal limit groups, that we denote $L_i(x,p,q,a)$.

With a limit group $L_i(x,p,q,a)$, viewed as a graded limit
group with respect to the parameter subgroup $<q>$, we associate its 
taut graded Makanin-Razborov
diagram (see section 2 in [Se4] for the construction of the taut diagram of a limit group).
With each resolution in the taut Makanin-Razborov diagram, we further associate
its singular locus, and the graded resolutions that are associated with
each of the strata in the singular locus. We conclude the first step of the
construction of the diagram, by associating the (graded) completion  with each of the
graded resolutions in our finite collection (of resolutions), that we denote, $Comp(z,x,p,q,a)$.

We continue to the construction of the second step of the diagram with each of the completions 
$Comp(z,x,p,q,a)$ in parallel. With each such completion we associate the collection of
values: $(x^0_2,z_0,x^0_1,p_0,q^0_1,q^0_2,a)$ ,for which:
\roster
\item"{(1)}" $(z_0,x^0_1,p_0,q^0_1,a)$ factors through the completion, $Comp(z,x,p,q,a)$, and
$(x^0_1,p_0,q^0_1,a)$ is rigid or strictly solid with respect to  one of the rigid or solid
limit groups in the graded diagram of $\Sigma(x,p,q,a)$.

\item"{(2)}" $(x^0_2,p_0,q^0_2,a)$ is a rigid or a strictly solid value of
 one of the rigid or solid limit groups in the graded Makanin-Razborov diagram of 
$\Sigma(x,p,q,a)$. In case it is strictly solid, it is the shortest in its strictly solid family.
\endroster

First, for each completion, $Comp(z,x,p,q,a)$, that is placed in the initial level of the diagram,
 we collect all its test sequences that extend to values that satisfy properties (1) and (2). By the techniques
that were used in constructing formal limit groups (section 3 in [Se2]), with these test sequences and their extended
values it is possible to associate (canonically) a finite (possibly empty) collection of graded limit groups that
have a similar structure as (graded) closures of
the completions $Comp(z,x,p,q,a)$.

With each of the completions that are placed in the initial level of the diagram, $Comp(z,x,p,q,a)$, 
we associate the collection of all the 
sequences: $$\{(x_2(n),z(n),x_1(n),p(n),q_1(n),q_2(n),a)\}_{n=1}^{\infty}$$
so that for each $n$, the corresponding value satisfies  conditions (1) and (2), and
the sequence: $\{(z(n),x_1(n),p(n),q_1(n),a)\}_{n=1}^{\infty}$ forms a (graded) test sequence with respect to
the given (graded) completion $Comp(z,x,p,q,a)$. In addition we require that for every index $n$, the lengths of the values
of fixed 
set of generators of the vertex groups in the abelian decompositions that are associated with  
all the levels of the completion, $Comp(z,x,p,q,a)$, except for its terminal level, are at least $n$ times longer than 
the lengths of the values $q_2(n)$.

By the techniques that are used to analyze graded formal limit groups,
that are presented in section 3 of [Se2], with this collection of sequences it is possible to canonically associate a finite
collection of limit groups that have the same structure as (graded) closures of the initial completion, $Comp(z,x,p,q,a)$, 
through which they all factor. These limit groups are obtained from the completion $Comp(z,x,p,q,a)$ by possibly adding
roots to abelian vertex groups in the abelian decompositions that are associated with the various levels of
the completion, $Comp(z,x,p,q,a)$, and replacing the terminal rigid or solid limit group of $Comp(z,x,p,q,a)$ (with respect
to the parameter subgroup $<q>$),  with
a rigid or solid limit groups with respect to the parameter subgroup $<q_1,q_2>$. 

\noindent
We will denote the finitely many  limit groups  
that are associated with all these sequences, 
$DQCl_i(x_2,z,x_1,p,q_1,q_2,a)$ (note that they  are graded with respect to the parameter 
subgroup $<q_1,q_2>$). 

We further look at all the values, $(x_2,z,x_1,p,q_1,q_2,a)$, that satisfy
conditions (1) and (2), and do not factor through any of the (finite set of) limit groups,
$DQCl_i(x_2,z,x_1,p,q_1,q_2,a)$. 
By section 5 in [Se1], with this collection of values we can
associate a canonical finite collection of maximal limit groups, that we denote $\{M_j(x_2,z,x_1,p,q_1,q_2)$\}, which we
view as graded limit groups with respect to the parameter subgroup $<q_1,q_2>$. 

Using the iterative procedure for the construction
of (quotient) resolutions, that is used in each step of the 
sieve method for quantifier elimination and presented in [Se6], we associate with
this collection of values and limit groups, $\{M_j\}$, finitely many multi-graded resolutions with respect to
the defining parameters $<q_1,q_2>$, and with each such graded resolution we
associate its finitely many  core resolutions, anvils, developing resolutions, and
(possibly) sculpted resolutions and carriers (see (the first step in) [Se6] for a detailed description of the
iterative construction of the multi-graded resolutions and the anvils and developing resolutions that are attached to them). 

\noindent
Note that in the sense of the sieve procedure that is presented in [Se6], each of the constructed anvils, 
has a smaller complexity than the completion, $Comp(x,z,p,q,a)$, that is associated with it (i.e., the completion with which
we have started this branch). 

At the vertices in the second level of the diagram we place the finite collection of anvils that were constructed from
the limit groups, $\{M_j\}$, and with each anvil we associate the (graded) completion of its developing resolution.
In the other vertices in the second level we place  those limit groups $DQCl_i$,
for which proper roots were added to (abelian vertex groups in) the completion, $Comp(z,x,p,q,a)$, from which they were constructed 
(by construction, each of the groups $DQCl_i$ has in particular a structure of a completion).
Note that with each vertex in the second level there is an associated completion, either one of the groups $DQCli_i$ or a completion of
the developing resolution of the associated anvil).
Each  vertex in the second level, is connected by a directed  edge that points to it and starts at a vertex in the first level of the
diagram in which the completion,
$Comp(z,x,p,q,a)$, that was used in its construction,  is placed. 

\smallskip
We continue iteratively. With each vertex in level $s$ there in an associated completion. This completion is either a limit group
that is obtained from a completion of level $s-1$ by adding proper roots to some of the abelian vertex groups that are associated 
with its various completions, or it is the completion of a developing resolution of an anvil that was constructed in step $s$ of
the procedure, according to the general step of the sieve procedure [Se6].

Given a completion that is placed in level $s$, and its associated developing resolution and anvil, 
we look at the collection of  sequences of
values: 
$$(x_{s+1}(n),w(n),p(n),q_1(n),\ldots,q_{s+1}(n),a)$$ 
for which the values $(x_{s+1}(n),p(n),q_{s+1}(n))$ are rigid or strictly solid with respect to one
of the (finitely many) rigid or solid limit groups that are associated with our given Diophantine set,
and the corresponding restricted values, $(w(n),p(n),q_1(n),\ldots,q_s(n),a)$,   form
a test sequence of the given completion.
Given all these sequences, we apply the techniques for the construction of (graded) formal limit groups, that are presented in section 3
of [Se2], and associate with the given completion a finite set of limit groups that are obtained from the given completion by 
(possibly) adding roots to some of the abelian vertex groups that are associated with its various levels, and replacing its terminal
limit group, with a rigid or solid limit groups that is graded with respect to the parameter subgroup, $q_1,\ldots,q_{s+1}$.
We denote these limit groups $DQCl$.
 As in the second step of the construction, in the $s+1$ level of the diagram, we place only those limit groups to which proper roots
were added to the abelian vertex groups that are associated with the various levels of their associated  completion from level $s$.

At this point we look at all the values:
$$(x_{s+1}(n),t(n),p(n),q_1(n),\ldots,q_{s+1}(n),a)$$ 
for which:
\roster
\item"{(1)}" the restricted values, $(x_{s+1}(n),p(n),q_{s+1}(n),a)$, are rigid or strictly solid values of  one
of the (finitely many) rigid or solid limit groups that are associated with our given Diophantine set. 

\item"{(2)}" with the completion in level $s$ there is an associated developing resolution, and an associated anvil.
 The  value, $(t(n),p(n),q_1(n),\ldots,q_s(n),a)$, is a value of the anvil that is associated with the completion.

\item"{(3)}" the value 
$(x_{s+1}(n),t(n),p(n),q_1(n),\ldots,q_{s+1}(n),a)$ restricts to a value:
$(x_{s+1}(n),w(n),p(n),q_1(n),\ldots,q_{s+1}(n),a)$, that further restricts to a value:
$(w(n),p(n),q_1(n),\ldots,q_{s+1}(n),a)$, which is a value of the completion from level $s$ that we have started with. We further assume
that the value:
$(x_{s+1}(n),w(n),p(n),q_1(n),\ldots,q_{s+1}(n),a)$, does not factor through any of the limit groups $DQCl$ that we have associated with 
the given completion from level $s$.
\endroster

By our standard techniques, that were presented in section 5 of [Se1], with this collection of 
values we associate its Zariski closure, and with it we canonically associate its dual 
 finite collection 
of (graded)
limit groups. Given these limit groups, 
we apply the construction that was used in the general step of  the sieve procedure and  presented in
[Se6], to construct a finite collection of multi-graded resolutions, with which there 
are associated
core resolutions, developing resolutions, (possible) sculpted resolution and carriers, and anvils.  
As in the second step of the procedure, we add to the vertices in level $s+1$, a finite collection of vertices, and in each such vertex
we place a completion of the developing resolution of an anvil that was constructed from one of the completions in level
$s$.

\vglue 1pc
\proclaim{Proposition 2.2} The iterative procedure that is associated with a Diophantine set
terminates after finitely many steps.
\endproclaim

\nfp To prove termination we use essentially the same argument that was used to prove 
the termination of the sieve procedure in ([Se6],22). Unfortunately, the sieve
procedure is long and technical, hence, we can not repeat even the definitions
of the objects that are constructed along it, and are used in proving
the termination. Therefore, for the rest of the proof we will assume that the
reader is familiar with the structure of the sieve procedure,  the objects
that are constructed along it, and the proof of its termination, that are all presented in [Se6].

Since our procedure is a
locally finite branching process, if it doesn't terminate it must contain
an infinite path. Given a completion that is placed in some vertex  in the diagram that we constructed,
we have associated with it finitely many limit groups that have a structure of a closure, hence,
there is a bound on the order of proper roots that we can add to the abelian vertex groups that are associated
with the various levels of the completion along the next steps of the procedure. Hence, given a completion that
is placed in a vertex of the diagram, we can start from it and continue along a path of the diagram that passes
only through limit groups that are obtained from the completion by adding proper roots to abelian vertex groups
that are associated with the completion for only finitely many steps, and then we must pass to an anvil that
was constructed from the last  closure of the completion according to the general step of the sieve procedure.

Since the  construction of the resolutions and the anvils
that we use is identical to the construction that is used in the sieve procedure, proposition 26 in [Se6] remains valid, 
i.e., given an infinite path of our procedure,
for each positive integer $m$, there exists a step $n_m$ and width $d_m$, so that the sculpted
and penetrated sculpted resolutions of width $d_m$ at step $n_m$ are all
eventual (i.e., they do not change along the rest of the infinite path),
and the number of $n_m$ sculpted resolutions of width $d_m$, $sc(n_m,d_m)$,
satisfies: $sc(n_m,d_m)=m$.

Therefore, as in theorem 27 in [Se6], to conclude the proof of theorem 2.2, 
i.e., 
to prove the termination of the procedure for the construction of the diagram
that is associated with a Diophantine set after finitely many steps, we need to
show the existence of a global bound on the number of eventual sculpted resolutions of the same
width that are associated with the anvils along an infinite path of the procedure. 

Our approach towards obtaining a bound on the number of eventual sculpted resolutions with
the same width along
an infinite path of the procedure, is essentially identical to the one that is used to prove theorem 27 in [Se6], 
and is based on the argument that was used to obtain a bound on the number of rigid and strictly solid families of
values of rigid and solid limit groups, that is presented in the first two sections of [Se3] (theorems
2.5,2.9 and 2.13 in [Se3]).

Recall that in proving theorem 27 in [Se6], we argued that if
there is no bound (independent of the width) on the number of eventual 
sculpted resolutions of the same
width that is associated with
an anvil along a given infinite path of the sieve procedure, then there must be two sequences of values:
of the same rigid or solid limit group: $\{(x_i(n),w(n),p(n),a)\}$ and $\{(x_j(n),w(n),p(n),a)\}$, for some $j>i$,
so that for every index $n$: $x_i(n)=x_j(n)$ in the rigid case, and $(x_i(n),p(n),a)$ belongs to the same
family of $(x_j(n),p(n),a)$ in the solid case. However, the values, $\{(x_j(n),(n),a)\}$,
are assumed to be either (extra) rigid or strictly solid, and the values,  $\{(x_i(n),p(n),a)\}$, are 
assumed to be flexible (or not strictly solid), and we got a contradiction, hence, 
we obtained a  bound on the number of eventual
sculpted resolutions of the same width.

Assume that our procedure for the construction of the diagram that is associated with a Diophantine set contains
an infinite path, and along this path there is no bound    
(independent of the width) on the number of eventual 
sculpted resolutions of the same
width that are associated with
an anvil along this given infinite path of the procedure.

First, we observe that in a  test sequence of each of the developing and sculpted resolutions that are associated with the
anvils along our iterative procedure, 
we may assume that the lengths of the parts of the
variables
$x_i(n)$, and the parameters $p(n)$, that do not belong to the terminal level of the
graded sculpted resolutions, are much bigger than the lengths of the values $q_i(n)$ (that are
assumed to be
part of the distinguished vertex group in the terminal level).
Therefore, in applying the argument that was used in proving theorem 27 in [Se6] to the sculpted resolutions
that were constructed along an infinite path of our procedure, in obtaining a sequence of compatible 
JSJ decompositions that are used in analyzing the
sequences, $\{(x_i(n),p(n),q_i(n),a)\}$ (see theorem 36 in [Se6] for the compatible JSJ decomposition), 
where these values are  restrictions of a test sequence of some developing resolution
along the infinite path, the subgroup $<q_i>$ remains  elliptic, until we 
reach the terminal level of the (eventual) sculpted resolution in question.

Hence, by applying the same argument that was used
to prove theorem 27 in [Se6], we obtain two sequences of rigid or strictly solid values of
the same rigid or solid limit group: $\{(x_i(n),p(n),q_i(n),a)\}$
and $\{(x_j(n),p(n),q_j(n),a)\}$, for some $j>i$, along a given infinite path,
that are compatible in all the levels except, perhaps, the terminal
level. 
This  contradicts  the assumption that along our iterative procedure, we collected 
rigid and strictly solid values that do not factor through the closures that were associated with 
completions of developing resolutions of anvils
that were constructed  in  previous steps of the procedure. Therefore, like in the proof of theorem 27 in [Se6] we get
a contradiction,  and hence we proved that there exists  a global bound on 
the number of eventual sculpted resolutions of the same width along an infinite path of our procedure (i.e., along a path in the
constructed diagram). This global bound contradicts
 the existence of an infinite path in the procedure for the construction
of the diagram that we associated with a Diophantine set, which finally  implies that the procedure for the construction of the diagram
terminates
after finitely many steps (see the proof of theorem 27 in [Se6] for a detailed description of the notions, constructions and arguments that we applied).

\line{\hss$\qed$}

Proposition 2.2 enables one to associate a finite diagram with a Diophantine family. The existence of
such diagram together with the existence of a global bound on the number of rigid and strictly solid values
of rigid and solid limit groups (for any given value of the parameter subgroup), that was proved in theorems
2.5, 2.9 and 2.13 in [Se3], enable us to conclude the equationality of Diophantine families.

Let $D(p,q)$ be a Diophantine family. Let $q_1,\ldots,q_n$ be a sequence of values of the parameters of
the family $D(p,q)$, for which the intersections $\cap_{j=1}^m D(p,{q_j})$, is a strictly decreasing sequence for
$m=1,\ldots,n$. 

\noindent
With the Diophantine family $D(p,q)$ we associate its diagram, that we denote $Diag_D$. 
The diagram is a finite forest in which with each vertex
we have associated, in particular, a completion. Let $depth_D$ be the depth of the diagram. By the global bounds on
the number of rigid and strictly solid families of values of rigid and solid limit groups, and since there are only finitely many
(graded) completions in the initial level of the diagram $Diag_D$, there exists a global bound on the number of fibers
that are associated with a value $q_1$ of the parameter group $<q>$, and the finitely many completions that are
placed in the initial level of $Diag_D$. We denote this bound $init_D$. By applying the same argument, the finiteness of
the completions that are placed in each level of the diagram $Diag_D$, together with the existence of a global
bound on the number of rigid and strictly solid families of values, given a fiber in a completion that is placed in level $m$
of the diagram, there is a global bound on the number of fibers that are associated with the finitely
many completions that are placed in level $m+1$ of the diagram $Diag_D$, and are further associated with the given fiber
(of a completion in level $m$ of the diagram), 
 and with a value of
the parameters $q_{m+1}$, where the global bound does not depend on the level $m$, the given fiber (in level $m$), or the 
value of the parameters $q_{m+1}$. We denote this global bound $width_D$.

By the construction of the diagram $Diag_D$, given the value $q_1$ there are at most $init_D$ fibers that are
associated with it and with completions that are placed in the initial level of $Diag_D$. Since by our assumptions,
$D(p,q_1) \cap D(p,q_2)$
is strictly contained in $D(p,q_1)$, In the set of fibers that are associated with $D(p,q_1) \cap D(p,q_2)$, at least
one of the fibers that are associated with $D(p,q_1)$ is replaced by at most $width_D$ fibers that are associated
with completions that are placed in vertices in the second level of $Diag_D$. Continuing iteratively, since the intersections:
$\cap_{j=1}^m D(p,{q_j})$
$m=1,\ldots,n$, are strictly decreasing, for each index $m$, at least one of the fibers that are associated with 
$\cap_{j=1}^m D(p,{q_j})$, is replaced by at most $width_D$ fibers that are associated with completions in the next level in
the diagram $Diag_D$, in the set of fibers that is associated with the intersection:
$\cap_{j=1}^{m+1} D(p,{q_j})$. In particular, since the diagram $Diag_D$ is finite, if a fiber that is associated with a 
 completion that is placed at a terminal vertex of $Diag_D$ and with
$\cap_{j=1}^m D(p,{q_j})$
 is replaced in 
$\cap_{j=1}^{m+1} D(p,{q_j})$, 
 then such a fiber is replaced by the empty set. 

Therefore, if the intersections:
$\cap_{j=1}^m D(p,{q_j})$, is a strictly decreasing sequence for
$m=1,\ldots,n$, then $n \leq 2 \cdot init_D \cdot (width_D)^{depth_d-1} $,
so the Diophantine family $D(p,q)$ is equational. 

\line{\hss$\qed$}

\vglue 1.5pc
\centerline{\bf{\S3. Duo Limit Groups}}
\medskip
In section 1 we have shown that in the minimal rank case Diophantine sets are
equational, and then used it to show that the sets $NR_s$ ($NS_s$), that indicate those values of the
parameter set $<p,q>$, for which a minimal rank rigid limit group $Rgd(x,p,q,a)$ (solid limit group
$Sld(x,p,q,a)$) admits
at least $s$ rigid (strictly solid families of) specializations, is in the Boolean algebra generated by
 equational sets  (theorems
1.5,1.7).

In the previous section, we have shown that general Diophantine sets are equational.
In the next section, we show that 
sets of the form $NR_s$ and $NS_s$ are stable.
In this section we present the main tool  that we are going to use in proving the stability of the
sets $NR_s$ and $NS_s$ (and afterwards the stability of general definable sets over a free group),
that we call  $duo$ $limit$ $groups$.

In section 4 of [Se3] we defined $configuration$ $limit$ $groups$ that are
associated with rigid and solid limit groups (definition 4.1 in [Se3]). Recall that given a positive integer $s$
and a rigid or solid limit group, $Rgd(x,p,a)$ or $Sld(x,p,a)$, a configuration limit group is obtained as a limit
of a convergent sequence of tuples, $\{x_1(n),\ldots,x_s(n),p_n,a)\}$, where for each index $n$, and every
index $i$, $1 \leq i \leq s$, the values, $(x_i(n),p_n,a)$, are rigid or strictly solid, and for different indices,
$1 \leq i_1 < i_2 \leq s$, the rigid or strictly solid values, $(x_{i_j}(n),p_n,a)$, $j=1,2$,
are distinct or belong to distinct strictly solid
families (see definition 4.1 in [Se3] for the exact definition).

We start by
presenting duo limit groups that are associated with configuration limit
groups of rigid and solid limit groups.
Then we prove the existence of a (universal) finite collection of duo limit groups that are associated with  a
configuration limit group,
that "covers" all the other duo limit groups that are associated with a rigid or a solid
limit group.
We conclude this section by proving a strong uniform bound for that covering property in the rigid case
(theorems 3.3), and leave the analogous statement for solid limit groups as an open question. 
We note, that the strong bound for the covering property is not needed
for proving stability in the sequel.

\vglue 1pc
\proclaim{Definition 3.1} Let $F_k$ be a non-abelian free group, and let
$Rgd(x,p,q,a)$ ($Sld(x,p,q,a)$) be a rigid (solid) limit group
 with respect to the parameter subgroup $<p,q>$. Let $s$ be a (fixed) positive integer, and let
$Conf(x_1,\ldots,x_s,p,q,a)$ be a configuration limit group that is associated with the
limit group $Rgd(x,p,q,a)$ ($Sld(x,p,q,a)$) (see definition 4.1 in [Se3] for configuration limit groups).

A $duo$ $limit$ $group$, $Duo(d_1,p,d_2,q,d_0,a)$ (shortened as $Duo$),
is a limit group that is obtained as an amalgamated free product
of two completions along the common distinguished vertex group in the abelian decompositions that are
 associated with
their terminal levels, and such that the amalgamated free product has  the following properties:
\roster
\item"{(1)}" with $Duo$ there exists an associated map:
$$\eta: Conf(x_1,\ldots,x_s,p,q,a) \to Duo.$$
For brevity, we denote $\eta(p), \eta(q), \eta(a)$ by $p,q,a$ in correspondence.

\item"{(2)}"
$\eta(F_k)=\eta(<a>) \, < \, <d_0>$, $\eta(<p>) < <d_1>$,
and $\eta(<q>) < <d_2>$.

\item"{(3)}" $Duo=Comp_1(d_1,p,a)*_{<d_0>} \, Comp_2(d_2,q,a)$, where $Comp_1(d_1,p,a)=<d_1>$
and $Comp_2(d_2,q,a)=<d_2>$, are
(graded) completions
with respect to the parameter subgroup $<d_0>$,
and the subgroup $<d_0>$ is the distinguished vertex group in the (two) abelian decompositions that are associated with the
terminal levels of the two completions.

\item"{(4)}" there exists a tuple of values, $(x_1^0,\ldots,x_s^0,p_0,q_0,a)$, which is a
specialization of  the configuration
limit group $Conf$, for which:

\itemitem{(i)} the corresponding  values, $(x_i^0,p_0,q_0,a)$, $i=1,\ldots,s$, are distinct
and rigid specializations of the rigid limit group, $Rgd(x,p,q,a)$ (strictly solid and belong
to distinct strictly solid families of $Sld(x,p,q,a)$).

\itemitem{(ii)} the
value, $(x_1^0,\ldots,x_s^0,p_0,q_0,a)$,
can be extended
to a specialization of the duo limit group $Duo$ (i.e., there exists a configuration
homomorphism that can be extended to a specialization of $Duo$).
\endroster
\endproclaim

With a duo limit group we naturally associate their $duo$-$families$, which are the (duo) analogue of a fiber
of a completion.

\vglue 1pc
\proclaim{Definition 3.2}
Let $Duo(d_1,p,d_2,q,d_0,a)$ be a duo limit group, so that
$Duo=Comp_1(d_1,p,a)*_{<d_0>} \, Comp_2(d_2,q,a)$.
We call a set of specializations of $Duo$  a $rectangle$, if there exists some value $d_0^0$ of the
variables $d_0$, and a fiber of the completion $Comp_1$ and a fiber of the completion $Comp_2$, that
are both associated with the value $d_0^0$, such that
the specializations in the rectangle are precisely all the specializations of $Duo$ that restrict
to values in the fibers of $Comp_1$ and $Comp_2$.

A sequence of specializations of the duo limit group $Duo$ is called a duo test sequence,
if it restricts to test sequences of
the completions, $Comp_1$ and $Comp_2$.
We say that a finite collection of duo limit groups, $Duo_1,\ldots,Duo_t$, covers a rectangle $rectangle$,
that is associated with some duo
limit group $Duo$,  if there exists
a finite collection of rectangles that are associated with the duo limit groups, $Duo_1,\ldots,Duo_t$,
such that every duo test sequence in the  given rectangle $rectangle$,
has a  subsequence, that  restricts to
a sequence of configuration homomorphisms (i.e., values that satisfy condition (i) in part (4) in definition 3.1),  
and the values of these configuration homomorphisms can
be extended to  values in one of the  rectangles from the (fixed) finite collection of
rectangles that are associated with $Duo_1,\ldots,Duo_t$.
\endproclaim

The procedure that was used to prove the equationality of Diophantine sets in the previous section,
enables one to prove the existence of a finite collection of duo limit groups, that cover
all the rectangles that are associated with a duo limit group that is associated with a given rigid or a solid
limit group. We note that the strong boundedness that is proved in theorem 3.3 only in the rigid case,
is not used in
proving stability in the sequel. However, the main diagram that is constructed in order to prove
the theorem, and its associated duo limit groups  (that generalize to solid limit groups as well), 
 are the main tools in our approach to  stability.

\vglue 1pc
\proclaim{Theorem 3.3} Let $F_k$ be a non-abelian free group, let $s$ be a positive
integer, and let
$Rgd(x,p,q,a)$  be a rigid  limit group
defined over $F_k$.
There exists a finite collection of duo limit groups that are associated
with configuration homomorphisms of $s$ distinct rigid homomorphisms of $Rgd$, $Duo_1,\ldots,Duo_t$,
and some global bound $b$,
so that every rectangle that is associated with a duo limit group $Duo$, that is associated
with configuration homomorphisms of $s$ distinct rigid homomorphisms of $Rgd$,
is covered by the given finite collection $Duo_1,\ldots,Duo_t$.
Furthermore, every rectangle that is associated with an arbitrary duo limit group $Duo$, is covered by
at most $b$ rectangles that are associated with the given finite collection,
$Duo_1,\ldots,Duo_t$.
\endproclaim

\nfp To
construct the (finite) universal collection of duo limit groups,
we apply the iterative procedure that was used  to prove the equationality
of Diophantine sets (theorem 2.1).

First, we associate with the given rigid limit group $Rgd(x,p,q,a)$  and
the given positive integer $s$, a finite collection of configuration limit groups (as we
did in section 4 of [Se3]). To do that we collect all the tuples of the
form $(x_1^0,\ldots,x_s^0,p_0,q_0,a)$, for which each value, $(x_i^0,p_0,q_0,a)$, is a rigid specialization
of the rigid limit group $Rgd(x,p,q,a)$, with respect to
the parameter subgroup $<p,q>$, and so that for each $i,j$, $1 \leq i < j \leq s$, the rigid
values,  $(x_i^0,p_0,q_0,a)$ and $(x_j^0,p_0,q_0,a)$ are distinct.
By the standard arguments that are presented in section 5 of [Se1], with this
collection of tuples, $\{(x_1^0,\ldots,x_s^0,p_0,q_0,a)\}$, we can canonically associate a finite
collection of maximal (configuration) limit groups, $Conf_i(x_1,\ldots,x_s,p,q,a)$, $1 \leq i \leq m$.

With each of the configuration limit groups, $Conf_i(x_1,\ldots,x_s,p,q,a)$, viewed as a graded limit
group with respect to the parameter subgroup $<q>$, we associate its
taut graded Makanin-Razborov
diagram (see proposition 2.5 in [Se4] for the construction of the taut Makanin-Razborov diagram).
With each resolution in the taut Makanin-Razborov diagram, we further associate
its singular locus, and the graded resolutions that are associated with
each of the strata in the singular locus. We conclude the first step of the
construction of the diagram, by associating the (graded) completion  with each of the
graded resolutions in our finite collection, that we denote, $Comp(z,p,q,a)$. Note that each of the constructed completions is
graded with respect to the parameter subgroup, $<q>$.

We continue to the construction of the second step of the diagram with each of the completions
$Comp(z,p,q,a)$ in parallel (note that the elements $x_1,\ldots,x_s \in Comp(z,p,q,a)$ can be expressed as words in the
generators $z$ of the completion $Comp(z,p,q,a)$). With each such completion we associate the collection of
tuples of values, $(y_1^0,\ldots,y_s^0,z_0,p_0,q_1^0,q_2^0,a)$, for which:
\roster
\item"{(1)}" $(z_0,p_0,q_1^0,a)$ factors through the completion, $Comp(z,p,q,a)$.
Each of the associated values (the restrictions), $(x_i^0,p_0,q_1^0,a)$, $1 \leq i \leq s$,
is a rigid specialization of
the given rigid
limit group $Rgd(x,p,q,a)$, and any two rigid specializations,
$(x_i^0,p_0,q_0,a)$ and $(x_j^0,p_0,q_0,a)$, are distinct for
$1 \leq i < j \leq s$.

\item"{(2)}" each of the values, $(y_i^0,p_0,q_2^0,a)$, $1 \leq i \leq s$,
is a rigid specialization of the rigid limit group
$Rgd(x,p,q,a)$.
Any two rigid specializations,
$(y_i^0,p_0,q_2^0,a)$ and $(y_j^0,p_0,q_2^0,a)$, are distinct  for
$1 \leq i < j \leq s$.
\endroster

With the completion, $Comp(z,p,q,a)$, we associate the collection of all the
sequences: $$\{(y_1(n),\ldots,y_s(n),z(n),p(n),q_1(n),q_2(n),a)\}_{n=1}^{\infty}$$
so that for each $n$, the corresponding tuple of values satisfies  conditions (1) and (2), and
the (restricted) sequence $\{(z(n),p(n),q_1(n),a)\}_{n=1}^{\infty}$ form a (graded) test sequence with respect to
the given (graded) completion $Comp(z,p,q,a)$. By the techniques that were used to analyze graded formal limit groups,
that are presented in section 3 of [Se2], with this collection of sequences it is possible to canonically associate a finite
collection of (graded) limit groups, that have the structure of closures of the completion, $Comp(z,p,q,a)$ (i.e., they
differ from the completion, $Comp(z,p,q,a)$, in additional roots that are possibly added to abelian vertex groups in the abelian
decompositions that are associated with the various levels of the completion, $Comp(z,p,q,a)$, and they also differ from the completion
in the limit group that is associated with their terminal level).
However, note
that the constructed limit groups are graded with respect to the parameter subgroup, $<q_1,q_2>$, and not with respect to
the parameter subgroup $<q>=<q_1>$ like the original completion, $Comp(z,p,q,a)$.
We will denote these limit groups, that we view and call graded
closures, $DQCl_i(s,z,p,q_1,q_2,a)$.

We continue by looking at all the tuples of values, $(y_1^0,\ldots,y_s^0,z_0,p_0,q_1^0,q_2^0,a)$, that satisfy
conditions (1) and (2), and do not factor through any of the (finite) closures,
$DQCl_i(s,z,p,q_1,q_2,a)$.
With this collection of tuples we can
associate a canonical finite collection of maximal limit groups, $M_j(y_1,\ldots,y_s,z,p,q_1,q_2)$, which we
view as graded limit groups with respect to the parameter subgroup $<q_1,q_2>$.

Using the construction of quotient resolutions, that is used in the general step of the sieve procedure [Se6],
we associate with
this collection of tuples of values, and with the finitely many graded limit groups, $M_j$, that are associated with their
Zariski closure,  finitely many multi-graded resolutions,
and with each such multi-graded resolution we
associate its (multi-graded) core resolutions, developing resolutions, anvils, and
(possibly) sculpted resolutions and  carriers (see [Se6] for a detailed description of the
iterative construction of these multi-graded resolutions and the finite collection of resolutions that are attached to them).

We continue iteratively, precisely as we did in proving the equationality
of Diophantine sets (theorem 2.1). We start each step with the completions that were constructed in the previous step, and
continue with each of them in parallel. We first look at all the test sequences of such a completion that can be extended to tuples
of values that satisfy the properties (1) and (2) above. With these collections of test sequences we associate finitely many
closures of the completions that were constructed in the previous step of the procedure. Then we consider all the specializations of the completions
that were constructed in the previous step of the procedure, that can be extended to tuples of values that satisfy properties (1) and (2),
and these tuples of values do not factor through any of the previously constructed closures (of the completions that were
constructed in the previous step). We analyze these tuples of values by applying the construction of quotient resolutions,
that was used in the general step of the sieve procedure for quantifier elimination [Se6]. This analysis associated with the given
collection of tuples of values finitely many multi-graded resolutions,  together with their core resolutions, anvils, developing
resolutions, and possibly sculpted resolutions and their carriers (all these are presented in detail in [Se6]).

Finally, like the
sieve procedure for quantifier elimination [Se6], and like the iterative procedure that was used in proving the equationality
of Diophantine sets in the previous section, the iterative procedure that we described terminates after finitely many steps.

\vglue 1pc
\proclaim{Proposition 3.4} The iterative procedure that is presented above terminates after finitely
many steps.
\endproclaim

\nfp Identical to the proof of proposition 2.2.

\line{\hss$\qed$}

When the iterative procedure terminates, we obtain a finite diagram that we denote, $Diag$. In each vertex
 of the diagram there is a completion. The completions that are placed in vertices in the initial level of the
diagram, $Diag$, are the completions of the resolutions in the graded taut Makanin-Razborov diagrams of the
maximal configuration limit groups that are associated with the given rigid limit group, $Rgd(x,p,q,a)$. Note
that the resolutions and their completions in the initial level are graded with respect to the parameter subgroup,
$<q_1>$.

The completions that are placed in vertices in the second level of the diagram, are either closures $DQCl_i$ in which proper 
roots were added to abelian vertex groups in the completion, $Comp(z,p,q,a)$, that they were constructed from, or  completions of the developing
resolutions of anvils that were constructed in the second step of the iterative procedure. These closures and developing
resolutions and their completions
are graded with respect to the parameter subgroup, $<q_1,q_2>$. Each completion in the initial
level of the diagram is connected by finitely many (possibly no)  directed edges to the closures and the completions of developing
resolutions of the anvils that were constructed from it in the second step of the iterative procedure.

The
completions that are placed in vertices in the next levels of the diagram are similar. The completions that
are placed in vertices in level $m$ of the diagram are either closures of completions in level $m-1$ in which proper roots were
added to abelian vertex groups of these completions (from level $m-1$), or the completions of developing resolutions of anvils that
were constructed in step $m$ of the iterative procedure. These closures, developing resolutions and their completions
are graded with respect to the parameter subgroup, $<q_1,\ldots,q_m>$.
A completion of a developing resolution in level $m-1$
is connected by finitely many (possibly no) directed edges to its closures and to the completions of developing resolutions of the
anvils that were constructed from it in level $m$ of the iterative procedure.

\medskip
To define the universal set of duo limit groups, that are claimed in theorem 3.3, we start
with the collection of completions that were constructed along the terminating iterative
procedure, and with each such completion we associate a finite collection  of duo
limit groups.

Given a (graded) completion that was constructed along the diagram $Diag$,
that we denote $Comp$,
we associate with it finitely
many duo limit groups.
To construct these duo limit groups, we fix a generating set of each vertex group in each of the abelian
decompositions that are associated with the various levels of the completion, $Comp$, and a generating
set of the parameter subgroup, $<q>$.
We look at the entire collection
of graded test sequences that factor through the given graded completion, $Comp$, for which the restrictions
of the values in these test sequences to the variables $p$, can be extended to configuration
homomorphisms of at least one of the (finitely many) maximal configuration limit groups that are associated with
the given rigid limit group, $Rgd(x,p,q,a)$.

We further
require that the $n$-th value in each of these test sequences, and its  extension to a configuration
homomorphism, will satisfy that the maximal length of the (restricted)
values of the fixed generating sets of each of
the non-distinguished vertex groups in the completion $Comp$, are at least $n$ times bigger than the
maximal length of the (restricted) values of a fixed generating set of the parameter subgroup, $<q>$.

With this entire collection of graded test sequences, and their extensions to configuration homomorphisms,
we associate a graded Makanin-Razborov diagram, precisely as we did in constructing
the formal graded Makanin-Razborov diagram in section 3 of [Se2].
By the construction of formal
graded Makanin-Razborov diagrams, the abelian decompositions that are associated with the various limit groups
that appear along the resolutions of the diagrams, are the graded abelian decompositions of these limit groups where
the parameter subgroup is taken to be the completion $Comp$, from the diagram $Diag$,
that we have started the construction with. Furthermore, by the analysis of graded formal resolutions, as it appears in
section 3 of [Se2], each of the resolutions in the constructed Makanin-Razborov diagrams terminates with a (graded) closure
of the graded completion, $Comp$, that we have started with.

The sequences of values that we analyze, are values from test sequences of the completion, $Comp$, together with extensions
of the values of the subgroup, $<p,q>$, to (non-degenerate) configuration homomorphisms of one of the finitely
many configuration limit
groups that are associated with the rigid limit group, $Rgd(x,p,q,a)$. Hence, each value in these sequences is obtained
from a value of the completion, $Comp$, a value of the parameters $q$, and $s$ rigid values of the rigid limit
group, $Rgd(x,p,q,a)$.

\noindent
We further required that the lengths
of the values of the variables $q$ is much smaller than the lengths of the values of the fixed generating
sets of the vertices in the abelian decompositions that are associated with the various levels of
the completion, $Comp$.

At this point we analyze the algebraic structure of a limit group that is obtained as a limit of a sequence of values that
we consider. We do that by looking at the limit tree to which a subsequence of such a sequence converges.
Since each value in the sequences we analyze  is obtained from a value of the completion, $Comp$, a value of the
parameters $q$, that is much shorter than the values of fixed generating sets of  the vertex groups in the abelian
decompositions that are associated with the various levels of $Comp$, except for the terminal level, by adjoining $s$ rigid
values of $Rgd(x,p,q,a)$, and the values of the completion $Comp$ form a test sequence of it, the abelian decomposition
of the obtained limit group that can be read from the limit tree, must have similar structure as the abelian decomposition that is
associated with the top level of the completion, $Comp$.

By going down through the levels of the completion $Comp$, the same argument implies that the obtained limit group is the
amalgamation of a quotient of the completion, $Comp$, with a $slow$ limit group that contains the subgroup $<q>$,
that are amalgamated along
 a quotient of the terminal level of
the completion $Comp$. Hence, the (formal) abelian JSJ decomposition of this obtained limit group (that is an
abelian JSJ decomposition with respect to the image of $Comp$) can be constructed
from a
graded abelian decomposition, $\Delta$,
of this $slow$ limit group with respect to a parameter subgroup which is the terminal level of the
completion, $Comp$, where the distinguished vertex in the graded abelian decomposition of $slow$, $\Delta$, is amalgamated with the
image of the completion, $Comp$, along the image of the terminal level of $Comp$, that is contained  in the distinguished vertex
of the abelian JSJ decomposition of the subgroup $slow$, $\Delta$.

By going through the levels of each of the resolutions in the (formal) Makanin-Razborov diagrams that we have
associated with the finitely many limit groups that are associated with the sequences of values that we consider,
the (formal) abelian JSJ decompositions that are associated with the limit groups that are placed along these
resolutions have a similar structure, i.e., they are obtained from graded JSJ decompositions of the corresponding
$slow$ subgroups with respect to the image of the terminal level of the completion, $Comp$, where the distinguished vertex
group
in each such abelian decomposition is replaced by a limit group that is obtained from it by an amalgamated
product with the image of the completion $Comp$, along an amalgamated subgroup which is the image of the
terminal level of $Comp$.

Therefore, the completion
of a resolution in the constructed Makanin-Razborov diagrams is the amalgamated
product of a graded closure of the completion, $Comp$, with another completion (that contains the subgroup
$<q>$ as a subgroup), that are amalgamated along the common distinguished vertex groups in the abelian
decompositions that are associated with the terminal levels of the two completions.



By the construction of the completions in these (formal) graded Makanin-Razborov diagrams,
there is also a natural map from a
(maximal) configuration limit group of the original rigid limit group, $Rgd(x,p,q,a)$, into it. The subgroup $<p>$
is mapped into the closure of the given completion $Comp$,
and the subgroup $<q>$
is mapped into the other completion. Hence, the obtained amalgamated product is a duo limit group.
We take the completions of the resolutions that appear in the entire finite collection
of Makanin-Razborov diagrams that are associated with the
various completions, $Comp$, that are placed in the various vertices of the diagram $Diag$,
to be the finite collection of (universal) duo limit groups, $Duo_1,\ldots,Duo_t$,
that is indicated in the statement of theorem 3.3.

\medskip
Let $Duo$ be a duo limit group that is associated with the given rigid  limit group, and
suppose that we are given a rectangle, $rectangle$, that factors through it, i.e.,
a rectangle that is associated with a given value, $d_0^0$, of the variables
$d_0$ in the duo limit group $Duo$. We need to show that the given rectangle, $rectangle$, that factors through the duo
limit groups $Duo$, is covered by a bounded collection of rectangles that factor through the
(universal) finite collection of duo limit
groups $Duo_1,\ldots,Duo_t$. Note that the bound on the number of rectangles of $Duo_1,\ldots,Duo_t$ is
supposed to be global and does not depend on the particular duo limit group $Duo$, or the rectangle $rectangle$.

By definition, the duo limit group, $Duo$, contains an image of a configuration limit group of the rigid limit
group, $Rgd(x,p,q,a)$. We denote this configuration limit group, $Conf$.
Given the rectangle, $rectangle$, we start with a value $q_1$ of the parameters $q$,
 that can be extended to  a value in  the rectangle, $rectangle$, that restricts to a
$configuration$ $homomorphism$ of the  configuration limit group $Conf$.
I.e.,
a value in the rectangle that satisfies property (4)
of a duo limit group (definition 3.1). With this value $q_1$ of the parameters $q$ we associate
the boundedly many fibers that are associated with it in the initial level of the diagram, $Diag$, that was
constructed iteratively from the sets of configuration homomorphisms.
We further
associate with the value $q_1$ the boundedly many rectangles of those duo limit groups, $Duo_1,\ldots,Duo_t$,
that were constructed from completions that appear in the initial level of the diagram $Diag$. Note that
by the construction of the diagram $Diag$, and the duo limit groups, $Duo_1,\ldots,Duo_t$,  there
is a global bound on the  number of
rectangles that are associated with any given value of the parameters $q$, and with those duo limit groups,
$Duo_1,\ldots,Duo_t$, that are associated with completions in the initial level of $Diag$.

We consider all the duo test sequences that factor through the given rectangle, $rectangle$,
and restrict to configuration homomorphisms of the configuration limit group $Conf$ (that is mapped into
the Duo limit group $Duo$ that covers the given rectangle $rectangle$).
I.e., duo test sequences of values in the rectangle $rectangle$,
that satisfy property (4) in definition 3.1.

Given this collection of duo test sequences of the given rectangle $rectangle$,
we look at those duo test sequences for which
their restrictions to configuration homomorphisms ([Se3],definition 4.1) can be extended to values of one
of the (boundedly many) rectangles that are associated with the value $q_1$ of the parameters $q$,
and with those duo limit groups,
$Duo_1,\ldots,Duo_t$, that were constructed from completions that are placed in the initial level of the
diagram, $Diag$.

For each such duo test sequence,
we extended the restrictions of the values in the duo test sequence
to configuration homomorphisms, to the shortest possible value in the (boundedly many) rectangles,
that are associated with $q_1$ and with the duo limit groups $Duo_1,\ldots,Duo_t$, that were constructed
from completions in the initial level of the diagram $Diag$.

By the techniques that are presented in section 3 of [Se2] (that constructs graded formal limit groups), with this
collection of duo test sequences and their extended values,
we can associate finitely many limit groups, that are all $duo$ $closures$ of the given
duo limit group that is dual to (i.e., the coordinate group of) the rectangle, $rectangle$, i.e.,
limit groups that are amalgamated products of closures of the two completions,
$Comp_1(d_1,p,a)$ and $Comp_2(d_2,q,a)$, that are associated with the rectangle, $rectangle$.

Furthermore, by the construction of the
duo closures, with
each such duo closure, there is an associated map from the limit group which is the dual to
one of the boundedly many rectangles that are
associated with the duo limit groups, $Duo_1,\ldots,Duo_t$, and with the value $q_1$, into the duo closure.

\noindent
Note that it can be that no sequence  of restrictions of duo test sequences in the rectangle,  $rectangle$,
to configuration homomorphisms, can be extended to values in the rectangles that are associated with $q_1$ and with the duo limit
groups, $Duo_1,\ldots,Duo_t$, that are associated with  the initial level of the diagram $Diag$. In this (empty) case, no duo closures
are associated with the rectangles that are associated with $q_1$ and with the duo limit groups, $Duo_1,\ldots,Duo_t$, that are associated
with the initial level of the diagram $Diag$.

We have associated finitely many (possibly none) duo closure with the given rectangle, $rectangle$. With each duo closure of
rectangle, there is a pair of associated closures of the two completions, $Comp_1(d_1,p,a)$ and $Comp_2(d_2,q,a)$. By
definitions 1.14 and 1.15 in [Se2], with a closure of a completion one naturally associates  with each abelian vertex group of
an abelian decomposition that is associated with one of the levels of the corresponding completion,
a coset of a finite index subgroup. Hence with each duo closure of $rectangle$, we associate
a coset of a finite index subgroup with each abelian vertex group  that is associated with one of the levels of $Comp_1$ and
$Comp_2$.

Since there are finitely many duo closures of $rectangle$, for each abelian vertex group that is associated with a level of
$Comp_1$ or $Comp_2$, we can take the intersection of the finitely many finite index subgroups that are associated with it. Hence,
with each abelian vertex group that is associated with one of the levels of $Comp_1$ or $Comp_2$ we associate a finite index
subgroup, and with each duo closure of the given $rectangle$ we can associate a finite set of collections of cosets of each
of these finite index subgroups.

We can place the (finite) set of all possible collections of cosets of the finite index subgroups that are associated with the
abelian
vertex groups in $Comp_1$ and $Comp_2$ in a planar diagram, where one axis is for collections of cosets of the finite index subgroups of
abelian vertex groups in $Comp_1$, and the second axis is for collections of cosets of finite index subgroups of abelian
vertex groups in $Comp_2$. The given set of duo closures of the given rectangle, $rectangle$, cover some (possibly none)
of the possible collections of cosets. To prove theorem 3.3, we show that even though the number of duo closures of
the given rectangle is finite and not necessarily bounded, and the indices of the finite index subgroups need not be bounded
either, it is possible to get a combinatorial bound on the form of the collections of cosets that are associated with
the duo closures that we constructed.

\vglue 1pc
\proclaim{Proposition 3.5} After possibly replacing the set of closures $\{cld_i\}$ and their associated
maps, $\{\eta_i\}$, and hence possibly changing the planar diagram that is associated with the set
of closures (as we may need to refine the collections of cosets of finite index abelian subgroups that
need to be considered), the points in the finite planar diagram that are associated with (the new)
collections of cosets of the finite index subgroups in the planar diagram that are associated with the finitely
many (new) duo closures,
that we constructed from boundedly many rectangles of the duo limit groups, $Duo_1,\ldots,Duo_t$, and from the given rectangle,
$rectangle$, are the union of boundedly many
product domains, where each such product domain is determined by a subset of rows and columns of the finite planar diagram.
Furthermore, the bound on the number of product domains depend only on the (universal) duo limit groups,
$Duo1,\ldots,Duo_t$.
\endproclaim

\nfp With the given rectangle, $rectangle$, and the boundedly many rectangles of $Duo_1,\ldots,Duo_t$, we have
associated finitely many (possibly none) duo closures of the duo limit group which is dual to $rectangle$,
where into each such closure there is a map from one of the rectangles of $Duo_1,\ldots,Duo_t$. The points in
the planar diagram are associated with these duo closures (with each closure we have associated finitely
many points in the planar diagram).

We fix one of the boundedly many chosen rectangles of $Duo_1,\ldots,Duo_t$, and denote it $Rectangle$. We denote
the duo limit group that is dual to $Rectangle$, $duo_R$, and the duo limit group that is dual to the
rectangle that we have started with, $rectangle$, we denote $duo_r$. With the rectangles, $Rectangle$ and
$rectangle$, we have associated finitely many  (possibly none) duo
closures of $duo_r$, that we denote $cld_1,\ldots,cld_m$,
and maps: $\eta_i: duo_R \to cld_i$, $i=1,\ldots,m$.

In order to prove the proposition our goal is to show that the points in the planar diagram that are associated with
the closures, $cld_1,\ldots,cld_m$, are a bounded union of product domains, where the bound on the number
of product domains depend only on the (universal) duo limit groups, $Duo_1,\ldots,Duo_t$, and not on the
given duo limit group, $Duo$, or its rectangle,  $rectangle$. Since we have chosen only boundedly many
(possibly none) rectangles of $Duo_1,\ldots,Duo_t$, a presentation of the points in the planar diagram that
are associated with one of these rectangles, $Rectangle$, as a bounded union of product domains, clearly
implies the statement of the proposition.

The duo limit groups that are dual to $rectangle$ and $Rectangle$ can be represented as amalgamated products over
the coefficient group $<a>=F_k$: $duo_r=Comp_1(d_1,p,a)*_{<a>} \, Comp_2(d_2,q,a)$ and
$duo_R=RComp_1(u_1,p,a)*_{<a>} \, RComp_2(u_2,q,a)$, and so is each of the closures of $duo_r$:
$cld_i=CComp^i_1(d^i_1,p,a)*_{<a>} \, CComp^i_2(d^i_2,q,a)$, $i=1,\ldots,m$. The maps $\eta_i:duo_R \to cld_i$
map the image of the configuration limit group in $duo_R$ onto the image of the configuration limit
group in $cld_i$. Hence, in particular, it maps the subgroups, $<p>$ and $<q>$ in $duo_R$, onto the
corresponding subgroups $<p>$ and $<q>$ in $cld_i$. However, it may be that $RComp_1$ is not mapped into
$CComp^i_1$ or $RComp_2$ is not mapped into $CComp^i_2$. To prove the proposition, we first replace the given set
of closures, $\{cld_i\}$, and their associated maps, $\{\eta_i\}$, by a different collection of closures and maps,
so that for the new set of closures, $RComp_1$ and $RComp_2$ are mapped into $CComp^i_1$ and $CComp^i_2$ in
correspondence.

\vglue 1pc
\proclaim{Lemma 3.6} It is possible to replace the given set of closures of $duo_r$, by a new (finite) set of
closures (still denoted $cld_i$), that cover the same collections of cosets of finite index subgroups of
the abelian vertex groups that appear in the various levels of the completions $Comp_1$ and $Comp_2$ as
the previous set of closures, so that for every new map $\eta_i:duo_R \to cld_i$, $RComp_1$ is mapped into
$CComp^i_1$ and $RComp_2$ is mapped into $CComp^i_2$.
\endproclaim

\nfp
The closure $cld_i$ can be written as an amalgamated product: ${CComp^i_1}_{<a>} \, CComp^i_2$.
Suppose that for one of the closures, $cld_i$, $i=1,\ldots,m$, either the completion $RComp_1$ is not mapped
by $\eta_i$ into $CComp^i_1$ or $RComp_2$ is not mapped by $\eta_i$ into
$CComp^i_2$.  Wlog we can assume that the image of $RComp_1$ is not in $CComp^i_1$.

As the subgroup $<p,a>$ is contained in $CComp^i_1$, and 
the image of $RComp_1$ is not in $CComp^i_1$, the image of $RComp_1$ in $cld_i$, $IRC_1$, inherits a non-trivial
graph of groups from the presentation of $cld_i$ as an amalgamated product: $cld_i={CComp^i_1}*_{<a>} \, CComp^i_2$.
By going through the various levels of the completion $CComp^i_2$ from top to bottom, there is a highest level of
$CComp^i_2$, that we denote level $h$, for which the inherited graph of groups $IRC_1$ is non-trivial. We denote
this inherited abelian decomposition $\Delta$. Since the 
decomposition that is associated with every level of the completion $CComp^i_2$ is an abelian decomposition, the 
graph of groups $\Delta$ is an abelian graph of groups of $IRC_1$.

The abelian graph of groups $\Delta$ of $IRC_1$ naturally extends to an
abelian graph of groups $\Delta'$ of the amalgamation of $IRC_1$ 
with the completion $RComp_2$ along the amalgamated (coefficient) subgroup $<a>=F_k$. By construction the
map $\eta_i:duo_R \to cld_i$ factors through that amalgamated subgroup. 

\noindent
Since the subgroup $<p,a>$ is contained in the distinguished vertex group of the abelian decomposition $\Delta$, 
and the subgroup $<q,a>$ is contained in $RComp_2$, the subgroup $<p,q,a>$ is contained in the distinguished
vertex group of the abelian decomposition $\Delta'$. The image of the configuration subgroup $Conf$ in the
amalgamation of $IRC_1$ and $RComp_2$, $<x_1,\ldots,x_s,p,q,a>$, is generated by the subgroup $<p,q,a>$ and
the elements $x_1,\ldots,x_s$, where each of the subgroups $<x_j,p,q,a>$ is rigid with respect to the parameter
subgroup $<p,q,a>$. Since the (parameter) subgroup $<p,q,a>$ is elliptic in $\Delta'$, and the subgroups 
$<x_j,p,q,a>$ are rigid, the 
abelian decomposition
that is inherited by the subgroup, $<x_1,\ldots,x_s,p,q,a>$, from $\Delta'$ has to be trivial, and so the
entire image of the configuration limit group $Conf$ is contained in the distinguished vertex group in
$\Delta'$.

With every test sequence of the closure, $cld_i$, we can associate a sequence of homomorphisms of $duo_R$ into the coefficient group (by precomposing
homomorphisms of $cld_i$ with the map $\eta_i$). We consider all the test sequences of the closure $cld_i$ and use the modular groups that are associated 
with the abelian decomposition, $\Delta'$, that act trivially on the image of the configuration limit
group $Conf$, to shorten the restrictions of these homomorphisms to homomorphisms of the completion, $RComp_1$. 
By the construction of formal limit groups (section 3 in [Se2]),
with the collection of all these (shortened) sequences we can associate a finite collection of closures of $cld_i$ (that are closures of $duo_r$),
that form a covering closure of $cld_i$ (see definition 1.16 in [Se2] for a covering closure), so that (by shortening the restrictions of the
homomorphisms of $duo_R$ to homomorphisms of $RComp_1$) the image of the completion $RComp_1$ in each of these closures inherit a trivial decomposition
from each of the abelian decompositions that is associated with the top $h$ levels of $CComp^{i'}_2$. 

By repeating this argument iteratively, we can replace the closure $cld_i$ by finitely many closures of it, so that the image of $RComp_1$ in each
of these closures is contained in $CComp^{i'}_1$. An identical argument proves the same for the image of $RComp_2$, and the lemma follows.  

\line{\hss$\qed$}

In the sequel we continue with the new set of closures, still denoted $\{cld_i\}$, $i=1,\ldots,m$, with the
properties that are claimed in lemma 3.6. Note that this new set of closures covers the same collections of
cosets of finite index subgroups of abelian vertex groups in the various levels of $Comp_1$ and $Comp_2$.
However, by replacing the set of closures, we may need to replace the finite diagram of collection of such
cosets, and in the sequel we continue with this new diagram.

In the duo limit groups, $duo_R$ (that is dual to the rectangle, $Rectangle$), there is an image of the
configuration limit group, $Conf(x_1,\ldots,x_s,p,q,a)$. We denote this image: $<x_1,\ldots,x_s,p,q,a>$.
In the rectangle, $Rectangle$, there is a value that restricts to a (non-degenerate)
value of the configuration limit group,
$Conf$, i.e., a value for which the values of the $x_j$'s are rigid and distinct values of the given rigid
limit group $Rgd(x,p,q,a)$ (rigid with respect to the parameter subgroup $<p,q>$).

Each of the elements $x_j \in duo_R$, can be written in a normal form with respect to the amalgamated
product: $duo_R={RComp_1}*_{<a>} \, RComp_2$. For each $j$ let this normal form be:
$x_j=f^j_1s^j_1 \ldots f^j_{r_j}s^j_{r_j}$, where $f^j_{\ell} \in RComp_1$ and $s^j_{\ell} \in RComp_2$, $\ell=1,\ldots,r_j$.
We now look at the completion $RComp_1(u_1,p,a)$ as a graded limit group with respect to the parameter subgroup
$<p,a>$, and at the completion $RComp_2(u_2,q,a)$ as a graded limit group with respect to the parameter subgroup
$<q,a>$. With each of these graded limit groups we associate its graded Makanin-Razborov diagram. Clearly, all
the values of the two completions, $RComp_1$ and $RComp_2$, factor through these two completions.

Suppose that $(x_1^0,\ldots,x_s^0,p_0,q_0,a)$ is a (non-degenerate) configuration homomorphism that extends to
a value in the rectangle, $Rectangle$, i.e., it extends to a value that factors through the duo limit group $duo_R$.
This extended value restricts to values of the two completions,
$(u_1^0,p_0,a)$ and $(u_2^0,q_0,a)$, and to values of the
elements $f^j_{\ell} \in Rcomp_1$ and $s^j_{\ell} \in RComp_2$, $j=1,\ldots,s$,
$\ell=1,\ldots,r_j$.

The values $(u_1^0,p_0,a)$ and $(u_2^0,q_0,a)$ factor through graded resolutions in the graded Makanin-Razborov diagrams
of $RComp_1$ and $RComp_2$ in correspondence. Since the values, $x_1^0,\ldots,x_s^0$, are rigid values of
$Rgd(x,p,q,a)$ (with respect to the parameter subgroup $<p,a>$), the elements,
$f^j_{\ell} \in Rcomp_1$ and $s^j_{\ell} \in RComp_2$, must belong to the distinguished vertex groups
(the vertex groups that
contain the parameter subgroups $<p,a>$ and $<q,a>$ in correspondence)
in all the abelian decompositions along the various levels of the two graded resolutions of
$Rcomp_1$ and $RComp_2$ (they must belong to the distinguished vertex groups, since otherwise at least one of the subgroups
$<x_j,p,q,a>$ inherits a non-trivial abelian splitting, a contradiction to the rigidity of the values $x_j^0$).
Therefore, by the bounds on the number of rigid and strictly solid families of rigid
and strictly solid limit groups (theorems 2.5 and 2.9 in [Se3]), for fixed values $p_0$ and $q_0$ of the variables $p$ and $q$,
there is a global bound on the possible values of the variables,
$f^j_{\ell} \in Rcomp_1$ and $s^j_{\ell} \in RComp_2$, that determine (non-degenerate) configuration homomorphisms, i.e., values
$(x_1^0,\ldots,x_s^0,p_0,q_0,a)$ for which the values $x_j^0$, $j=1,\ldots,s$, are rigid and distinct. Furthermore, this
global bound depends only on the (universal) duo limit groups, $Duo_1,\ldots,Duo_t$, and not on the specific rectangle that
is associated with it.

\medskip
We now deduce the conclusion of proposition 3.5 from the universal bounds on the values of the variables $x_j$ for given values
of the variables $p$ and $q$. Let $cld_1,\ldots,cld_{m'}$ be the closures of the duo limit group $duo_r$
that are associated with one of the
boundedly many rectangles, $Rectangle$,  of the duo limit groups:
$Duo_1,\ldots,Duo_t$. We may assume that these closures satisfy the conclusion of lemma 3.5.
With each such duo closure there is an associated collection of cosets of finite index subgroups of the abelian vertex groups that
are associated with the various levels of the completions, $Comp_1$ and $Comp_2$, that are part of the duo limit group $duo_r$,
that is dual to the given rectangle, $rectangle$.

Given the closures, $cld_1,\ldots,cld_{m'}$, we construct a new closure $ucld$ of $duo_r$. We construct $ucld$ to be a closure for
which the finite index subgroups of the abelian groups that are associated with the various levels of $Comp_1$ and $Comp_2$, are
the intersections of the finite index subgroups that are associated with these abelian groups in the set of closures,
$cld_1,\ldots,cld_{m'}$. By construction, each of the closures, $cld_i$, is embedded in $ucld$.

Since the duo limit group  $duo_R$ is mapped by $\eta_i$ into each of the closures, $cld_i$,
the elements $f^j_{\ell}$ and $s^j_{\ell}$ are
mapped by $\eta_i$ into $cld_i$. Since $cld_i$ is mapped into the closure $ucld$,  the elements $f^j_{\ell}$ and $s^j_{\ell}$ are
mapped into $ucld$ via the composition of $\eta_i$ with this embedding, that we denote $\nu_i$. Because there is a global
bound (that depends only on $Duo_1,\ldots,Duo_t$) on the number of distinct values of the elements $f^j_{\ell}$ and $s^j_{\ell}$
for a given  value of the variables $p$ and $q$, there is a global bound (that depends only on $Duo_1,\ldots,Duo_t$) on
the distinct images of the set of elements $f^j_{\ell}$ and $s^j_{\ell}$ under the maps $\nu_i$, $i=1,\ldots,m'$.

We divide the images under the maps $\eta_i$ of the completions, $RComp_1$ and $RComp_2$, into the closures $cld_i$ into boundedly many
equivalence classes, according to the image under the map $\nu_i$  of the subsets of elements $f^j_{\ell}$ and $s^j_{\ell}$
(in correspondence) in
the closure $ucld$.

Suppose that $cld_{i_1}$ and $cld_{i_2}$ are two closures for which the maps of both $Rcomp_1$ and $Rcomp_2$ into $cld_{i_1}$ and
$cld_{i_2}$ belong to the same equivalence classes. Let $CComp_1^{i_1},CComp_2^{i_1},CComp_1^{i_2},CComp_2^{i_2}$ be the completions
that are associated with the two closures, $cld_{i_1}$ and $cld_{i_2}$, in correspondence. Then the groups:
$CComp_1^{i_1}*_{<a>} \, CComp_2^{i_2}$ and
$CComp_1^{i_2}*_{<a>} \, CComp_2^{i_1}$ are also closures of $duo_r$. Furthermore, each of the elements $x_j \in duo_r$
can be represented as:
$x_j=f^j_1s^j_1 \ldots f^j_{r_j}s^j_{r_j}$ in these two closures of $duo_r$. Hence there are values that factor through these
two closures
that restrict to (non-degenerate) configuration homomorphisms.

Therefore, the set of points in the planar diagram that was
associated with the set of closures that satisfy the conclusion of lemma 3.6, and for which the two completions,
$CCom^ip_1$ and $CComp^i_2$, of these closures
belong to the same equivalence classes, form a product domain. Since there are boundedly many equivalence classes of the
completions $CComp^i_1$ and $CComp^i_2$, the collections of points in the diagram that are associated with one of the
closures, 4$ld_1,\ldots,cld_{m'}$, that are all associated with the same rectangle, $Rectangle$, of $Duo_1,\ldots,Duo_t$,
is a bounded union of product domains. Since there is a bound on the number of rectangles that are associated with $Duo_1,\ldots,Duo_t$
that are associated with the initial level of the diagram $Diag$, and with the value $q_1$ of the variables $q$, the set of
points in the diagram that are associated with the entire set of closures, $cld_1,\ldots,cld_m$, is a bounded union of product
domains. Furthermore, the bound depends only on the duo limit groups, $Duo_1,\ldots,Duo_t$, and hence it depends only on the
given rigid limit groups, $Rgd(x,p,q,a)$, that we have started with.

\line{\hss$\qed$}


Suppose that the given (bounded) set of product domains that are associated with the closures of the the rectangle $rectangle$, that
were constructed from boundedly many  rectangles that are associated
with  the duo limit groups, $Duo_1,\ldots,Duo_t$,
 that were constructed in the initial
level of the diagram, $Diag$, and are associated with the value $q_1$, do not cover all the duo test sequences of the given rectangle, $rectangle$.
I.e., there
are still duo test sequences of values in $rectangle$ that do not have subsequences, so that the restrictions of these
subsequences to configuration homomorphisms can not be extended to values that
factor through the (chosen) boundedly many rectangles of $Duo_1,\ldots,Duo_t$.

In this case we look at the planar diagram that has finitely many points, and contain the product domains that appear in the statement
of proposition 3.5. One of the axis of this diagram has finitely many collections of cosets of finite index subgroups of abelian
vertex groups that
appear in the various levels of
the completion, $Comp_1$, and the other axis has collections of cosets of finite index subgroups of abelian vertex groups that appear
in the various levels of $Comp_2$.

\noindent
The boundedly many product domains in the diagram naturally define a stratification of the axis that is associated with $Comp_2$. Two collections
of finite index subgroups of abelian vertex groups in $Comp_2$, are set to be in the same stratum, if they  appear in the projection of the
same product domains. Since there are boundedly many product domains (by proposition 3.5), 
the number of strata in the constructed stratification is bounded (where the bound depends only on the universal duo limit
groups, $Duo_1,\ldots,Duo_t$).

For  each stratum in the stratification of the axis that is associated with $Comp_2$,
 we choose a value $q_2$ of the defining parameters $q$, that extends to a value $d_2^0$ of the variables $d_2$ (i.e., a
specialization of $Comp_2(d_2,q,a)$), with the
following properties, if such a value exists:
\roster
\item"{(i)}" $d_2^0$ belongs to a collection of cosets of finite index subgroups of abelian vertex groups in $Comp_2$ which is in the
specified stratum of the stratification (of the axis that is associated with $Comp_2$).

\item"{(ii)}"  for any product domain that its projection contains the given stratum,
there exists a test sequence of the completion, $Comp(d_1,p,a)$
(which
is the completion  that contain the subgroup $<p>$ in   the given duo limit group $Duo$), for which
the sequence:
$\{(d_1(n),p(n),d_2^0,q_2,d_0^0,a)\}$ is contained in the rectangle, $rectangle$, and restricts to configuration homomorphisms,
 that further extend to values in the closure of $rectangle$ that is associated with that product domain.

\item"{(iii)}" let $\{(d_1(n),p(n),d_0^0,a)\}$ be an arbitrary test sequence of specializations of
$Comp(d_1,p,a)$,
for which
the sequence:
$\{(d_1(n),p(n),d_2^0,q_2,d_0^0,a)\}$ is contained in the rectangle, $rectangle$, and restricts to configuration homomorphisms.
Suppose further that the values in the sequence,
$\{(d_1(n),p(n),d_0^0,a)\}$, belong to a fixed collection of cosets of finite index subgroups of the abelian vertex groups in the completion
$Comp_1$, and this collection is not in the projection to the axis that is associated with $Comp_1$, of any of the product domains that
its projection to the axis that is associated with $Comp_2$ contains the given stratum.

Then no (infinite) subsequence of the  sequence $\{d_1(n),p(n),d_2^0,q_2,a)\}$ (which is a
sequence of values in $duo$) restricts to configuration homomorphisms that can  be extended to
values in  the (boundedly many) rectangles that are associated with the value $q_1$, and with those duo
limit groups,
$Duo_1,\ldots,Duo_t$, that were associated with completions in the initial level of the diagram $Diag$.
\endroster

Since there are boundedly many strata in the stratification of the axis that is associated with $Comp_2$, we have chosen at most
boundedly many values $q_2$ of the parameters $q$.
We continue with all the boundedly many pairs of values, $(q_1,q_2)$, where $q_1$ is the value of the parameters $q$ that was chosen
for the initial level of the diagram $Diag$, and $q_2$ are all the boundedly many
values of the parameters $q$ that were chosen for the various strata in the stratification of the axis that is associated
with $Comp_2$.

With each such pair, $(q_1,q_2)$, we associate the boundedly many rectangles that are associated with it, and with the
duo limit groups, $Duo_1,\ldots,Duo_t$, that are associated with the second level of the diagram $Diag$. As there are boundedly
many pairs, $(q_1,q_2)$, and with each pair there are at most boundedly many associated rectangles, we have altogether associated
boundedly many rectangles with the second level of the diagram $Diag$.

At this stage we repeat what we did with the rectangles of the duo limit groups, $Duo_1,\ldots,Duo_t$, that are associated
with the initial level of the diagram $Diag$, and analyze the (bounded) collection of rectangles of the
duo limit groups, $Duo_1,\ldots,Duo_t$, that appear in the first two levels of the diagram $Diag$. We first associate with
collection of rectangles a finite collection of duo closures of the given duo limit group that is associated with the
rectangle, $rectangle$. With this collection of duo closures we associate a finite planar diagram with axes that consists of
collections of cosets of finite index subgroups of abelian vertex groups that appear in the various levels of the completions,
$Comp_1$ and $Comp_2$. In this diagram we indicate all the collections of cosets that are covered by the closures that were
constructed from the rectangles of duo limit groups that are associated with the first two levels of the diagram $Diag$. By
proposition 3.5, the collections that are covered by these closures are the union of boundedly many product domains. These product
domains give rise to a stratification of the axes that is associated with collection of cosets of abelian vertex groups in
$Comp_2$, and in this stratification there are boundedly many strata.
As we did in the first step of the diagram $Diag$, with each stratum of this stratification we associate a value
$q_3$ of the parameters $q$ that satisfy properties (i)-(iii).

We continue iteratively. At level $m$ of the diagram $Diag$, we look at all the boundedly many $m$-tuples of values of
the parameters $q$, $(q_1,\ldots,q_m)$, that were chosen in the previous $m-1$ levels. With each such $m$-tuple, we associate
the boundedly many  rectangles that are associated with it, and with the duo limit groups, $Duo_1,\ldots,Duo_t$, that are associated
with the $m$-th level of the diagram $Diag$. Given the boundedly many rectangles that are associated with the chosen
values of the parameter subgroup $q$, and with the duo limit groups,
$Duo_1,\ldots,Duo_t$, that appear in all the first $m$ levels of the diagram $Diag$, we construct finitely many
duo closures of the duo limit group that is dual to the given rectangle, $rectangle$.

\noindent
As we did in the first two steps of the diagram $Diag$, with this collection of duo closures we associate a finite planar
diagram.
In this diagram we indicate all the collections of cosets that are covered by the closures that were
constructed from the rectangles of duo limit groups that are associated with the first $m$ levels of the diagram $Diag$. By
proposition 3.5, the collections that are covered by these closures are the union of boundedly many product domains. These product
domains give rise to a bounded stratification of the axes that is associated with collection of cosets of abelian vertex groups in
$Comp_2$.
As we did in the first step of the diagram $Diag$, with each stratum of this stratification we associate a value
$q_{m+1}$ of the parameters $q$ that satisfy properties (i)-(iii).

5with the specialization $q$ in the next (second) level of that diagram.
The iterative process that we presented terminates with the duo limit groups, $Duo_1,\ldots,Duo_t$, that are associated
with the terminal level of the (finite) diagram $Diag$. By construction, when the process terminates we have associated boundedly
many rectangles (that are all associated with the duo limit groups, $Duo_1,\ldots,Duo_t$), with the given rectangle $rectangle$.
From the universality of the diagram $Diag$,
we obtain the covering property, that concludes the proof of theorem 3.3.

\vglue 1pc
\proclaim{Proposition 3.7} The bounded collection of rectangles of the universal duo limit
groups, $Duo_1,\ldots,Duo_t$, that were constructed iteratively by going through the levels of the
universal diagram $Diag$,
covers the given rectangle $duo$.
\endproclaim

\nfp The completions that appear in the initial level of the diagram, $Diag$, are completions of the resolutions in
the graded Makanin-Razborov diagrams of the maximal configuration limit groups of the given rigid limit group, $Rgd(x,p,q,a)$,
with respect to the parameter subgroup $<q>$. Hence, by the universality of the Makanin-Razborov diagrams and the maximal
configuration limit groups, given a value $q_1$ of the parameters $q$, the boundedly many fibers of the completions that appear
in the initial level of the diagram $Diag$, and are associated with the value $q_1$, restrict to all the possible
values $p_0$ of the variables $p$, so that the pair $(p_0,q_1)$ can be extended to a configuration homomorphism that is
associated with the given rigid limit group $Rgd(x,p,q,a)$, i.e., a configuration homomorphism of one of the maximal
configuration limit groups that are associated with $Rgd(x,p,q,a)$.

The duo limit group that is dual to the given rectangle, $rectangle$, is an amalgamated product of two completions,
$<d_1>=Comp_1(d_1,p,a)$ and $<d_2>=Comp_2(d_2,q,a)$, that are amalgamated along the coefficient group $F_k=<a>$.
The value $q_1$ of the parameters $q$ was chosen so that it extends to a value in the rectangle $rectangle$, that
restricts to a (non-degenerate) configuration homomorphism. Hence, $q_1$ extends to a value $d_2^0$ of the variables $d_2$,
such that from every test sequence of the completion $Comp_1$  it is possible to pass to a subsequence, $\{d_1(n)\}$, so that
all the combined values, $\{(d_1(n),d_2^0)\}$, restrict to (non-degenerate) configuration homomorphisms.

\noindent
Therefore, from every test sequence of $Comp_1$ it is possible to extract a subsequence, $\{d_1(n)\}$, such that the restrictions
of the values $d_1(n)$ to the variables $p$ extend to values in the boundedly many fibers that are associated with the
completions that appear in the initial level of the diagram
$Diag$, and with the value $q_1$.

With the rectangle $rectangle$, and the boundedly many rectangles that are associated with those duo limit groups,
$Duo_1,\ldots,Duo_t$, that are associated with the completions in the initial level of $Diag$, and with the value $q_1$, we
have associated a finite planar diagram. The planar diagram and the boundedly many product domains in it (see proposition 3.5),
give a bounded stratification of the axis of the planar diagram that in which there collections of cosets of finite
index subgroups of abelian vertex groups in the various levels of $Comp_2$. In each of the boundedly many strata we chose
an element $q_2$, that satisfy properties (i)-(iii) above.

Therefore, $q_2$ extends to a value $d_2^0$ of the variables $d_2$, such that
from every test sequence of $Comp_1$ that restrict to values of the abelian groups in the various levels
of $Comp_1$, that do not belong to a collection of cosets of finite index subgroups of these abelian groups that is in the
projection of
a planar domain that projects to the stratum of $q_2$, it is possible to extract a subsequence, $\{d_1(n)\}$, so that
all the combined values, $\{(d_1(n),d_2^0)\}$, restrict to (non-degenerate) configuration homomorphisms.
Furthermore, the
restrictions
of the values $d_1(n)$ to the variables $p$ extend to values in the boundedly many fibers that are associated with the
completions that appear in the second level of the diagram
$Diag$, and with the pair $(q_1,q_2)$.

We continue by applying this argument iteratively. At each level the restrictions to the variables $p$ of the values in the
fibers that are associated with the
completions in level $m$ of the diagram $Diag$, and with a tuple, $(q_1,ldots,q_m)$, contain all the restrictions
to the variables $p$ of test sequences of $Comp_1$ for which the
restrictions of the values in the test sequence of $Comp_1$ to values of the abelian groups in the various levels
of $Comp_1$,  do not belong to a collection of cosets of finite index subgroups of these abelian groups that is in the
projection of
a planar domain that projects to the stratum of a fixed extension of $q_m$  to a specialization of $Comp_2$.

Since the diagram $Diag$ is finite by proposition 3.4, when we get to the last level of the diagram,
and there is no level to continue to, the
boundedly many product domains  that are associated with the (boundedly many) rectangles that
are associated with all the duo limit groups, $Duo_1,\ldots,Duo_t$, and the boundedly many values
$(q_1,\ldots,q_i)$, $i=1,\ldots,m$, cover the entire planar diagram that we constructed from the
collections of cosets of finite index subgroups of abelian vertex groups in $Comp_1$ and $Comp_2$.
This proves that these boundedly many rectangles of $Duo_1,\ldots,Duo_t$ cover the given rectangle, $rectangle$,
 of the given duo limit group $Duo$.

\line{\hss$\qed$}

\smallskip
{\rm Remark:} For a given rigid limit group, $Rgd(x,p,q,a)$, and a positive integer $s$, theorem 3.3
proves the existence of finitely many universal duo limit groups, so that every rectangle that is
associated with the corresponding set $NR_s$, is covered by boundedly many rectangles of the universal duo limit
groups. The diagram $Diag$ that is used in the proof of theorem 3.3 generalizes to solid limit groups,
and so is the construction of a finite collection of universal duo limit groups that are associated with it.
Given an arbitrary 
 duo limit group that is associated with a solid limit group, and a rectangle of this duo limit group, it
is not difficult to show that there are finitely many rectangles of the universal duo limit groups that cover
that rectangle (see definition 3.2 for these notions). However, it remains open if there exists a global
bound on the required number of the covering rectangles. Furthermore, 
using the notion of $duo$ $envelopes$ of
a general definable set, that is presented in section 1 in [Se9], it is possible to generalize the statement
of theorems 3.3 to rectangles in general definable sets  (over a free or a torsion-free hyperbolic group). The
validity of the statement for general definable sets remains open as well.

\vglue 1.5pc
\centerline{\bf{\S4. Rigid and Solid Values}}
\medskip
In section 1 we have shown that in the minimal (graded) rank case Diophantine sets are 
equational, and then used it to show that the sets $NR_s$ ($NS_s$), that indicate those values of the 
parameter set $<p,q>$, for which a minimal (graded) rank rigid (solid) limit group $Rgd(x,p,q,a)$ 
($Sld(x,p,q,a)$) admits
at least $s$ rigid (strictly solid families of) values, are in the Boolean algebra generated by
 (minimal rank) equational sets  (theorems
1.5 and 1.7). 

\noindent
In section 2 we have shown that Diophantine sets are equational in the
general case, omitting the minimal (graded)  rank assumption. 
In this section
we combine the equationality of general Diophantine sets with the concept of duo limit
groups that is presented in the previous section, 
to show that the sets $NR_s$ and $NS_s$ that are associated with  general rigid and solid limit groups
are stable.

\vglue 1pc
\proclaim{Theorem 4.1} Let $F_k=<a_1,\ldots,a_k>$ be a non-abelian free group, and let
$Rgd(x,p,q,a)$ ($Sld(x,p,a)$) be a rigid (solid) limit group, 
 with respect to the parameter subgroup $<p,q>$.
Let $s$ be a positive integer, and let $NR_s$ ($NS_s$)
be the set of values of the
defining parameters $<p,q>$ for which the rigid (solid) limit group, $Rgd(x,p,q,a)$ ($Sld(x,p,a)$), has
at least $s$ rigid (strictly solid families of) values. Then the set $NR_s$ ($NS_s$) is stable. 
\endproclaim

\nfp To prove the stability of the set $NR_s$ ($NS_s$), we  
 bound  the length of a sequence of couples of values, $(p_1,q_1),\ldots,(p_n,q_n)$,
for which the formula that is associated with the set $NR_s$ ($NS_s$)
 defines a linear order,
i.e., for which $(p_i,q_j) \in NR_s$ if and only if $i<j$.

We start with the construction of the  diagrams that are needed in order
to get the bound on the lengths of linearly ordered sequences of couples. First, we 
associate with
the set $NR_s$ ($NS_s$) the finite diagram $Diag$ that was constructed
in proving theorem 3.3 (the construction of the diagram $Diag$ that is presented in the rigid case
in theorem 3.3, generalizes in a straightforward way to the solid case).
Recall, that in each step of the diagram we collected all the values:
$$(\{x^i_1,\ldots,x^i_s\}_{i=1}^{\ell},p_0,q_1,\ldots,q_{\ell},a)$$ for which for all indices $i$, 
$1 \leq i \leq \ell$, the values: 
$(x^i_1,\ldots,x^i_s,p_0,q_i,a)$ are rigid (strictly solid) and distinct (belong to distinct strictly 
solid families). We further apply the construction of (quotient) resolutions that is presented and used
in the general
step of the sieve procedure [Se6], to analyze these values and associate finitely many completions, anvils, developing resolutions, and
possibly carriers and sculpted resolutions with them. 
By proposition 3.4 the construction of the diagram terminates after finitely many
steps, and we obtain a finite diagram, that we denote $Diag$. 
The obtained diagram is a finite directed forest,
where at each vertex of the forest we place a (graded) completion, that is either a closure of a completion in the
previous level or it is a completion of  the developing resolution of an anvil that
was constructed along the iterative procedure
(see the detailed
construction of the diagram and its description in the proof of theorem 3.3). The graded completions 
in the diagram are graded with respect to the parameter subgroups $<q_1>$ (completions 
in the first level of the diagram),
$<q_1,q_2>$ (in the second level), and $<q_1,\ldots,q_m>$ for completions  in the $m$-th level of the diagram.

With each of the graded completions in the diagram $Diag$ 
we associate a finite collection of duo limit groups, precisely as we associated duo limit groups with
the completions that were constructed in  the proof of theorem 3.3 (the construction that is presented in the
proof of theorem 3.3 is in the rigid case, and precisely the same construction works in  the solid case).
 Hence, with the entire set of completions 
in the diagram $Diag$, we associate a finite collection of 
universal duo limit groups: $Duo_1,\ldots,Duo_t$, that are precisely the universal duo limit groups that
appear in the statement of theorem 3.3.

The diagram $Diag$ that we associated with $NR_s$ ($NS_s$) is a directed graph for which in each vertex we
place a closure of a completion in the previous level or the completion of the developing resolution of an anvil that was constructed at that level 
of the corresponding branch of the iterative procedure that constructs the diagram.
We set 
$depth_{NR_s}$ ($depth_{NS_s}$) to be the depth (or the number of levels) 
of the directed graph that is associated with the diagram.

The parameter subgroup of the completions that appear in level $m$ of the diagram is denoted, $<q_1,\ldots,q_m>$.
With each value of (the generators of) 
this parameter subgroup one associated boundedly many fibers of the completions that
are placed in level $m$ of the diagram $Diag$. With each such fiber (in level $m$), 
and a value of the variables $q_{m+1}$,
there are boundedly many fibers of the completions that appear in level $m+1$ of the diagram and are
associated with them.
We further  set $width_{NR_s}$ ($width_{NS_s}$) to be the maximal number 
of fibers of  completions that 
are placed in level $m+1$ of the diagram, and are associated with the same fiber of a completion in level $m$
of the diagram, and the same value of the variables $q_{m+1}$
(where the maximum is over all the possible levels $m$ of the
diagram $Diag$, including level 0, in which case there are no fibers in a previous level,
and with a value of the variables $q_1$ there are
at most  boundedly
many associated fibers of completions that are placed in the initial level of the diagram $Diag$).

By the existence of a global bound on the numbers of rigid values of a rigid limit group, and strictly
solid families of a solid limit group (theorems 2.5 and 2.9 in [Se3]), there exists a bound on the maximal number
of rectangles that are associated with one
of the duo limit groups, $Duo_1,\ldots,Duo_t$, and with a fixed fiber of a completion that is placed 
in a vertex in the diagram $Diag$.
We set $rec_{NR_s}$ ($rec_{NS_s}$) to be this bound.

\medskip
Let $Duo$ be one of the (universal) duo limit groups, $Duo_1,\ldots,Duo_t$, that are associated
with $NR_s$ ($NS_s$), and suppose that $Duo=<d_1,p,a>*_{<d_0>} \, <d_2,q,a>$. 
We view $Duo$ as a graded
limit group with respect to the parameter subgroup $<d_2,q>$. Every value of $Duo$ restricts to a value
of the associated configuration limit group $Conf$, $(x_1,\ldots,x_s,p,q,a)$ (see
definition 3.1 for the properties of a duo limit group. Note that the elements $x_1,\ldots,x_s$ can be written as words
in the elements $d_1$ and $d_2$). 
With $Duo$ we further associate the Diophantine condition that forces the
associated restriction  of the configuration limit group, $Conf$, not to be a configuration
homomorphism (to be degenerate). I.e., either one of the values $(x_i,p,q,a)$ is flexible (not strictly
solid), or two rigid specializations $(x_i,p,q,a)$ and $(x_j,p,q,a)$, $i<j$, coincide (belong to the
same strictly solid family. See definition 1.5 in [Se3] for this Diophantine condition in the solid case).  Note
that this degeneracy condition is a Diophantine condition on specializations of the duo limit group $Duo$, and we call it 
the $degenerating$ Diophantine condition.

By theorem 2.1 Diophantine sets are equational. Hence, given a Diophantine set $D(p,q)$ there exists a global
bound on any strictly decreasing sequence of intersections: $\cap_{i=1}^m \, D(p,q_i)$.
Therefore, starting with the duo limit group $Duo$, viewed as a graded limit group with respect to the parameter
subgroup $<d_2,q,a>$, and the specializations that factor through it,
 there exists a global bound on the length of sequences of values:
 $d_2(1),\ldots,d_2(u)$, of the elements
$d_2$ in the duo limit group $Duo$, for which the  sets of values of the variables $d_1$, $D1_r$, $1 \leq r \leq u$, for which
these values together with the corresponding values $d_2(1),\ldots,d_2(r)$, $1 \leq r \leq u$, 
extend to specializations of $Duo$, and the combined
specializations of $Duo$ satisfy the degenerating Diophantine condition, strictly decreases for
$1 \leq r \leq u$.
We set $length_{NR_S}$
($length_{NS_S}$) to be the
maximum 
of these bounds, 
where the maximum is taken over all the universal duo limit groups $Duo_1,\ldots,Duo_t$.

To get a bound on the cardinality of sets of values $\{(p_i,q_i)\}$ that can be ordered by the
sets $NR_s$ and $NS_s$ we need another invariant of the universal duo limit groups, $Duo_1,\ldots,Duo_t$. 
Let $Duo$ be one of these duo limit groups. By the properties of duo limit groups (definition 3.1)
$Duo=Comp_1(d_1,p,a)*_{<d_0,a>} \, Comp_2(d_2,q,a)$, and there is a map from a (maximal) configuration limit
group $Conf$, that is associated with $Rgd(x,p,q,a)$ ($Sld(x,p,q,a)$), into $Duo$. We denote the image of $Conf$
in $Duo$, $<x_1,\ldots,x_s,p,q,a>$.

Each of the elements $x_{\ell}$, ${\ell}=1,\ldots,s$, can be written in a  normal form with respect to the amalgamated
product, $Comp_1*_{<d_0>} \, Comp_2$. Let $x_{\ell}=u^{\ell}_1v^{\ell}_1 \ldots u^{\ell}_{r_{\ell}}v^{\ell}_{r_{\ell}}$, 
$u^{\ell}_e \in Comp_1$,
$v^{\ell}_e \in Comp_2$, ${\ell}=1,\ldots,s$, be such  normal forms. 

\noindent
We continue by viewing $Comp_1(d_1,p,a)$ as a graded limit group with respect to the parameter subgroup 
$<d_0,p,a>$,
and $Comp_2(d_2,q,a)$ as a graded limit group with respect to the parameter subgroup $<d_0,q,a>$. With 
$Comp_1$ and $Comp_2$, viewed as graded limit groups with respect to $<d_0,p,a>$ and $<d_0,q,a>$ in
correspondence, we associate their graded Makanin-Razborov diagrams. 

Suppose first that we are given a rigid
limit group $Rgd(x,p,q,a)$ and its associated set $NR_s$. 
With each value of the elements  $x_1,\ldots,x_s$, that generate the image of the configuration limit
group $Conf$ in $Duo$, there are associated values of the 
elements $u^{\ell}_e,v^{\ell}_e$, $1 \leq {\ell} \leq s$, $1 \leq e \leq r_{\ell}$. If the values of $x_1,\ldots,x_s$ are rigid
values of $Rgd(x,p,q,a)$, then extension of these values to values of $d_1$ and $d_2$ must factor
through graded resolutions in the graded Makanin-Razborov diagrams of $Comp_1(d_1,p,a)$ and $Comp_2(d_2,q,a)$
with respect to the parameter subgroups $<d_0,p,a>$ and $<d_0,q,a>$ in correspondence, in which
the elements $u^{\ell}_e,v^{\ell}_e$ are contained in the distinguished vertex group in all the abelian decompositions 
along the graded resolutions (i.e., the vertex group that contains the subgroups $<d_0,p,a>$ and $<d_0,q,a>$
in correspondence). 

A graded resolution terminates in either a rigid or a solid limit group. By theorems 2.5 and 2.9 in [Se3],
given a value of the parameter subgroups $<d_0,p>$ or $<d_0,q>$, it may extend to only boundedly many values that
a fixed set of generators of the distinguished vertex group in the abelian decomposition that is associated
with the terminal level of one of the graded resolutions in these Makanin-Razborov diagrams. The number of
these values of a fixed generating set of the distinguished vertex group is bounded by the number of rigid
or families of strictly solid families of values, that extend a given value of the subgroup $<d_0,p>$ or
$<d_0,q>$.
We set $excep_{NR_s}$ to be the sum of the bounds on the number of rigid or strictly solid families of values, that are
associated with a given value of $<d_0,p>$ and $<d_0,q>$), where the sum is taken over the terminal
rigid or solid limit groups of all the graded resolutions that appear in
the graded Makanin-Razborov diagrams of $Comp_1$ and $Comp_2$,  for all the duo limit groups $Duo_1,\ldots,Duo_t$.

Suppose  that we are given a solid
limit group $Sld(x,p,q,a)$ and its associated set $NS_s$. 
The duo limit group $Duo$ admits a free product with amalgamation: 
$Duo=Comp_1(d_1,p,a)*_{<d_0>} \, Comp_2(d_2,q,a)$.  
Given a resolution in the graded Makanin-Razborov
diagram of $Comp_1(d_1,p,a)$ with respect to the parameter subgroup $<p,d_0>$, and a resolution in the
graded Makanin-Razborov diagram of $Comp(d_2,q,a)$ with respect to the parameter subgroup $<q,d_0>$,
we take the completions of these two graded resolutions, and then the (finitely many) maximal
limit quotients of the amalgamated product of these two
completions. This amalgamation is a duo limit group that we denote $PQDuo$ (we used the same
construction and notation in the proof of proposition  3.10).  

The image of the configuration limit group in $Duo$, $<x_1,\ldots,x_s,p,q,a>$, is naturally mapped into $PQDuo$.
This image of the configuration limit group in $PQDuo$ restricts to $s$ images of
 the solid  limit group that we have started with, $Sld(x,p,q,a)$, into $PQDuo$. A non-degenerate homomorphism
of the configuration limit group, restricts to $s$ strictly solid values, $(x_i,p,q,a)$, $i=1,\ldots,s$,
of the solid
limit group $Sld(x,p,q,a)$. 

Hence, if such a non-degenerate homomorphism extends to a value of
the duo limit group $Duo$, and that value factors through $PQDuo$, then in the $s$ maps of $Sld(x,p,q,a)$ into
$PQDuo$, the image of every rigid vertex group, every edge group, and every subgroup that is generated by edges
that are adjacent to an abelian vertex group in the graded abelian decomposition that is associated with the
solid limit group $Sld(x,p,q,a)$, must be elliptic in all the abelian decompositions that are associated with
the various levels of $PQDuo$ (i.e., in all the abelian decompositions that are associated with 
the the two completions from which $PQDuo$ is composed).

Therefore, like in the rigid case, and by the global bounds on the number of rigid and strictly solid families
of values of rigid and strictly solid limit groups (theorems 2.5 and 2.9 in [Se3]), those elements  
$u^{\ell}_e,v^{\ell}_e$, $1 \leq {\ell} \leq s$, $1 \leq e \leq r_{\ell}$, that appear in the normal form
of the elements, $x_1,\ldots,x_s$, that generate rigid vertex groups, edge groups, or the subgroups that
are generated by
the edge groups that are adjacent to an abelian vertex group in the graded abelian decomposition of
$Sld$, admit only boundedly many   values for every possible value of the defining parameters $<d_0,p>$ and 
$<d_0,q>$
(in correspondence),
for each terminal rigid or solid limit group of a resolution in the graded Makanin-Razborov diagram of
$Comp_1$ or $Comp_2$.
We set $excep_{NS_s}$ to be the sum of the bounds on the number of such rigid and strictly solid families, 
where the sum is over all the terminal rigid and solid limit groups of graded resolutions that appear in
the graded Makanin-Razborov diagrams of $Comp_1$ and $Comp_2$, for all the duo limit groups $Duo_1,\ldots,Duo_t$.

\vglue 1pc
\proclaim{Proposition 4.2} With the notation of theorem 4.1, let: 
$(p_1,q_1),\ldots,(p_n,q_n)$ be a sequence of values of the defining
parameters $p,q$.
Suppose that  $(p_i,q_j) \in NR_s$ if and only if $i<j$. Then $n<M$ where:
$$M=(1+width_{NR_s})^{(depth_{NR_s} \cdot L_1)} \ ;  \ L_1=(t \cdot rec_{NR_s})^{t \cdot rec_{NR_s} \cdot L_2}$$  
$$L_2={excep_{NR_s}}^{L_3} \ ; \
  L_3={excep_{NR_s}}^{L_4} \ ; \ L_4 = length_{NR_s}+2$$
and a similar statement  holds for the sets $NS_s$, if we replace the constants for $NR_s$ with those for
$NS_s$.
\endproclaim

\nfp We prove the proposition for a set $NR_s$ (that is associated with a rigid
limit group). The proof for the sets $NS_s$ (that are associated with sold limit groups) is identical. Let  $n \geq M$ and:
$(p_1,q_1),\ldots,(p_n,q_n)$ be a sequence of values of the parameters $p,q$, for which
 $(p_i,q_j) \in NR_s$ if and only if $i<j$. By the definition of the set $NR_s$, for every $i<j$, 
there exists an $s$-tuple of values: 
$x^{i,j}=(x_1^{i,j},\ldots,x_s^{i,j})$, so that for every $1 \leq m \leq s$,
$(x^{i,j}_m,p_i,q_j,a)$ is a rigid value  of the
given rigid limit group $Rgd(x,p,q,a)$, and for $1 \leq m_1 < m_2 \leq s$, the corresponding
rigid values are distinct. For the rest of the argument, with each couple $(p_i,q_j)$, $i<j$,
we further associate such an $s$-tuple of values  $x^{i,j}$.

We iteratively filter the tuples $(x^{i,j},p_i,q_j)$, and then apply a simple pigeon-hole principle. 
We start with $q_n$. By the construction of the diagram $Diag$, that is associated
with the rigid limit group, $Rgd(x,p,q,a)$,  
at least $\frac {1} {width_{NR_s}}$ of the
values, $\{(x^{i,n},p_i,q_n,a)\}_{i=1}^{n-1}$, belong to the same fiber that is associated with
$q_n$ in the initial level of the diagram $Diag$. We proceed only with those indices $i$ for which the values, 
$(x^{i,n},p_i,q_n,a)$, belong to that fiber.

We continue with the largest index $i$, $i<n$, for which the tuple, 
$(x^{i,n},p_i,q_n,a)$, belongs to that fiber. We denote that index $i$, $u_2$. By the structure of the diagram
$Diag$,  at least $\frac {1} {1+width_{NR_s}}$ of the values
$\{(x^{i,u_2},x^{i,n},p_i,q_{u_2},q_n,a)\}$, for those indices $i<u_2$ that remained after the first filtration,
belong to either a closure of the same fiber in the initial level, or to one of the fibers in the second level of the
constructed diagram.
We proceed only with those indices $i$ for which the value, 
$\{(x^{i,u_2},x^{i,n},p_i,q_{u_2},q_n,a)\}$, 
belong to either a closure of the same fiber in the initial level, or to the same fiber in the second level of the
diagram $Diag$.

We proceed this filtration process iteratively. The diagram $Diag$ is finite and has depth, $depth_{NR_s}$.
At each step we remain with at least $\frac {1} {1+width_{NR_s}}$ of the values that we have
started the step with, and either we stay with the same fiber that we reached in the previous step,
or we continue to a fiber of a completion that is placed 
in the next level of the diagram. Since we have started with $n \geq M$
pairs of values, $\{(p_i,q_i)\}$, there must exist a subsequence (still denoted) $\{(p_i,q_i)\}_{i=1}^{L_1}$,
for which:
\roster
\item"{(i)}" $(p_i,q_j) \in NR_s$ if and only if $i<j$.

\item"{(ii)}" there exists a fiber of one of the completions that is placed in a vertex of the diagram
$Diag$, so that for $i<j<L_1$ the value $(p_i,q_j)$ extends to a (non-degenerate) configuration homomorphism: 
$(x^{i,j},p_i,q_j,a)$, that further extends to a value of one of the boundedly many  rectangles that are
associated with the fixed fiber and with one of the universal duo limit groups, $Duo_1,\ldots,Duo_t$.  
\endroster

We continue with the subsequence of pairs $\{(p_i,q_i)\}_{i=1}^{L_1}$ that satisfy properties (i) and (ii)
and further filter it. By construction there are $t$ duo limit groups, and with any given fiber of one
of the completions in the diagram $Diag$ and one of the duo limit groups, $Duo_1,\ldots,Duo_t$, there are
at most $rec_{NR_s}$ associated rectangles. Hence the 
sequence of values: $(p_i,q_j)$, $1 \leq i < j \leq L_1$, extends to (non-degenerate) configuration homomorphisms,
$(x^{i,j},p_i,q_j,a)$, $1 \leq i < j \leq L_1$, that further extend to values of at most $t \cdot rec_{NR_s}$
rectangles in the duo limit groups, $Duo_1,\ldots,Duo_t$.

By filtering the sequence of values, $\{(p_i,q_i)\}$, $1 \leq i < j \leq L_1$, according to the rectangle 
that contains the extended configuration homomorphism, $(x^{i,j},p_i,q_j,a)$, using a similar filtration as was
used to filter the subsequence that satisfies properties (i) and (ii), we get a new subsequence (still
denoted) $\{(p_i,q_i)\}_{i=1}^{L_2}$, for which:
\roster
\item"{(1)}" $(p_i,q_j) \in NR_s$ if and only if $i<j$.

\item"{(2)}" there exists a rectangle that is associated with one of the duo limit groups,
$Duo_1,\ldots,Duo_t$, 
so that for $i<j<L_2$ the value $(p_i,q_j)$ extends to a (non-degenerate) configuration homomorphism: 
$(x^{i,j},p_i,q_j,a)$, that further extends to a value of that given rectangle.
\endroster

The duo limit group that is associated with the rectangle in part (2) is an amalgamated
product: $Duo=Comp_1*_{<d_0,a>} \, Comp_2$. Viewing the completions, $Comp_1(d_1,p,a)$ and $Comp_2(d_2,q,a)$,
as graded limit groups with respect to the parameter subgroups, $<d_0,p,a>$ and <$<d_0,q,a>$ in correspondence,
we have associated graded Makanin-Razborov diagrams with $Comp_1$ and $Comp_2$, and each graded resolution in
these diagrams terminates in either a rigid or a solid limit group. Each value of the variables
$p$ and $q$, extend to at most $excep_{NR_s}$ rigid or families of strictly solid values of the terminal
rigid and solid limit groups of the graded resolutions in the graded Makanin-Razborov diagrams of $Comp_1$
and $Comp_2$.

Recall that we denote the image of the configuration limit group $Conf$ in the duo
limit group $Duo$, $<x_1,\ldots,x_s,p,a,a>$. Each of the elements $x_1,\ldots,x_s$ can be written in a normal form
as a word in elements, 
$u^{\ell}_e,v^{\ell_e}$, $1 \leq \ell \leq s$, $1 \leq e \leq r_{\ell}$, 
where $u^{\ell}_e \in Comp_1$ and $v^{\ell}_e \in Comp_2$. 

By filtering the sequence of values, $\{(p_i,q_i)\}$, $1 \leq i < j \leq L_2$,  according to the boundedly many
possible extensions of the values $q_i$ to a rigid or a strictly solid (family of) values of a 
terminal rigid or solid limit group of one of the
finitely many graded   resolutions in the graded Makanin-Razborov diagram of $Comp_2(d_2,q,a)$ with respect to
the parameter subgroup $<d_0,q>$, we are left with a sequence (still denoted),  
$\{(p_i,q_i)\}$, $1 \leq i < j \leq L_3$,  that satisfy properties (1) and (2). Furthermore for each  pair of indices,
$1 \leq i < j \leq L_3$, the associated (non-degenerate) configuration homomorphism, $(x^{i,j},p_i,q_j,a)$, restricts to
values of the elements $v^{\ell}_e$  that depend only on the index $j$, and not on the index $i$, i.e., these values can
be associated with the values $q_j$.

By further filter the sequence of values, $\{(p_i,q_i)\}$, $1 \leq i < j \leq L_3$,  according to the boundedly many
possible extensions of the values $p_i$ to a rigid or a strictly solid (family of) values of a terminal 
rigid or solid limit group of one of the
finitely many graded   resolutions in the graded Makanin-Razborov diagram of $Comp_1(d_1,p,a)$ with respect to
the parameter subgroup $<d_0,p>$, we are left with a sequence (still denoted),  
$\{(p_i,q_i)\}$, $1 \leq i < j \leq L_4$,  that satisfy properties (1) and (2). Furthermore for each  pair of indices,
$1 \leq i < j \leq L_4$, the associated (non-degenerate) configuration homomorphism, $(x^{i,j},p_i,q_j,a)$, restricts to
values of the elements $v^{\ell}_e$  that depend only on the index $j$, and not on the index $i$,
and values of the elements $u^{\ell}_e$  that depend only on the index $i$, and not on the index $j$. I.e.,
these values of the elements $u^{\ell}_e$ and $v^{\ell}_e$ can be associated with the values $p_i$ and $q_j$ in
correspondence.

Finally $L_4=length_{NR_s}+2$. For the last sequence,
$\{(p_i,q_j)\}$, $1 \leq i < j \leq L_4$, $(p_i,q_j) \in NR_s$ if and only if $i<j$, and  with each value $p_i$ we can associate
a fixed value of the elements $u^{\ell}_e$, and with each value $q_j$ we can associate a fixed value of $v^{\ell}_e$, hence,
with each value $p_i$ we can associate a fixed value of $Comp_1$, and with each value of $q_j$ we can associate
a fixed value of $Comp_2$.
Therefore, since $(p_i,q_j) \in NR_s$ if and only if $i<j$, starting with the duo limit group $Duo$, 
viewed as a graded limit group with respect to the parameter
subgroup $<d_2,q,a>$, we obtained a sequence of values:
 $d_2(1),\ldots,d_2(length_{NR_s}+1)$, of the elements
$d_2$ in the duo limit group $Duo$ (the fixed generators of $Comp_2$), for which the  sets of values of the variables $d_1$, 
$D1_r$, $1 \leq r \leq length_{NR_s}+1$, for which
these values together with the corresponding values $d_2(1),\ldots,d_2(r)$, $1 \leq r \leq length(NR_s)+1$, 
extend to values of $Duo$, and the combined
values of $Duo$ satisfy the degenerating Diophantine condition, strictly decreases for
$1 \leq r \leq length_{NR_s}+1$. This contradicts the choice of $length_{NR_s}$ to be a global bound on the length of such
strictly decreasing sequences of values of the variables $d_2$ for all the rectangles in all the duo limit groups,
$Duo_1,\ldots,Duo_t$.

\line{\hss$\qed$}

Proposition 4.2 proves the stability of the sets $NR_s$ and $NS_s$.

\line{\hss$\qed$}

Theorem 4.1 proves the stability of the sets $NR_s$ and $NS_s$, i.e., 
sets of values of the defining parameters for which a rigid or solid limit group have at least $s$ rigid or strictly solid families of
values are stable. Since stable sets are closed under Boolean operations, this proves
that sets of values of the defining parameters for which there are precisely $s$ rigid or strictly solid families of values (of
a given rigid or solid limit group) are stable.
As we did in the minimal rank case (theorem 1.8), in order to prove that the theory of a free group is stable,
i.e., that a general definable set over a free group is stable, 
we need to analyze the (definable) 
set of values of the defining
parameters for which a given (finite) collection of covers of a graded resolution forms a covering
closure  (see definition 1.16 in [Se2] for a covering closure).

\vglue 1pc
\proclaim{Theorem 4.3 (cf. theorem 1.8)} Let $F_k=<a_1,\ldots,a_k>$ be a non-abelian free group, 
let $G(x,p,q,a)$ be a graded
limit group (with respect to the parameter subgroup $<p,q>$), and let $GRes(x,p,q,a)$ be a well-structured
graded resolution
of $G(x,p,q,a)$ that terminates in the rigid (solid) limit group, $Rgd(x,p,q,a)$ ($Sld(x,p,q,a)$). 

Let $GCl_1(z,x,p,q,a),\ldots,GCl_v(z,x,p,q,a)$ be a given set of graded closures of $GRes(x,p,q,a)$. Then the
set of specializations of the parameters $<p,q>$ for which the given set of closures forms a covering closure
of the graded resolution $GRes(x,p,q,a)$, $Cov(p,q)$, is stable.
\endproclaim

\nfp The proof is based on the arguments that were used to prove theorems 1.8 and 4.1.
We start with the construction of the  diagrams that are needed in order
to get the bound on the lengths of linearly ordered sequences of couples for $Cov(p,q)$.

We begin with the  construction of a diagram that is similar in nature to the diagram $Diag$ that was
constructed in analyzing the sets 
$NR_s$ and $NS_s$ (in proving theorems 4.1 and 3.3). The construction  starts with the same collection of
values as we did in analyzing the sets $Cov(p,q)$ in the minimal
(graded) rank case (theorem 1.8).

We look at the entire collection of values: 
$$(x_1,\ldots,x_s,y_1,\ldots,y_m,r_1,\ldots,r_s,p,q,a)$$ for
which (cf. the proof of theorem 1.8):
\roster
\item"{(i)}" for the tuple $p,q$ there exist precisely $s$ rigid (strictly solid families of) values
of the rigid (solid) limit group, $Rgd(x,p,q,a)$ ($Sld(x,p,q,a)$), and at least (total number of) 
$m$ distinct rigid 
and strictly solid families of values of the terminal (rigid and solid) limit groups of the
closures: 
$GCl_1(z,x,p,q,a),\ldots,GCl_v(z,x,p,q,a)$.

\item"{(ii)}" in case the terminal limit groups of $GRes$ is rigid, the $x_i$'s are the distinct rigid 
values of $Rgd(x,p,q,a)$. In case the terminal limit group of $GRes$ is solid, the $x_i$'s  belong to 
the $s$ distinct strictly solid families of $Sld(x,p,q,a)$. 

\item"{(iii)}" the $y_j$'s are either distinct rigid values or belong to distinct strictly solid 
families of values of the terminal (rigid or solid) limit groups of the closures: $GCl_1,\ldots,GCl_v$.

\item"{(iv)}" the $r_i$'s are variables that are  added only in case the terminal limit group
of $GRes$ is solid. In this case the $r_i$'s  demonstrate that the (ungraded) resolutions
that are associated with the given closures and the values, $y_1,\ldots,y_m$, form a covering closure of
the (ungraded) resolutions that are associated with the resolution $GRes$ and the values $x_1,\ldots,x_s$.   
These include primitive roots of the values of all the non-cyclic abelian groups, and edge groups,
in the abelian decomposition that is associated with the terminal solid limit group of $GRes$, $Sld(x,p,q,a)$, and
variables that demonstrate that multiples of these primitive roots up to the least common multiples
of the indices of the finite index subgroups of abelian vertex groups along the resolution  $GRes$ 
that are associated with the graded closures, $GCl_1,\ldots,GCl_v$,
factor through the ungraded resolutions that are associated with the values $y_1,\ldots,y_m$ and their
corresponding closures (cf. section 1 of [Se5] in which we added similar variables to form valid
proof statements, that initialize the sieve procedure).
\endroster

We look at the collection of such values that satisfy properties (i)-(iv) for all the possible values of $s$ and $m$ (note that 
$s$ and $m$ are bounded, since the number of rigid values of a rigid limit group and 
the number of strictly solid families of values of a solid limit group that are associated with a given
value of the defining parameters are globally bounded by theorems 2.5 and 2.9 in [Se3]).

For each fixed $s$ and $m$ we associate with the collection of the values that satisfy
properties (i)-(iv) 
its Zariski closure. With the Zariski closure we associate its dual,  i.e., a canonical finite collection of maximal limit groups, 
that we view as graded
with respect to the parameter subgroup $<q>$. With these graded limit groups we associate
the (graded) resolutions that appear in their taut graded Makanin-Razborov diagrams, and the
resolutions that are associated with the various strata in the singular loci of the diagrams. Given the resolutions
in the collections of the taut Makanin-Razborov diagrams for all the possible values of $s$ and $m$, we
iteratively construct a diagram in a similar way to the construction of the diagram  $Diag$ that is associated
with the sets $NR_s$ and $NS_s$, in proving theorems 4.1 and 3.3. This construction terminates after finitely 
many steps (for precisely the same reasons that the construction of the diagram $Diag$ that is associated with 
the sets $NR_s$ and $NS_s$ terminates after finitely many steps - see proposition 3.4), that finally gives us the first diagram
that is associated with the set $Cov(p,q)$, that we denote $Diag_1$. 

Like the diagram $Diag$ that was constructed in proving theorems 3.3 and 4.1, the diagram $Diag_1$ is a directed forest, so that
in each vertex we further place a graded completion of either a closure of a completion in the previous level or
of the developing resolution of the anvil that was constructed at that
step (and branch) of the iterative procedure that constructed the diagram $Diag_1$. With every graded completion that is placed
in one of the finitely many vertices of the diagram $Diag_1$, we further associate a finite collection of duo limit groups
by applying the same construction that associates duo limit groups with the completions that are placed in the vertices of
the diagram $Diag$ in the proofs of theorems 3.3 and 4.1. We denote the union of the collections of  duo limit groups that are associated with
all the vertices in $Diag_1$,
$Duo^1_1,\ldots,Duo^1_{t_1}$.

As we did in proving theorem 4.1, we set 
$depth^1_{Cov}$ to be the depth of the directed graph associated with the  diagram $Diag_1$, and
$width^1_{Cov}$ to be the maximal number of fibers (of completions) in  level $m+1$ of the diagram $Diag_1$, to which one
continues to from a given fiber of a completion in level $m$ of the diagram $Diag_1$, and a given additional value of the parameters $q$,
where the maximum is taken over all the possible levels $m$, all the completions in these levels, all their fibers,
and all the possible values of 
the parameter subgroups (by the finiteness of the diagram $Diag_1$, and the bounds on the number of rigid and families of
strictly solid families of rigid and solid limit groups (theorems 2.5 and 2.9 in [Se3]) there is a global bound on this maximum).

Given each (graded) completion that appears along the diagram $Diag_1$, we associated with it its collection
of (universal) duo limit groups. By definition 3.1,  each of the finitely many associated duo limit group can be written as an
amalgamated product: $Duo=Comp_1(d_1,p,a)*_{<d_0,a>} \, Comp_2(d_2,q,a)$. 
As we did in the proof of proposition 4.2, we view the completion, $Comp_1(d_1,p,a)$, as
a graded limit group with respect to the parameter subgroup $<p,d_0,a>$, and the completion $Comp_2(d_2,q,a)$
as  graded limit group
with respect to the
parameter subgroup $<q,d_0,a>$. With $Comp_1$ and $Comp_2$, viewed as graded limit groups, we associate their graded 
Makanin-Razborov diagrams (with respect to the above two subgroups of parameters). 
By theorems 2.5 and 2.9 in [Se3], there exist global bounds on the number of rigid and strictly solid
families of values (having the same specialization of the parameter subgroup),
for each of the rigid and solid limit groups in these graded Makanin-Razborov diagrams. For each duo limit group,
$Duo^1_1,\ldots,Duo^1_{t_1}$, we look at the sum of these bounds for all the rigid and solid limit groups that appear 
along the two graded Makanin-Razborov diagrams that are associated with the corresponding two completions,
$Comp_1$ and $Comp_2$. We set $excep^1_{Cov}$ to be the maximum of these sums, where the maximum is taken over all the duo
limit groups, $Duo^1_1,\ldots,Duo^1_{t_1}$.

By the construction of the duo limit group $Duo$, in it there is a subgroup:
$<x_1,\ldots,x_s,y_1,\ldots,y_m,r_1,\ldots,r_s,p,q,a>$ that we denote $Wit$. 
Each of the fixed set of generators of this subgroup can be written in  a normal form with
respect to the amalgamated product: 
$Duo=Comp(d_1,p,a)*_{<d_0,a>} \, Comp(d_2,q,a)$.

With a pair of resolutions in the Makanin-Razborov diagrams  of $Comp_1(d_1,p,a)$ and $Comp_2(d_2,q,a)$,
with respect to the 
parameter subgroups $<d_0,p,a>$ and $<d_0,q,a>$ in correspondence, we construct finitely many duo limit groups
by taking the maximal limit quotients of the amalgamation of their completions along the subgroup $<d_0,a>$.
We denote an obtained duo limit group $PQDuo$.
If a specialization of the subgroup  $Wit$ satisfies the properties (i)-(iv), 
 then an extension of this value to values of $d_1$ and $d_2$ must factor through 
one of the duo limit groups $PQDuo$.
The group $Wit$
contains  $s$ images of the rigid or solid limit group that we have started with, as well as images of the terminal
rigid or solid limit groups of the given graded cover resolutions.  Since the specialization
of $Wit$ satisfies the properties (i)-(iv),  those elements  in $Wit$ that are contained in a 
rigid vertex group, or an edge group, or  in the group that is generated by the edge groups that are adjacent to
an abelian vertex group in the abelian decompositions that are associated with the various rigid and solid limit
groups that are mapped into $Wit$, must be contained in 
rigid vertex groups, or in edge groups, or in subgroups that are
generated by edge groups in abelian vertex groups, in all the abelian decompositions 
along the duo limit groups $Wit$ through which specializations of $Wit$ that satisfy properties (i)-(iv) factor
(i.e., the modular groups that are associated with these abelian decompositions 
do not change their conjugacy class).

Hence, by theorems 2.5 and 2.9 in [Se3] and using our notation,
a value of the parameter subgroups $<d_0,p>$ or $<d_0,q>$, may extend to at most $excep^1_{Cov}$ families of specializations
of the subgroup $Wit$ that satisfy the properties (i)-(iv), and hence to at most $excep^1_{Cov}$ families of specializations
of the elements that appear in a fixed normal forms of a fixed set of generators of $Wit$, and are contained in
rigid vertex groups, edge groups, or subgroups that are generated by edge groups that are adjacent to an
abelian vertex group in one of the rigid or solid limit groups that are mapped into $Wit$.

By the construction of the duo limit groups that are associated with the completions that appear
in the diagram $Diag_1$, and as we did in the proof of proposition 4.2, given a fiber of one of the 
completions that appear along the diagram, $Diag_1$, there are at most boundedly many
rectangles of the duo limit groups, $Duo^1_1,\ldots,Duo^1_{t_1}$, that are associated with that fiber. We set
$rec^1_{Cov}$ to be that bound on the number of rectangles (of $Duo^1_1,\ldots,Duo^1_{t_1}$) 
that can be associated with a fiber (of a completion in $Diag_1$).

Let $Duo$ be one of the (universal) duo limit groups $Duo^1_1,\ldots,Duo^1_{t_1}$ that are associated with
the completions in the diagram $Diag_1$.
Suppose that $Duo=Comp_1(d_1,p,a)*_{<d_0,a>} \, Comp_2(d_2,q,a)$. We now view $Duo$ as a graded
limit group with respect to the parameter subgroup $Comp_2(d_2,q,a)$. With each specialization of $Duo$ there
exists an induced value of the generators of the subgroup $Wit$:  
$(x_1^0,\ldots,x_s^0,y_1^0,\ldots,y_m^0,r_1^0,\ldots,r_s^0,p_0,q_0,a)$ (see
the construction of the diagram $Diag_1$). In the construction of the diagram $Diag_1$, and its duo limit
groups, $Duo^1_1,\ldots,Duo^1_{t_1}$, these values were assumed to satisfy (the non-degeneracy) properties (i)-(iv) above.

With the specializations of $Duo$ we further associate the Diophantine condition that 
forces at least one of the conditions (ii)-(iv), that was imposed on the restriction of these specializations
 to values
of the generators of the subgroup $Wit$ 
to fail. I.e., this Diophantine condition either forces
the value of one of the subgroups, $(x_i,p,q,a)$, to be non-rigid (not strictly solid), or
 the value of two of the subgroups, $(x_i,p,q,a)$ and $(x_j,p,q,a)$, $i<j$, to coincide 
(belong to the
same strictly solid family), or it forces the same type of degenerations for the values of the subgroups $<y_i,p,q,a>$,
or one of the values $r_i$, that was assumed to be a primitive element (an element with no proper roots), has a root with an order that 
divides the least common multiple of the indices of the finite index subgroups that are associated with the given set of graded
closures :$GCl_1,\ldots,GCl_v$.

By the equationality of Diophantine set (theorem 2.1),
there exists a 
global bound on the length of a sequence of values, $d_2(1),\ldots,d_2(u)$, of the variables
$d_2$ (that generate $Comp_2(d_2,q,a)$), for which the intersections of the Diophantine sets (of values of $d_1$) 
that are associated with
the prefixes, $d_2(1),\ldots,d_2(m)$, strictly decrease for $1 \leq m \leq u$. 
We set $length^1_{Cov}$
to be the
maximum of these bounds for all the universal duo limit groups $Duo^1_1,\ldots,Duo^1_{t_1}$.

\medskip
After constructing the first diagram that is associated with $Cov(p,q)$, $Diag_1$, and its duo limit groups,
$Duo^1_1,\ldots,Duo^1_{t_1}$, 
we continue with each of the (universal)
duo limit groups, $Duo^1_i$, that is associated with $Cov(p,q)$,
and construct a second diagram, that is similar to the first one.
Let $Duo$ be one of the duo limit groups, 
$Duo^1_1,\ldots,Duo^1_{t_1}$. By the structure of a duo limit group, $Duo$ can be presented as the
amalgamated product: $Duo=Comp_1(d_1,p,a)*_{<d_0,a>} \, Comp_2(d_2,q,a)$. 
We start the second diagram that is associated with $Duo$, that we denote $Diag_2$,
by collecting all the tuples of values:
$$(x_0,p_0,q_0,d^0_1,d^0_0,d^0_2)$$ 
for which:
\roster
\item"{(1)}" the value: 
$(d^0_1,d^0_0,d^0_2)$ is a specialization of the duo limit group $Duo$. The value $d^0_2$ restricts
to $q_0$, and the value $d^0_1$ restricts to $p_0$.

\item"{(2)}" the value, $(d^0_1,d^0_0,d^0_2)$ restricts to a value of the elements:
$$(x_1,\ldots,x_s,y_1,\ldots,y_m,r_1,\ldots,r_s,p_0,q_0,a).$$ These (restricted) values satisfy 
(the non-degeneracy) properties (ii)-(iv) that
are listed in the construction of the first diagram, $Diag_1$, that is associated with $Cov(p,q)$.

\item"{(3)}" the value:
$(x_0,p_0,q_0,a)$ is a rigid or a strictly solid specialization of the terminal rigid or solid
limit group of the
graded resolution $GRes$ that we have started with, and it is an extra rigid or
an extra strictly solid specialization of that terminal rigid or solid limit group. I.e., it does not coincide 
with any rigid value or does not belong to
the strictly solid family of one of the strictly solid values, that are obtained as the restriction of the 
value $(d_1^0,d_0^0,d_2^0)$ to the elements $(x_i^0,p_0,q_0,a)$, $i=1,\ldots,s$.
\endroster

We continue to the next steps of the construction by collecting values in the
same form, 
precisely as we did  in the construction of the first diagram, $Diag_1$, that is associated with $Cov(p,q)$,
where the subgroup $<d_1>$ plays the role of the parameter subgroup $<p>$ (in the construction of $Diag_1$),
and the subgroup $<d_2>$ plays the role of the parameter subgroup $<q>$ in the construction of $Diag_1$.
By the same argument that implies the termination of the construction of the diagram $Diag_1$ (see proposition 3.4),
the constructions of the diagrams that are associated with the 
various duo limit groups, $Duo^1_i$, terminate after finitely many steps.  
We denote each of the diagrams that are
associated with the duo limit groups $Duo^1_i$, $Diag^2_i$.

Recall that like the diagram $Diag_1$, each of the diagrams, $Diag^2_i$, is a finite directed forest, so that at
each vertex of the forest there is a graded completion. With every graded completion that is placed at a vertex
in a diagram, $Diag^2_i$, we further associate a finite collection of duo limit groups, precisely as we associated
the duo limit groups,
$Duo^1_1,\ldots,Duo^1_{t_1}$,
with the graded completions in the diagram $Diag_1$
(see the proof of theorem 3.3 for the construction of these associated duo limit groups).
We denote the finite collection of duo limit groups that are associated with the graded completions in all
the diagrams, $Diag^2_i$, $i=1,\ldots,t_1$,
$Duo^2_1,\ldots,Duo^2_{t_2}$. 

As we did with the first diagram that is associated with $Cov(p,q)$, we set 
$depth^2_{Cov}$ to be the maximal depth of the directed forests that are  associated with the constructed 
(second) diagrams, $Diag^2_i$. We set
$width^2_{Cov}$ to be 
the maximal number of fibers (of completions) in  level $m+1$ of any of the diagrams $Diag^2_i$, to which one
continues to from a given fiber of a completion in level $m$ in that diagram, and a given additional value of (the parameters) $d_2$,
where the maximum is taken over all the diagrams, $Diag^2_i$, all the possible levels $m$, all the completions in these levels, 
all their fibers,
and all the additional possible values of 
the parameter subgroup $d_2$.

Let $Duo$ be one of the duo limit groups, $Duo^2_1,\ldots,Duo^2_{t_2}$. $Duo$ can be written as an amalgamated
product: $Duo=Comp_1(e_1,d_1,a)*_{<e_0,a>} \, Comp_2(e_2,d_2,a)$.
As we did in the proof of proposition 4.2 and with the duo limit groups that are associated with $Diag_1$,
 we view the completion, $Comp_1(e_1,d_1,a)$, as
a graded limit group with respect to the parameter subgroup $<d_1,e_0,a>$, and the completion $Comp_2(e_2,d_2,a)$
as  graded limit group
with respect to the
parameter subgroup $<d_2,e_0,a>$. With $Comp_1$ and $Comp_2$, viewed as graded limit groups, we associate their graded 
Makanin-Razborov diagrams (with respect to the above two subgroups of parameters). 
In each of these graded Makanin-Razborov diagrams there are finitely many rigid and solid limit groups.
By theorems 2.5 and 2.9 in [Se3], there exist global bounds on the number of rigid and strictly solid
families of values (having the same specialization of the parameter subgroup),
for each of the rigid and solid limit groups in these graded Makanin-Razborov diagrams. 
For each duo limit group,
$Duo^2_1,\ldots,Duo^2_{t_2}$, we look at the sum of these bounds for all the rigid and solid limit groups that appear 
along the two graded Makanin-Razborov diagrams that are associated with the corresponding two completions,
$Comp_1$ and $Comp_2$ that are associated with that duo limit group. We set $excep^2_{Cov}$ to be the maximum of these sums, 
where the maximum is taken over all the duo
limit groups, $Duo^2_1,\ldots,Duo^2_{t_2}$.  

Given a fiber of one of the 
completions that appear along one of the diagrams, $Diag^2_i$, there are at most boundedly many
rectangles of the duo limit groups, $Duo^2_1,\ldots,Duo^2_{t_2}$, that are associated with that fiber. We set
$rec^2_{Cov}$ to be that bound on the number of rectangles (of $Duo^2_1,\ldots,Duo^2_{t_2}$) 
that can be associated with a fiber (of a completion in any of the diagrams $Diag^2_i$).

Let $Duo$ be one of the  duo limit groups, $Duo^2_1,\ldots,Duo^2_{t_2}$, that are associated with
the completions in the diagrams $Diag^2_i$.
Suppose that $Duo=Comp_1(e_1,d_1,a)*_{<e_0,a>} \, Comp_2(e_2,d_2,a)$. We view $Duo$ as a graded
limit group with respect to the parameter subgroup $Comp_2(e_2,d_2,a)$. 
With each value of $Duo$ there
exists an associated 
extra value, $(x_0,p_0,q_0,a)$, of the terminal rigid or solid limit group, $Rgd(x,p,q,a)$ or $Sld(x,p,q,a)$, of the 
given graded resolution $GRes$.
In the construction of the diagrams $Diag^2_i$, and its duo limit
groups, $Duo^2_1,\ldots,Duo^2_{t_2}$, these values were assumed to be either rigid or strictly solid, and to
satisfy (the non-degeneracy) property (3) above.

With $Duo$ we further associate the Diophantine condition that forces the
associated extra specialization of the rigid or solid terminal limit group of $GRes$ to be either flexible 
(not rigid or not strictly solid),
or to coincide with one of the rigid values, or to belong to one of the strictly solid families
of values, that appear in the corresponding (induced) value of generators of the subgroup $Wit$: 
$(x_1^0,\ldots,x_s^0,y_1^0,\ldots,y_m^0,r_1^0,\ldots,r_s^0,p_0,q_0,a)$. I.e.  we add a Diophantine condition that
forces the collection of values not to satisfy property (3) in the definition of the collection of
values that are collected in each step of the construction of the  diagrams $Diag^2_i$, and their associated duo limit groups:
$Duo^2_1,\ldots,Duo^2_{t_2}$. 

By the equationality of Diophantine set (theorem 2.1),
there exists a 
global bound on the length of a sequence of values, $e_2(1),\ldots,e_2(u)$, of the variables
$e_2$ (that generate $Comp_2(e_2,d_2,a)$, for which the intersections of the Diophantine sets (of values of $e_1$) 
that are associated with
the prefixes, $e_2(1),\ldots,e_2(m)$, strictly decrease for $1 \leq m \leq u$. 
We set $length^2_{Cov}$
to be the
maximum of these bounds for all the  duo limit groups $Duo^2_1,\ldots,Duo^2_{t_2}$.

\vglue 1pc
\proclaim{Proposition 4.4} With the notation of theorem 4.3, let: 
$(p_1,q_1),\ldots,(p_n,q_n)$ be a sequence of couples of values of the defining
parameters $p,q$
for which $(p_i,q_j) \in Cov(p,q)$ if and only if $i<j$. Then $n<M$ where:
$$M=(1+width^1_{Cov})^{(depth^1_{Cov} \cdot L_1)} \ ; \ L_1=( rec^1_{Cov})^{(rec^1_{Cov} \cdot L_2)}$$  
$$L_2={(excep^1_{Cov})}^{L_3} \ ; \
  L_3={(excep^1_{Cov})}^{L_4} \ ; \ L_4=2^{2L_5}$$
$$L_5=length^1_{Cov}+2+(1+width^2_{Cov})^{(depth^2_{Cov} \cdot L_6)} \ ; \ L_6=( rec^2_{Cov})^{(rec^2_{Cov} \cdot L_7)}$$  
$$L_7={(excep^2_{Cov})}^{L_8} \ ; \
  L_8={(excep^2_{Cov})}^{L_9} \ ; \ 
  L_9 = length^2_{Cov}+2$$
\endproclaim

\nfp The argument that we use is a strengthening of the argument that was used to prove proposition 4.2.
Let  $n \geq M$ and let:
$(p_1,q_1),\ldots,(p_n,q_n)$ be a sequence of values of the parameters $p,q$, for which
 $(p_i,q_j) \in Cov(p,q)$ if and only if $i<j$. By the definition of the set $Cov(p,q)$, for every $i<j$, 
there exists a tuple: 
$$(x^{i,j}_1,\ldots,x^{i,j}_{s_{i,j}},y^{i,j}_1,\ldots,y^{i,j}_{m_{i,j}},r^{i,j}_1,\ldots,r^{i,j}_{s_{i,j}},p_i,q_j,a)$$ 
that satisfies properties (i)-(iv) (which are the properties that values from which we construct the first diagram,
$Diag_1$, that is associated with the set $Cov(p,q)$, and its associated duo limit groups,
$Duo^1_1,\ldots,Duo^1_{t_1}$,
have to satisfy). In the sequel we denote the subgroup that is generated by these elements in the completions that are
placed in $Diag_1$, and its associated duo limit groups, $Wit$.

We iteratively filter the tuples that are associated with the couples $(p_i,q_j)$, 
in a similar way to what we did in proving proposition 4.2. 
We start with $q_n$. By the construction of the first diagram $Diag_1$,
at least $\frac {1} {width^1_{Cov}}$ of the
specializations  of the subgroup $Wit$ that are associated with the values: $(p_i,q_n)$, $1 \leq i \leq n-1$, belong to the same fiber 
that is associated with
$q_n$ in one of the completions that are placed in the initial level of the diagram $Diag_1$. 
We proceed only with those indices $i$ for which the specializations of the subgroup $Wit$ that are associated with the values:
$(p_i,q_n)$
belong to that fiber.

We proceed as in the proof of proposition 4.2. We continue with the largest index $i$, $i<n$, for which the specialization of $Wit$ that
is associated with the tuple 
$(p_i,q_n)$
belongs to that fiber. We denote that largest index $i$, $u_2$. By the structure of $Diag_1$,
at least $\frac {1} {1+width^1_{Cov}}$ of the specializations of $Wit$ that are associated with the values:
$(p_i,q_{u_2})$ and $(p_i,q_n)$,
 for those indices $i<u_2$ that remained after the first filtration,
belong to either the same fiber in the initial level of $Diag_1$, or to a fixed  fiber of a completion that is placed in the second level of the
diagram $Diag_1$.
We proceed only with those indices $i$ for which the specializations of $Wit$  that are associated with the pairs, $(p_i,q_{u_2})$ and $(p_i,q_n)$,
belong either to the initial fiber or to a fixed fiber of a completion in the second level of $Diag_1$.

We proceed this filtration process iteratively (as in the proof of proposition 4.2). 
Since the diagram $Diag_1$ is finite and has depth, $depth^1_{Cov}$,
and since at each step we remain with at least $\frac {1} {1+width^1_{Cov}}$ of the tuples that we have
started the step with, and since $n$, the number of tuples that we started with, satisfies $n \geq M$, 
after we iteratively apply the  filtration process we must obtain
a  subsequence,
 (still denoted) $\{(p_i,q_i)\}_{i=1}^{L_1}$,
for which:
\roster
\item"{($\hat i$)}" $(p_i,q_j) \in Cov(p,q)$ if and only if $i<j$.

\item"{($\hat {ii}$)}" there exists a fiber of one of the completions that is placed in a vertex of the diagram
$Diag_1$, so that for $i<j<L_1$ the value $(p_i,q_j)$ extends to a (non-degenerate) specialization of the subgroup $Wit$
(i.e., a specialization of $Wit$ that satisfies properties (i)-(iv)),
 that further extends to a value of one of the boundedly many  rectangles that are
associated with the fixed fiber and with one of the universal duo limit groups, $Duo^1_1,\ldots,Duo^1_{t_1}$.  
\endroster

We continue as in the proof of proposition 4.2, and further filter the subsequence of pairs $\{(p_i,q_i)\}_{i=1}^{L_1}$ 
(that satisfy properties ($\hat i$) and ($\hat {ii}$)).
By construction there are $t_1$ duo limit groups that are associated with the diagram $Diag_1$. With a given fiber of one
of the completions in the diagram $Diag_1$, and one of the duo limit groups, $Duo_1,\ldots,Duo_t$, there are
at most $rec^1_{Cov}$ associated rectangles. Hence the 
sequence of values: $(p_i,q_j)$, $1 \leq i < j \leq L_1$, extends to (non-degenerate) specializations of the subgroup $Wit$ (i.e., 
specializations that
satisfy properties (i)-(iv)),
that further extend to values of at most $rec^1_{Cov}$
rectangles in the duo limit groups, $Duo^1_1,\ldots,Duo^1_{t_1}$.

By filtering the sequence of values, $\{(p_i,q_j)\}$, $1 \leq i < j \leq L_1$, according to the rectangle 
that contains the extended specializations of the subgroup $Wit$,
using a similar filtration as was
used to filter the subsequence that satisfies properties ($\hat i$) and ($\hat {ii}$), 
we get a new subsequence (still
denoted) $\{(p_i,q_i)\}_{i=1}^{L_2}$, for which:
\roster
\item"{($\hat 1$)}" $(p_i,q_j) \in Cov$ if and only if $i<j$.

\item"{($\hat 2$)}" there exists a rectangle that is associated with one of the duo limit groups,
$Duo^1_1,\ldots,Duo^1_{t_1}$, 
so that for $i<j<L_2$ the value $(p_i,q_j)$ extends to a (non-degenerate) specializations of the subgroup $Wit$ (a value that satisfies properties
(i)-(iv)), 
 that further extends to a value in that given rectangle.
\endroster

The duo limit group that is associated with the rectangle in part ($\hat 2$) is an amalgamated
product: $Duo=Comp_1*_{<d_0,a>} \, Comp_2$. Viewing the completions, $Comp_1(d_1,p,a)$ and $Comp_2(d_2,q,a)$,
as graded limit groups with respect to the parameter subgroups, $<d_0,p,a>$ and $<d_0,q,a>$ in correspondence,
we have associated graded Makanin-Razborov diagrams with $Comp_1$ and $Comp_2$, and each graded resolution in
these diagrams terminates in either a rigid or a solid limit group. Each value of the variables
$p$ and $q$, extends to at most $excep^1_{Cov}$ rigid or families of strictly solid values of the terminal
rigid and solid limit groups of the graded resolutions in the graded Makanin-Razborov diagrams of $Comp_1$
and $Comp_2$.

Recall that given a pair of resolutions,
one in the graded Makanin-Razborov diagrams of $Comp_1(d_1,p,a)$ (with respect to the parameter subgroup
$<p,d_0>$), and a resolution in the graded Makanin-Razborov diagram of $Comp_2(d_2,q,a)$ (with respect to the
parameter subgroup $<q,d_0>$), we constructed from them finitely many duo limit groups, $PQDuo$, that are
the maximal limit quotients of the amalgamation of the completions of the two given resolutions along the
amalgamated subgroup $<d_0,a>$. 

Given  a specialization of the subgroup  $Wit$ that satisfies the (non-degenerate) properties (i)-(iv), 
 the extension of this value to values of $d_1$ and $d_2$ must factor
through duo limit groups $PQDuo$, 
in which
elements in $Wit$ that are contained in the image (in $Wit$)
of rigid vertex groups, edge groups, or subgroups generated
by edge groups that are adjacent to abelian vertex groups, in the abelian decompositions of the given rigid
or solid limit group $Rgd$ ($Sld$) or of the terminal rigid and solid limit groups of the given finite
set of closures,  
are contained in rigid vertex groups, or in edge groups, or in subgroups that are
generated by edge groups that are adjacent to  abelian vertex groups, in all the abelian decompositions 
along the various levels of the duo limit group $PQDuo$.
(i.e., the modular groups that are associated with these abelian decompositions 
do not change their conjugacy class). 

For $1 \leq i < j \leq L_2$, the specializations of the subgroup $Wit$
that are associated with the pairs $(p_i,q_j)$ satisfy properties (i)-(iv). The modular groups that are
associated with the various levels of a duo limit group $PQDuo$ do not change the families of the restrictions
of the corresponding specializations of $Wit$ to values of its associated rigid and solid limit groups, hence, do not
change the fact that such a specialization of $Wit$ satisfies properties (i)-(iv). Hence, we may
assume that the values of the variables $d_1$ and $d_2$ that are associated with the pairs $(p_i,q_j)$
$1 \leq i < j \leq L_2$, are values of the terminal rigid or solid limit groups of the two resolutions
of $Comp(d_1,p,a)$ and  of $Comp(d_2,q,a)$, with respect to the parameter subgroups, $<d_0,p>$ and $<d_0,q>$, in
correspondence.

By filtering the sequence of values, $\{(p_i,q_i)\}$, $1 \leq i < j \leq L_2$,  according to the boundedly many
possible extensions of the values $q_i$ to a rigid or a strictly solid (family of) values of a 
terminal rigid or solid limit group of one of the
finitely many graded   resolutions in the graded Makanin-Razborov diagram of $Comp_2(d_2,q,a)$ with respect to
the parameter subgroup $<d_0,q>$, we are left with a sequence (still denoted),  
$\{(p_i,q_i)\}$, $1 \leq i < j \leq L_3$,  that satisfy properties ($\hat 1$) and ($\hat 2$). 

By further filtering the sequence of values, $\{(p_i,q_i)\}$, $1 \leq i < j \leq L_3$,  according to the boundedly many
possible extensions of the values $p_i$ to a rigid or a strictly solid (family of) values of a terminal 
rigid or solid limit group of one of the
finitely many graded   resolutions in the graded Makanin-Razborov diagram of $Comp_1(d_1,p,a)$ with respect to
the parameter subgroup $<d_0,p>$, we are left with a sequence (still denoted),  
$\{(p_i,q_i)\}$, $1 \leq i < j \leq L_4$,  that satisfy properties ($\hat 1$) and ($\hat 2$).

Furthermore, for each  pair of indices,
$1 \leq i < j \leq L_4$, the values of the pairs $(p_i,q_j)$, extend to  values of the duo limit group $Duo$, hence,
to values of the two completions, $Comp_1$ and $Comp_2$, from which $Duo$ is composed. 
By the filtration that we used, the associated values of the elements $d_1$ (the generators of the completion $Comp_1(d_1,p,a)$), that we may
assume to be  values of one of the rigid or solid limit groups in the graded Makanin-Razborov diagram of $Comp_1$ with respect
to the parameter subgroup $<p,d_0,a>$,  depend only on the index $i$
and not on the index $j$.  
The associated values of the elements $d_2$ (the generators of the completion $Comp_2(d_2,q,a)$), that we may
assume to be  values of one of the rigid or solid limit groups in the graded Makanin-Razborov diagram of $Comp_2$ with respect
to the parameter subgroup $<q,d_0,a>$,  depend only on the index $j$
and not on the index $i$.

\medskip
For the rest of the argument we continue with the sequence of values that we filtered, that we still
denote $(p_i,q_i)$, $1 \leq i \leq L_4$. With each pair of values from this sequence, $(p_i,q_j)$,
there is an associated specialization of
the subgroup $Wit$, and for $1 \leq i < j \leq L_4$, these values satisfy properties (i)-(iv), that testify that
the corresponding pairs, $(p_i,q_j)$, are contained in $Cov(p,q)$. Furthermore, these specializations of the subgroup $Wit$,
extend to values in a fixed rectangle in one of the duo limit groups, $Duo^1_1,\ldots,Duo^1_{t_1}$, that are 
associated with the diagram $Diag_1$. The extensions to values in the rectangle restrict to values of the two
completions, $Comp_1(d_1,p,a)$ and $Comp_2(d_2,q,a)$,  from which the rectangle (or its dual duo limit group) is composed. 
The sequence $(p_i,q_i)$, $1 \leq i \leq L_4$, and its associated specializations of the subgroup $Wit$,
were filtered so that the values of the elements $d_1$ and $d_2$ that extend the corresponding 
specializations of the subgroup $Wit$, were chosen so that the values of $d_1$ depends only on the index $i$, and the value
of $d_2$ depends only on the index $j$. 

We denote the value of the elements $d_1$ that is associated with $p_i$, $d_1(i)$,
and the value of the elements $d_2$ that is associated with $q_j$, $d_2(j)$. The pairs of values $(d_1(i),d_2(j))$ were
filtered from values of the variables $d_1$ and $d_2$ that are associated with pairs with indices, $1 \leq i < j \leq L_4$.
However,  every pair $(d_1(i),d_2(j))$, $1 \leq i,j \leq L_4$, is in the rectangle that is associated with the 
given sequence of values, $\{(p_i,q_j)\}$, and as a value in the rectangle it restricts to a specialization of the subgroup $Wit$.
For indices $1 \leq i < j \leq L_4$, these specializations of $Wit$ satisfy the properties (i)-(iv). 
For indices, $1 \leq j \leq i \leq L_4$, the pairs $(p_i,q_j)$ are not in $Cov(p,q)$, hence, the specializations
of the subgroup $Wit$ that are associated with the corresponding value, $(d_1(i),d_2(j))$, do not satisfy at least
one of the properties (i)-(iv). Therefore, for the last pairs of indices, $1 \leq j \leq i \leq L_4$,
one of the two following properties must hold for each of the associated specialization
of the subgroup $Wit$:
\roster
\item"{(a)}" the value $(d_1(i),d_2(j))$ restricts to a degenerate specialization of the subgroup $Wit$, i.e., to a specialization
of $Wit$ that doesn't satisfy one of the properties (ii)-(iv). In this case, the failure of each of the properties
(i)-(iv) can be translated to a Diophantine condition that the corresponding specialization of $Wit$ has to satisfy, precisely
as the degeneration in the corresponding value of the configuration limit $Conf$, was translated into
a Diophantine condition in the proof of theorem 4.1.

\item"{(b)}" the value 
$(d_1(i),d_2(j))$ restricts to a specialization of the subgroup $Wit$ that doesn't satisfy property (i).
In this case, for the corresponding specialization of the subgroup $Wit$
there exists some extra
rigid or strictly solid value $(x^{i,j}_0,p_i,q_j,a)$ of the 
terminal rigid or solid limit group of the graded resolution $GRes$ that we have started with, and this extra rigid
or strictly solid value does not coincide with a rigid value and does not belong to any strictly
solid family which is a part of the corresponding specialization of the subgroup $Wit$.
\endroster

We continue by filtering the set of values, $\{(p_i,q_i)\}$, $1 \leq i \leq L_4$ according to the two possibilities
(a) and (b).
We start with $q_1$. 
At least half of the
specializations  of the subgroup $Wit$ that are associated with the values: $(p_i,q_1)$, $1 \leq i \leq L_4$, satisfy the same
property, which is either (a) or (b). 
We proceed only with those indices $i$ for which the specializations of the subgroup $Wit$ satisfy that property.
We proceed as in the proof of proposition 4.2. We continue with the smallest index $i$, $1<i$, that satisfy the property
that the majority of the specializations of the subgroup $Wit$ that
are associated with the tuples $\{p_i,q_1)\}$ satisfy. 
We denote that smallest index $i$, $u_2$. 
At least half of the
specializations  of the subgroup $Wit$ that are associated with the values: $(p_i,q_{u_2})$, for those indices, $u_2 \leq i \leq L_4$,
that remained after the initial filtration (the filtration of the pairs $(p_i,q_1)$), satisfy the same
property, which is either (a) or (b). 
We proceed only with those indices $i$ for which the specializations of the subgroup $Wit$  that are associated with the pairs, $(p_i,q_{u_2})$,
satisfy the same property ((a) or (b)), and the specializations of the subgroup $Wit$ that are associated with the pairs, $(p_i,q_n)$,
satisfy the same property (a) or (b).

We proceed this filtration process iteratively (as in the proof of proposition 4.2). $L_4=2^{2L_5}$, and at each 
step we are left with at least half of the pairs that existed in the previous step. Hence, when the iterative filtration
terminates we are left with at least $L_5$ pairs, (still denoted) $\{(p_i,q_i)\}$, $1 \leq i \leq L_5$, so that for every
 pair $(p_i,q_j)$, $1 \leq j \leq i \leq L_5$, the specializations of the subgroup $Wit$ that are associated with these pairs
either all satisfy property (a) or they all satisfy property (b). 

Suppose that the specializations of the subgroup $Wit$ that are associated with the pairs, 
$(p_i,q_j)$, $1 \leq j \leq i \leq L_5$, do all
satisfy property (a), i.e., they all do not satisfy at least one of the properties (ii)-(iv). The failure of
the properties (ii)-(iv) translates to a Diophantine condition that the specializations of $Wit$ need to satisfy, hence,
it translates to a Diophantine condition that the pairs of associated values, $(d_1(i),d_2(j))$, $1 \leq i \leq j \leq L_5$,
need to satisfy ($d_1$ and $d_2$ are the generators of the completions $Comp_1$ and $Comp_2$, in correspondence, that
together generate the duo limit group that is associated with the rectangle that is associated with the
sequence, $\{(p_i,q_i)\}$).
 
Therefore, like in the end of the proof of proposition 4.2, starting with the duo limit group $Duo$ that is
associated with the sequence $\{(p_i,q_i)\}$, $1 \leq i \leq L_5$,
viewed as a graded limit group with respect to the parameter
subgroup $<d_2,q,a>$, we obtained a sequence of values:
 $d_2(1),\ldots,d_2(L_5)$, of the elements
$d_2$ in the duo limit group $Duo$ (the fixed generators of $Comp_2$), for which the  sets of values of the variables $d_1$, 
$D1_r$, $1 \leq r \leq L_5$, for which
these values together with the corresponding values $d_2(1),\ldots,d_2(r)$, $1 \leq r \leq L_5$, 
extend to values of $Duo$, and the combined
value of $Duo$ satisfy the (degenerating) Diophantine condition, that is equivalent to the failure of at least one of
the properties (ii)-(iv) for the corresponding specializations of the subgroup $Wit$, strictly decreases for
$1 \leq r \leq L_5$. Since we assumed that $L_5 \geq length^1_{Cov}+1$, 
this contradicts the choice of $length^1_{Cov}$ 
to be a global bound on the length of such
strictly decreasing sequences of values of the variables $d_2$ for all the rectangles in all the duo limit groups,
$Duo^1_1,\ldots,Duo^1_{t_1}$.

Hence, for the rest of the argument, we may assume that the specializations of the subgroup $Wit$ that are associated with the pairs, 
$(p_i,q_j)$, $1 \leq j \leq i \leq L_5$, do all
satisfy property (b), i.e., that they do not satisfy property (i). The failure of property (i) implies that with each 
specialization of the subgroup $Wit$ that is associated with such a pair, $(p_i,q_j)$,   
there exists some extra
rigid or strictly solid value $(x^{i,j}_0,p_i,q_j,a)$ of the 
terminal rigid or solid limit group of the graded resolution $GRes$ that we have started with, and this extra rigid
or solid value does not coincide with a rigid value and does not belong to any strictly
solid family which is a part of the corresponding specialization of the subgroup $Wit$.

At this point we analyze the sequence of values, 
$(p_i,q_j)$, $1 \leq j \leq i \leq L_5$, and their associated values $d_1(i)$ and $d_2(j)$, 
$1 \leq j \leq i \leq L_5$, by exactly the same argument that was used to prove proposition 4.2, in 
a reverse order (starting with $q_1$ instead of starting with $q_n$).

\medskip
We start with the set of values: 
$(x^{i,j}_0,p_i,d_1(i),q_j,d_2(j),a)$, $1 \leq j \leq i \leq L_5$, so that for every for every pair 
$1 \leq j \leq i \leq L_5$,
the following properties hold:
\roster
\item"{($\hat a$)}" $(d_1(i),p_i,d_2(j),q_j)$ is a value in a fixed rectangle (independent of $i$ and
$j$) that is associated with  one of the duo limit groups, $Duo^1_1,\ldots,Duo^1_{t_1}$  
(that we will denote $Duo^1$ in the sequel). Furthermore, this value restricts to a specialization of the subgroup
$Wit$, that satisfies the non-degeneracy properties (ii)-(iv).

\item"{($\hat b$)}" for each $1 \leq j \leq i \leq L_5$, the value $(x^{i,j}_0,p_,q_j,a)$ is a rigid
or a strictly solid value of the terminal rigid or solid limit group of the given resolution $GRes$, and
this value is distinct from all the $s$ rigid values (not in the same strictly solid families) that are
part of the restriction of the value
$(d_1(i),p_i,d_2(j),q_j)$ to the subgroup $Wit$. I.e., the values
$(x^{i,j}_0,p_i,q_j,a)$ demonstrate that (the non-degeneracy) condition (i) fails for the restriction of the values
$(d_1(i),p_i,d_2(j),q_j)$ to the subgroup $Wit$. 
\endroster

As in the proof of proposition 4.2, we start by iteratively filter the tuples $(x^{i,j}_0,p_i,q_j)$. 
We start with $q_1$. By the construction of the diagram $Diag^2_i$, that is associated with the duo limit group
$Duo^1$,
at least $\frac {1} {width^2_{Cov}}$ of the
values, $(x^{i,1}_0,p_i,q_1)$, $1 \leq i \leq n$, belong to the same fiber 
that is associated with
$q_1$ in one of the completions that are placed in the initial level of the diagram $Diag^2_i$. 
We proceed only with those indices $i$ for which the values, 
$(x^{i,1}_0,p_i,q_1)$, $1 \leq i \leq n$, belong to that fiber. 

We proceed this filtration process iteratively (as in the proof of proposition 4.2). 
Since the diagram $Diag^2_i$ is finite and has depth bounded by  $depth^2_{Cov}$,
and since at each step we remain with at least $\frac {1} {1+width^2_{Cov}}$ of the tuples that we have
started the step with, and since  the number of tuples that we started with is (at least) $L_5$, 
after we iteratively apply the  filtration process we must obtain
a  subsequence,
 (still denoted) $\{(p_i,q_i)\}_{i=1}^{L_6}$,
for which:
\roster
\item"{($\tilde i$)}" $(p_i,q_j) \in Cov(p,q)$ if and only if $i<j$.

\item"{($\tilde {ii}$)}" there exists a fiber of one of the completions that is placed in a vertex of the diagram
$Diag^2_i$, so that for $1 \leq j \leq i \leq L_6$, the value 
$(x^{i,j}_0,p_i,d_1(i),q_j,d_2(j),a)$, that satisfies properties ($\hat a$) and ($\hat b$),  
extends to a value 
  of one of the boundedly many  rectangles that are
associated with the fixed fiber and with one of the universal duo limit groups, $Duo^2_1,\ldots,Duo^2_{t_2}$.  
\endroster

By filtering the sequence of values, $\{(p_i,q_j)\}$, $1 \leq j \leq i  \leq L_6$, according to the rectangle 
that contains the  values that extend the associated values:
$(x^{i,j}_0,p_i,d_1(i),q_j,d_2(j),a)$,   
using a similar filtration as was
used to filter the subsequence that satisfies properties ($\tilde i$) and ($\tilde {ii}$), 
we get a new subsequence (still
denoted) $\{(p_i,q_i)\}_{i=1}^{L_7}$, for which:
\roster
\item"{($\tilde 1$)}" $(p_i,q_j) \in Cov$ if and only if $i<j$.

\item"{($\tilde 2$)}" there exists a (fixed) rectangle that is associated with one of the duo limit groups,
$Duo^2_1,\ldots,Duo^2_{t_2}$, 
so that for $1 \leq j \leq i \leq L_7$ the value, 
$(x^{i,j}_0,p_i,d_1(i),q_j,d_2(j),a)$,   
 extends to a  value of the
 fixed rectangle.
\endroster

The duo limit group that is associated with the rectangle in part ($\tilde 2$) is an amalgamated
product: $Duo^2=Comp_1(e_1,d_1,p)*_{<e_0,a>} \, Comp_2(e_2,d_2,a)$. Viewing the completions, $Comp_1(e_1,d_1,a)$ and $Comp_2(e_2,d_2,a)$,
as graded limit groups with respect to the parameter subgroups, $<e_0,d_1,a>$ and $<e_0,d_2,a>$ in correspondence,
we have associated graded Makanin-Razborov diagrams with $Comp_1$ and $Comp_2$, and each graded resolution in
these diagrams terminates in either a rigid or a solid limit group. Each value of the variables
$d_1$ and $d_2$, extends to at most $excep^2_{Cov}$ rigid or families of strictly solid values of the terminal
rigid and solid limit groups of the graded resolutions in the graded Makanin-Razborov diagrams of $Comp_1$
and $Comp_2$.

As in the proof of proposition 4.2, we continue by filtering the sequence of values, 
$(x^{i,j}_0,p_i,d_1(i),q_j,d_2(j),a)$,   
 $1 \leq j \leq i \leq L_7$,  
according to the boundedly many
possible extensions of the values $d_2(j)$ to a rigid or a strictly solid (family of) values of a 
terminal rigid or solid limit group of one of the
finitely many graded   resolutions in the graded Makanin-Razborov diagram of $Comp_2(e_2,d_2,a)$ with respect to
the parameter subgroup $<e_0,d_2>$, we are left with a sequence (still denoted),  
$(x^{i,j}_0,p_i,d_1(i),q_j,d_2(j),a)$,   
 $1 \leq j \leq i \leq L_8$,  
 that satisfy properties ($\tilde 1$) and ($\tilde 2$). 

By further filtering the sequence of values, 
$(x^{i,j}_0,p_i,d_1(i),q_j,d_2(j),a)$,   
 $1 \leq j \leq i \leq L_8$,  
according to the boundedly many
possible extensions of the values $d_1(i)$ to a rigid or a strictly solid (family of) values of a terminal 
rigid or solid limit group of one of the
finitely many graded   resolutions in the graded Makanin-Razborov diagram of $Comp_1(e_1,d_1,a)$ with respect to
the parameter subgroup $<e_0,d_1>$, we are left with a sequence (still denoted),  
$(x^{i,j}_0,p_i,d_1(i),q_j,d_2(j),a)$,   
 $1 \leq j \leq i \leq L_9$,  
that satisfy properties ($\tilde 1$) and ($\tilde 2$).

Furthermore, for each  pair of indices,
 $1 \leq j \leq i \leq L_9$,  
the values of the pair, $(d_1(i),d_2(j))$, extends to a value of the duo limit group $Duo^2$, hence,
to values of the two completions, $Comp_1$ and $Comp_2$, from which $Duo^2$ is composed. 
By the filtration that we used, the associated values of the elements $e_1$ (the generators of the completion $Comp_1(e_1,d_1,a)$), that we may
assume to be  values of one of the rigid or solid limit groups in the graded Makanin-Razborov diagram of $Comp_1$ with respect
to the parameter subgroup $<d_1,e_0,a>$,  depend only on the index $i$
and not on the index $j$.  
The associated values of the elements $e_2$ (the generators of the completion $Comp_2(e_2,d_2,a)$), that we may
assume to be  values of one of the rigid or solid limit groups in the graded Makanin-Razborov diagram of $Comp_2$ with respect
to the parameter subgroup $<d_2,e_0,a>$,  depend only on the index $j$
and not on the index $i$.

Finally, $L_9=length^2_{Cov}+2$. For the last sequence,
$(x^{i,j}_0,p_i,d_1(i),q_j,d_2(j),a)$,   
 $1 \leq j \leq i \leq L_9$, we have associated a value $e_1(i)$ with each index $i$, $1 \leq i \leq L_9$,
that is independent of the index $j$, and a value $e_2(j)$ with each value $d_2(j)$, $1 \leq j \leq L_9$,
which is independent of the index $i$. Since   
a pair $(p_i,q_j) \in Cov(p,q)$ if and only if $i<j$, and for $i \geq j$, the restriction of the value
$(d_1(i),d_2(j))$ to the $Wit$ subgroup satisfies the non-degeneracy properties (ii)-(iv), the values
$x_0^{i,j}$, and hence the values $(e_1(i),e_2(j))$, must satisfy the Diophantine condition that demonstrates
that the value $x_0^{i,j}$ is either not rigid or not strictly solid or that it coincides or in the
same strictly solid family of a rigid or a strictly solid value that is one of the $s$ rigid or
strictly solid values that are part of the corresponding specialization of the subgroup $Wit$.
Therefore, starting with the duo limit group $Duo^2$, 
viewed as a graded limit group with respect to the parameter
subgroup $<e_2,d_2,a>$, we obtained a sequence of values:
 $e_2(1),\ldots,e_2(length^2_{Cov}+1)$, of the elements
$e_2$ in the duo limit group $Duo^2$ (the fixed generators of $Comp_2$), 
for which the  sets of values of the variables $e_1$, 
$E1_r$, $1 \leq r \leq length^2{Cov}+1$, for which
these values together with the corresponding values $e_2(1),\ldots,e_2(r)$, $1 \leq r \leq length^2_{Cov}+1$, 
extend to values of $Duo^2$, and the combined
values of $Duo^2$ satisfy the degenerating Diophantine condition, strictly decreases for
$1 \leq r \leq length^2_{Cov}+1$. This contradicts the choice of $length^2_{Cov}$ to be a global bound on 
the length of such
strictly decreasing sequences of values of the variables $e_2$ for all the rectangles in all the duo limit groups,
$Duo^2_1,\ldots,Duo^2_{t_2}$.

\line{\hss$\qed$}

Proposition 4.4 proves the stability of the sets $Cov(p,q)$ and theorem 4.3 follows.

\line{\hss$\qed$}

\vglue 1.5pc
\centerline{\bf{\S5. Stability}}
\medskip
In the previous section we have shown that  
the sets $NR_s(p,q)$, $NS_s(p,q)$, and $Cov(p,q)$, that indicate those values of the 
parameter set $<p,q>$, for which a rigid limit group $Rgd(x,p,q,a)$ admits
at least $s$ rigid values, a solid limit group admits at least $s$ 
strictly solid families of values(theorem 4.1), and a given 
finite set
of (graded) closures forms a covering closure of a given graded resolution
(theorem 4.3),
are stable. 
In this section we combine these theorems, with the arguments that were
used in proving theorem 1.9, to prove that a general definable set over a free group is stable.

\vglue 1pc
\proclaim{Theorem 5.1} The elementary theory of a non-abelian free
group is stable.
\endproclaim

\nfp The argument that we use is a rather straightforward modification of the argument that was used in
the minimal (graded) rank case (theorem 1.9). Let $L(p,q)$ be a definable set over a non-abelian free group $F_k$.
As in the proof of theorem 1.9, we need to use the objects and terminology that is used in the sieve procedure
that finally leads to quantifier elimination, and is presented in [Se5] and [Se6]. 
The exact definitions of these objects is long and involved, and
we refer the reader to section 1 in [Se5] for a detailed presentation of them. 

The sieve procedure, that is used to prove quantifier elimination, is much more difficult in the general case
[Se6], in comparison with the minimal rank case (section 1 in [Se5]). Still, the overall strategy for quantifier
elimination, and the output of the sieve procedure in the general case and in the minimal rank case are similar.

Recall that as for minimal rank definable sets,  
with a (general) definable set, $L(p,q)$, the sieve procedure associates a finite collection of graded PS 
resolutions,
and with each such graded PS resolution it associates a finite collection of graded
closures of these resolutions that contains Non-Rigid, Non-Solid, Left, Root, Extra PS, and 
collapse extra PS resolutions
(see definitions  1.25-1.30 of [Se5] for the exact definitions of these resolutions). 

Let $PSRes_i$, $i=1,\ldots,r$, be  the finitely many PS resolutions that are associated with the given definable set
$L(p,q)$. For each index $i$, $i=1,\ldots,r$, let $Rgd_i(x,p,q,a)$ ($Sld_i(x,p,q,a)$) be the terminal rigid (solid) limit group of $PSRes_i$.
With the PS resolution $PSRes_i$ and its terminal rigid or solid limit group $Rgd_i$ or $Sld_i$, we associate the definable
set, $NR^i_1(p,q)$ or $NS^i_1(p,q)$, that defines those values of the defining parameters $p,q$ for which $Rgd_i$ ($Sld_i$)
extends to a rigid or a strictly solid value of $Rgd_i$ or $Sld_i$. By theorem 4.1 the sets $NR^i_1$ and $NS^i_1$ are stable.

With each of the PS resolutions, $PSRes_i$, 
the sieve procedure  associates a finite collection of graded closures of
it that contains Non-Rigid, Non-Solid, Left, Root, Extra PS, and 
collapse extra PS resolutions. With the graded resolution $PSRes_i$, and its given set of closures, we associate
a definable set $Cov_i(p,q)$, that defines those values of the defining parameters $p,q$ for which the associated
fibers of $PSRes_i$ that are associated with the value $p,q$ are covered by the fibers that are associated with the
given finite set of closures of it and with the value of $p,q$.  By theorem 4.3 $Cov_i(p,q)$ is stable.

By the sieve procedure (cf. the proof of theorem 1.9),  the definable set $L(p,q)$ is the finite union:
$$cup_{i=1}^r \, NR^i_1(p,q) \, (NS^i_1(p,q)) \, \setminus \, Cov_i(p,q)$$
In particular, $L(p,q)$ is a Boolean combination of the sets $NR^i_1$ ($NS^i_1$) and $Cov_i$. Since by theorems 4.1 and 4.3, 
the sets, $NR^i_1$ ($NS^i_1$) and $Cov_i(p,q)$, are stable, and the collection of stable sets is closed under Boolean operations,
so is their Boolean
combination $L(p,q)$, hence, the theory of a free group is stable.

\line{\hss$\qed$}

According to [Se8], a definable set over a non-elementary, torsion-free hyperbolic group can be analyzed
using the same sieve procedure as the one constructed over a free group  (see section 6 in [Se8]). 
As a corollary, like over a free group, every definable set over a torsion-free hyperbolic group is a Boolean combination
of sets of the form: 
$NR^i_1$ ($NS^i_1$) and $Cov_i(p,q)$, where these sets are defined precisely as they are defined over free groups in
theorems 4.1 and 4.3.

As the sieve procedure generalizes to torsion-free hyperbolic groups, the argument that proves that Diophantine sets are 
equational over free groups (theorem 2.1) generalizes to every torsion-free hyperbolic group. 
The definitions of  rigid and solid limit groups generalize to torsion-free hyperbolic groups, and the global boundedness of the number
of rigid and strictly solid families of values of a rigid or a solid limit group for any given value of the defining
parameters that holds over a free group holds over every torsion-free hyperbolic group (see section 3 in [Se8]).

Hence, one can define configuration limit groups that are associated with a rigid or a solid limit group over a torsion-free
hyperbolic group, and the construction of  Duo limit groups that is presented in section 3 generalizes to torsion-free hyperbolic 
groups as well. Finally, the 
 arguments that were used in proving theorems
4.1 and 4.3 over a free group, generalize to torsion-free group, so the sets $NR_S$, $NS_S$, and $Cov$, that are proved to be stable over  
a free group are stable over any torsion-free hyperbolic group. As any definable set over a torsion-free hyperbolic group is a Boolean
combination of sets of the form $NR_s$, $NS_s$, and $Cov$,  and the collection of stable sets is closed under Boolean operations,
every definable set over a torsion-free hyperbolic group is stable.

\vglue 1pc
\proclaim{Theorem 5.2} The elementary theory of a non-elementary (torsion-free) hyperbolic 
group is stable.
\endproclaim


\vskip 2em

\end